\documentclass[10pt]{amsart}

\usepackage{amsmath,amsthm,amsfonts,amssymb,url,constants}
\usepackage[francais]{babel}
\newconstantfamily{c}{symbol=c}

\usepackage[T1]{fontenc}
\usepackage{xypic}

\newcommand{\GQ}{{\Gal(\bar \Q/\Q)}}
\newcommand{\p}{\mathfrak{p}}

\newcommand{\Q}{\mathbb{Q}}
\newcommand{\Z}{\mathbb{Z}}

\newcommand{\R}{\mathbb{R}}

\newcommand{\End}{\mathrm{End}}
\newcommand{\Spec}{\mathrm{Spec\,}}
\newcommand{\spec}{\mathrm{Spec\,}}
\newcommand{\Frob}{{\rm{Frob}}}

\newcommand{\GL}{\mathrm{GL}}
\newcommand{\GSP}{\mathrm{GSP}}
\newcommand{\Sl}{\mathrm{SL}}
\newcommand{\U}{\mathrm{U}}

\newcommand{\Gl}{\mathrm{GL}}

\newcommand{\Gal}{\mathrm{Gal}}

\newcommand{\F}{\mathbb{F}}

\newcommand{\tr}{\mathrm{tr}\,}

\newcommand{\Td}{{\widehat{T}}}
\newcommand{\bigo}[1]{{O\left(#1\right)}}

\newtheorem{question}{Question}

\newtheorem{conjecture}{Conjecture}

\newtheorem{theoremerep}{Théorème}
\newtheorem{theoreme}{Th\'eor\`eme}
\newtheorem{prop}{Proposition}
\newtheorem{cor}{Corollaire}
\newtheorem{lemme}{Lemme}

\theoremstyle{definition}
\newtheorem{definition}{Définition}
\newtheorem{remarque}{Remarque}
\newtheorem{remarques}{Remarques}
\newtheorem{exemple}{Exemple}

\renewenvironment{pf}
        {\noindent {\it Démonstration --- \ }}
        {\hfill\nobreak $\Box$ \par\bigbreak}

\newcommand{\Id}{\text{Id}}
\newcommand{\un}{{\bf 1}}

\newcommand{\Li}{\text{Li}\,}

\newcommand{\NN}[1]{|\!| #1 |\!|}
\newcommand{\NA}[1]{|\!| #1 |\!|_{A(G)}}
\newcommand{\NAA}[2]{|\!| #1 |\!|_{A( #2 )}}

\newcommand{\herm}[1]{\langle #1 \rangle}
\newcommand{\Gd}{{\widehat{G}}}
\newcommand{\Hd}{{\widehat{H}}}
\newcommand{\Ind}{{\text{Ind}}}

\newcommand{\supp}{{\text{Supp}}}

\newcommand{\reg}{{\text{reg}}}
\newcommand{\nreg}{{\text{nreg}}}
\newcommand{\nr}{{\text{nr}}}
\newcommand{\diag}{{\text{diag}}}
\newcommand{\Cb}{{\mathbb{C}}}
\newcommand{\cc}{\C[c]}
\newcommand{\ccr}{\Cr}
\newcommand{\ccl}{\Cl[c]}

\renewcommand{\Lc}{{\mathcal{L}}}
\newcommand{\Bc}{{\mathcal{B}}}
\newcommand{\oui}{{\underline{\text{oui}}}}
\newcommand{\non}{{\underline{\text{non}}}}

\newcommand{\ta}{{\tr=a}}
\newcommand{\tz}{{\tr=0}}
\newcommand{\dr}{{\rm{dr}}}
\date{}

\title{Théorème de Chebotarev et complexité de Littlewood}

\begin{document}
\baselineskip 15pt

\address{Jo\"el Bella\"iche\\Brandeis University\\
415 South Street\\Waltham, MA 02454-9110\\U.S.A}
\email{jbellaic@brandeis.edu}
\author{Jo\"el Bella\"iche}
\maketitle

\tableofcontents

\par \bigskip
\par \bigskip

\section{Introduction}

\subsection{Objectifs}

Sous l'hypothèse de Riemann généralisée pour les fonctions $L$ d'Artin, que nous noterons comme d'habitude (GRH), et sous la conjecture d'Artin, 
le théorème de Chebotarev admet une preuve simple et naturelle et une forme effective élégante et forte (cf. \cite[page 148]{ik})\footnote{La preuve de cette forme du théorème de Chebotarev
apparaît pour la première fois dans l'article \cite{effective} de Murty, Murty, et Saradha, et est reprise dans le livre \cite{eff} de Murty et Murty mais dans les deux cas le résultat n'est énoncé que sous une forme affaiblie. La
forme forte du théorème est énoncée et prouvée par Iwaniec et Kowalski dans leur livre \cite{ik}. Il y est également affirmé qu'on peut prouver le même résultat sans supposer la conjecture d'Artin, mais cette affirmation est erronée, comme E. Kowalski me l'a confirmé.}, mais qui n'a semble-t-il jamais encore été utilisée.
Cette forme effective fait apparaître un invariant d'un sous-ensemble invariant $D$ d'un groupe fini $G$ que nous appelons sa {\it complexité de Littlewood} $\lambda_G(D)$. Le but de cet article
est de commencer une étude détaillée de  cet invariant, de prouver un certain nombre de corollaires de ce théorème de Chebotarev et de donner des applications arithmétiques "concrètes" de ces résultats concernant par exemple le plus petit nombre premier qui engendre le groupe cyclique 
$(\Z/\ell \Z)^\ast$; le plus petit nombre premier modulo lequel un polynôme $P\in\Z[X]$ irréductible fixé a une racine (ou bien a au moins deux racines; ou n'en a aucune, ou encore reste irréductible); la conjecture de Lang-Trotter pour les courbes elliptiques et ses généralisations, notamment aux courbes  et aux variétés abéliennes; la forme effective du critère d'isogénie de Faltings; la conjecture d'uniformité de Serre. 

Ces résultats et applications se divisent en trois familles, qui sont les suivantes: premièrement le théorème de Chebotarev proprement dit (i.e. concernant la densité des nombres premiers dont le  Frobenius  dans une extension finie de $\Q$ de groupe de Galois $G$ est dans un certain ensemble invariant par conjugaison $D \subset G$) et les propriétés de la complexité de Littlewood $\lambda_G(D)$; deuxièmement les théorèmes permettant de majorer le plus petit nombre premier dont le Frobenius dans $G$ appartient à $D$; troisièmement,  les généralisations du théorème de Chebotarev au cas d'une extension infinie, ou d'un système infini d'extensions -- ce sont ces généralisations qui sont utiles pour la conjecture de Lang-Trotter par exemple. Dans le reste de cette introduction, nous discutons en détail ces trois familles de résultats  et les idées qui les sous-tendent. 

\subsection{Le théorème de Chebotarev et la complexité de Littlewood}

Soit $L$ une extension finie galoisienne\footnote{Pour ne pas alourdir les notations, nous ne considérons dans cet article que des extensions de corps de nombres  
dont le corps de base est $\Q$. Il est en principe toujours possible de se ramener à ce cas.} de $\Q$ de groupe de Galois $G$. Soit $M$ le produit des nombres premiers ramifiés dans $L$. Pour $p$ un nombre premier ne divisant pas $M$, on note $\Frob_p$, (ou $\Frob_{p,G}$ quand il y a une ambiguité sur l'extension considérée) la classe de conjugaison de l'élément de Frobenius associé à $p$ dans $G$. Pour $f$ une fonction centrale (i.e. invariante par conjugaison) sur $G$ à valeurs complexes, posons $$\pi(f,x) = \sum_{p < x, p \nmid M} f(\Frob_p).$$
Associons à $f$ deux invariants: le premier,
$$\mu(f) = \mu_G(f)=  \frac{1}{|G|} \sum_{g \in G} f(g)$$ est simplement la valeur moyenne de $f$. Le second est
\begin{eqnarray} \label{defNL} \lambda(f) = \lambda_G(f) =  \sum_{\pi} | \hat{f}(\pi) | \dim \pi,\end{eqnarray}
où $\pi$ parcourt l'ensemble $\Gd$ des classes d'équivalences de représentations complexes irréductibles de $G$, et 
$\hat{f}: \Gd \rightarrow \Cb$ est la {\it transformée de Fourier} de $f$, définie par $\hat{f}(\pi)=\frac{1}{|G|} \sum_g \tr \overline{f(g)} \pi(g)$.
 Nous appellerons  $\lambda(f)$ la {\it norme de Littlewood} de $f$. 
 
 L'ensemble des représentations irréductibles $\pi$ de $G$ tels
que $\hat{f}(\pi) \neq 0$ est le {\it support spectral} de $f$. Comme d'habitude $\Li(x)$ est le logarithme intégral: $\Li(x)=\int_2^x \frac{1}{\log t}dt$.

\begin{theoreme}[Chebotarev effectif] \label{chef} Supposons vraie (GRH) et la conjecture d'Artin pour les
 fonctions $L$ d'Artin associées aux représentations irréductibles de $\Gal(L/\Q)$ qui appartiennent au support spectral de $f$. Il existe une constante absolue $\ccl{ct}>0$ telle que pour $x \geq 3$, on ait:
\begin{eqnarray} \label{chefe} \left| \pi(f,x) - \mu(f) \Li(x)  \right| < \ccr{ct} x^{1/2} \lambda(f) (\log x+\log M+\log |G|) \end{eqnarray}
\end{theoreme}

 Comme cet énoncé diffère, quoique très légèrement, de la forme donnée dans \cite{ik}, nous expliquons comment le déduire en \S\ref{densitefrob}.

Le cas le plus important est  celui où $f$ est la fonction indicatrice $\un_D$ d'un sous-ensemble $D$ de $G$ invariant par conjugaison. 
On note alors $\pi(D,x)$ pour $\pi(f,x)$, et ce nombre est le nombre d'éléments plus petits que $x$ de l'ensemble $\tilde D$ des nombre premiers $p$ ne divisant pas $M$ tels que $\Frob_p \in D$. Un ensemble de la forme $\tilde D$ sera appelé un ensemble Frobénien. Nous noterons
 $\lambda(D)$ pour $\lambda(\un_D)$, et nous appellerons ce nombre réel positif la {\it complexité de Littlewood} de $D$. De la même façon, 
 on appellera {\it support spectral} de $D$ le support spectral de $\un_D$. 
 On a évidemment $\mu(D)=\frac{|D|}{|G|}$ si bien que le théorème de Chebotarev effectif ci-dessus prend la forme: 
 \begin{theoreme} \label{thmchefD} Supposons vraie (GRH) et la conjecture d'Artin pour les
 fonctions $L$ d'Artin associées aux représentations irréductibles de $\Gal(L/\Q)$ qui appartiennent au support spectral de $D$.
  Pour $x \geq 3$, 
 \begin{eqnarray} \label{chefD} \left | \pi(D,x) - \frac{|D|}{|G|} \Li(x) \right| < \Cr{ct} x^{1/2} \lambda(D) (\log x +\log M +\log |G| ).\end{eqnarray}
\end{theoreme}

\par \medskip
Comme nous l'avons dit, cette forme précise du théorème de Chebotarev effectif n'a jamais été utilisée jusqu'ici
(mentionnons tout de même que dans le livre \cite{kowalski} de Kowalski, une variante pour $L/\Q$ remplacée par une extension de corps de fonctions l'est). Dans la littérature, on utilise à la place des formes affaiblies, à savoir celles qu'on obtient en remplaçant $\lambda(D)$ par $|D|$, ce qui revient à appliquer la majoration triviale $\lambda(D) \leq |D|$ -- c'est le théorème employé dans \cite{serre}, qui est une légère amélioration du théorème de Lagarias and Odlyzko, \cite{lo}, que nous appellerons ici {\it la version de Lagarias-Odlyzko-Serre}, ou bien $\lambda(D)$ par $\sqrt{|D|}$ (ce qui revient à appliquer la majoration dite "de Cauchy-Schwarz" $\lambda(D) \leq \sqrt{|D|}$) -- c'est le théorème employé dans \cite{effective} et \cite{eff}, que nous appellerons {\it la version de Murty-Murty-Saradha}.

Pour obtenir une application du théorème de Chebotarev effectif ci-dessus qui ne soit pas directement conséquence de la version de Murty-Murty-Saradha, 
il faut donc, dans des cas particuliers intéressants, prouver une meilleure majoration de $\lambda(D)$ que la borne $\sqrt{|D|}$. Dans la plupart des applications, on travaille avec une famille $(L_\nu)_{\nu \geq 1}$ de corps de nombres galoisiens sur $\Q$ de groupes de Galois $G_\nu$, dont l'ordre $|G_\nu|$ tend vers l'infini, et avec une famille
de sous-ensembles $D_\nu$ de $G_\nu$ invariants par conjugaison, et il importe d'obtenir une majoration de $\lambda(D_\nu)$ qui soit {\it asymptotiquement} meilleure que la majoration $\lambda(D_\nu) \leq |D_\nu|^{1/2}$, par exemple une majoration de la forme $\lambda(D_\nu) \leq  |D_\nu|^{\alpha}$ avec $\alpha<1/2$. Afin d'aider à l'obtention de telles majorations, nous développons une "boîte à outils" permettant de manipuler cet invariant $\lambda(D)$ (et plus généralement $\lambda(f)$).  Certains de ces outils sont des transpositions de
techniques utilisées par Serre dans \cite{serre} ou par Murty-Murty-Saradha dans \cite{effective} pour obtenir de meilleurs résultats que ce qu'ils obtiendraient 
par une application directe de leurs versions respectives du théorème de Chebotarev. Cependant, alors que pour Serre et Murty-Murty-Saradha ces techniques
supposaient de changer, souvent plusieurs fois, l'extension de corps de nombres considérée et de traduire, de manière parfois un peu délicate,  les résultats pour une extension en termes d'une autre, elles deviennent dans le langage de la complexité de Littlewood des lemmes élémentaires sur les représentations des groupes finis. Ainsi, la complexité de Littlewood permet non seulement d'obtenir dans certains cas de meilleurs résultats, mais aussi de retrouver plus simplement
des résultats de \cite{serre} ou \cite{effective}.

Quand la table des caractères d'un groupe fini $G$ est connue, on peut en théorie calculer la complexité de Littlewood $\lambda_G(D)$ de n'importe quel
ensemble invariant par conjugaison $D$, et la comparer avec la majoration de Cauchy-Schwarz $\sqrt{|D|}$. La table des caractères est connue dans plusieurs cas particulièrement importants dans les applications: celui où $G$ est abélien, 
celui où $G$ est le groupe symétrique $S_n$, pour lequel on dispose d'une description combinatoire de la table des caractères, ou encore celui où $G$ est un groupe fini de type de Lie, où la théorie de Deligne-Lusztig et les travaux subséquents de Lusztig donnent une description complète des caractères. "En théorie" seulement, car le calcul se révèle souvent extrêmement difficile, même dans le cas abélien. Nous donnons dans cet article un certain nombre de tels calculs dans des cas où nous avons pu les mener à bien: intervalles dans $\Z/n\Z$ (la borne de Cauchy-Schwartz n'est pas du tout optimale), intervalles dans 
$(\Z/n\Z)^\ast$ (la borne de Cauchy-Schwarz est essentiellement optimale), générateurs de $\Z/n\Z$ (la borne de Cauchy-Schwartz n'est pas optimale),
dérangements dans $S_n$ (la borne de Cauchy-Schwartz n'est pas optimale), matrices de traces fixées dans $\Gl_n(\F_\ell)$ (la borne de Cauchy-Schwarz est essentiellement optimale) ou dans $\GSP_4(\F_\ell)$ (elle ne l'est pas). Comme cette liste le suggère, il reste beaucoup à faire, en particulier une étude systématique en utilisant la théorie de Deligne-Lusztig de la complexité de Littlewood de l'ensemble des matrices de trace donnée dans un sous-groupe de type de Lie de $\Gl_n(\F_\ell)$.

Notons par ailleurs qu'il n'est pas difficile de voir (cf. \S \ref{egaliteNA}) que $\lambda(f)$ est égale à la {\it norme d'algèbre $\NA{f}$ de $f$}, une notion introduite par Pierre Eymard en 1964 
pour une fonction $f$ sur un groupe localement compact (\cite{eymard}), et bien plus ancienne dans le cas d'un groupe abélien localement compact $G$, puisque dans ce cas $\NA{f}$ est la norme $L^1$ 
de la transformée de Fourier $\hat{f}$ de $G$. Que cette notion, fondamentale en analyse harmonique, soit importante an arithmétique est évident 
depuis les travaux de Hardy et Littlewood. Rappelons à ce sujet la conjecture de Littlewood (d'où les noms "norme de Littlewood", et "complexité de Littlewood") selon laquelle, pour $D$ un sous-ensemble fini de $G=\Z$, on a $\NA{\un_D} \gg \log |D|$.
La conjecture de Littlewood (maintenant un théorème) a fait l'objet de nombreuses recherches et généralisations, notamment aux cas des groupes abéliens localement compacts. Tout récemment, sa généralisation au cas d'un groupe fini non nécessairement abélien (c'est-à-dire au cas que nous considérons ici)
a fait l'objet d'un article de Tom Sanders \cite{sanders}. 
Puisque pour la conjecture de Littlewood  et ses variantes, on cherche à minorer $\NA{\un_D}$, tandis que nous cherchons au contraire ici à le majorer, 
nos résultats sont essentiellement disjoints, et en fait, complémentaires, de ceux de \cite{sanders}. Notons qu'il ne semble pas que les applications au théorème de Chebotarev fassent partie des motivations de Tom Sanders.

\subsection{Le plus petit nombre premier $p$ tel que $\Frob_p \in D$ et l'invariant $\varphi_G(D)$.}

On reprend les notations du paragraphe précédent: $L/\Q$ une extension finie de groupe de Galois $G$, ramifiée exactement aux diviseurs premiers de $M$, et $D \subset G$ un sous-ensemble stable par conjugaison. A ces données on attache l'ensemble frobénien $\tilde D$ des nombres premiers $p$ ne divisant pas $M$
tels que $\Frob_{p,G} \in D$. On considère le problème consistant à {\bf majorer le plus petit nombre premier $p$ de l'ensemble frobénien $\tilde D$.} 

Notons tout d'abord qu'on obtient une telle majoration comme un corollaire direct du théorème de Chebotarev effectif (\ref{chefD}): dès que $x$ est assez grand pour que le terme principal $\frac{|D|}{|G|} \Li x$ soit supérieur à la valeur absolue du terme d'erreur, par exemple dès que $x / \log^2 x > \cc \left( \frac{|D|\lambda(D)}{|G|} (\log M +\log |G| ) \right)^2$,  on sait évidemment qu'il existe un $p \leq x$ appartenant à $\tilde D$. La même méthode permet
bien entendu d'obtenir une majoration du plus petit $p$ de $\tilde D$ en partant de la version de Lagarias-Odlyzko-Serre, ou de celle de Murty-Murty-Saradha 
du théorème, et la majoration que nous venons d'obtenir est meilleure exactement dans les mêmes circonstances, discutées plus hauts, que le 
théorème~\ref{thmchefD} est meilleur que les autres versions. 

En fait, dans la version de Lagarias-Odlyzko de la majoration du plus petit nombre premier de $\tilde D$, (cf. \cite[Théorème 5]{serre} et la remarque qui le suit), comme dans celle de Murty-Murty-Saradha, et dans celle donnée ci-dessus une astuce permet d'améliorer le facteur logarithmique de la borne obtenue. 
C'est un détail, que nous ne discuterons pas plus avant dans cette introduction. 

Il y a en revanche une possibilité d'amélioration bien plus importante, basée sur l'idée très simple suivante:  au lieu d'utiliser l'estimation (\ref{chefD}) de $\pi(D,x)$, on peut utiliser l'estimation (\ref{chefe}) de $\pi(f,x)$ pour une fonction $f: G \rightarrow \R$ centrale adéquate. Par adéquate, nous voulons dire  que la moyenne $\mu_G(f)$ est strictement positive et que $f$ ne prenne aucune valeur strictement positive en dehors de $D$. Ainsi, dès que $x$ est suffisamment grand pour que le terme principal $\mu_G(f) \Li(x)$ de (\ref{chefe}) soit supérieur à la valeur absolue du terme d'erreur de (\ref{chefe}), on sait qu'il existe un $p \leq x$ tel que $f(\Frob_p)>0$, donc tel que $p \in \tilde D$. On peut alors rechercher, parmi toutes les fonctions $f$ {\it adéquates} celle qui donne de cette manière la
meilleure  borne sur $p$. 
Cette idée conduit à définir l'invariant suivant.

\begin{definition} \label{defphig} Soit $D$ un sous-ensemble non-vide  stable par conjugaison de $G$. Nous noterons
$$\varphi_G(D) = \inf_f \frac{\lambda_G(f)}{\mu_G(f)},$$ quand $f$ parcourt l'ensemble des fonctions sur $G$ à valeurs réelles,
qui satisfont les deux conditions suivantes
\begin{itemize} \item[(a)] si $f(g)>0$, alors $g \in D$.
\item[(b)] $\mu_G(f) >  0$.
\end{itemize} 
\end{definition}

Cet invariant $\varphi_G(f)$ s'avère beaucoup plus difficile à étudier que $\lambda_G(f)$. On a évidemment 
\begin{eqnarray} \label{majvarpphi} \varphi_G(f) \leq \lambda_G(D)  |G| / |D| \end{eqnarray} car $f=\un_D$ satisfait (a) et (b), mais nous verrons des cas intéressants où cette inégalité est loin d'être une égalité. Nous montrerons, suivant l'idée esquissée ci-dessus:
\begin{theoreme} \label{thmpluspetit} Pour toute extension finie galoisienne $L/\Q$, $G=\Gal(L/\Q)$, 
$D \subset G$ non-vide et invariant par conjugaison, $M=\prod_{p \text{ ramifié dans } L} p$ comme ci-dessus, supposant GRH et la conjecture d'Artin vraie pour les fonctions $L$ d'Artin attachées aux représentations de $G$, le plus petit nombre premier
$p$ de $\tilde D$ vérifie:
$$p < \ccl{cpp} \varphi_G(D)^2 \log^2 M,$$
où $\ccr{cpp}$ est une constante absolue.
\end{theoreme} 

Bien entendu, on aurait pu utiliser la même idée (utiliser $\pi(f,x)$ au lieu de $\pi(D,x)$ pour une fonction $f$ bien choisie) pour obtenir une majoration du plus petit nombre premier de $\tilde D$
en utilisant, au lieu de (\ref{chefe}) la version de Murty-Murty-Saradha (resp. de Lagarias-Odlyzko-Serre) de cette estimation, i.e. celle qu'on obtient en remplaçant $\lambda_G(f)$ par $\NN{f}_2$, resp. $\NN{f}_1$. Cependant, cette plus grande généralité ne donnerait pas un résultat meilleur, la fonction $f$ {\it la mieux choisie} étant dans ce cas toujours $\un_D$ (cf. remarque~\ref{remarqueppm} ci-dessous). Autrement dit, c'est la souplesse que donne l'emploi de la norme de Littlewood qui permet de choisir un $f$  meilleur que $\un_D$. Insistons sur le fait que la borne du théorème~\ref{thmpluspetit} peut être meilleure que celle obtenue par Murty-Murty-Saradha (et a fortiori, celle obtenue par Lagarias-Odlyzko) pour deux raisons: parce que $\lambda_G(D)$ peut être strictement plus petit que $\sqrt{|D|}$ d'une part, et parce que 
$\varphi_G(f)$ peut être plus petit que  $\lambda_G(D)  |G| / |D|$, autrement dit que la borne inférieure définissant $\varphi_G(D)$ n'est pas atteinte par $\un_D$. Parfois ces deux raisons se cumulent, comme dans les calculs menant à l'application suivante.

\begin{theoreme}
Soit $P$ un polynôme unitaire à coefficients entiers, de degré $n \geq 1$. Soit $M$ le produit des nombres premiers divisant le discriminant de $P$. Supposons GRH et la conjecture d'Artin vraies pour les fonctions $L$ du corps de décomposition de $P$. Alors,
il existe un nombre premier $p  < 4 \ccr{cpp}  (n^2+n)^2 (\log M + n \log n)^2 $ ne divisant pas $M$ tel que le polynôme $P(X) \pmod{p}$ n'admette aucune racine dans $\F_p$.
\end{theoreme}
La borne qu'on obtiendrait avec la version de Murty-Murty-Saradha ou de Laragias-Odlyzko serait super-exponentielle en $n$, au lieu de polynomiale. Pour d'autres applications du même genre, voir \S\ref{redpolint}

\subsection{Cas des extensions infinies, ou des systèmes infinis d'extensions}

Souvent on est amené à considérer une extension $L/\Q$ algébrique infinie de groupe de Galois $G$, non ramifiée en dehors des nombres premiers divisant un entier $M \geq 1$ fixé. Si $D$ est une partie fermée, le problème de déterminer un équivalent de $\pi(D,x)$, le nombre de $p<x$, $p \nmid M$ tels que $\Frob_{p,G} \in D$ 
s'avère beaucoup plus délicat que dans le cas d'une extension finie (sauf bien sûr si $D$ est à la fois fermée et ouverte, où l'on se ramène facilement au cas
d'une extension finie). Comme Serre l'a montré dans \cite{serre}, le théorème de Chebotarev effectif (appliqué à des sous-extensions finies $L_\nu$ de $\Q$, de groupes de Galois $G_\nu$, avec $G = \projlim G_\nu$) permet néanmoins de montrer des majorations non triviales
de $\pi(D,x)$. Une variante de la situation ci-dessus est celle où l'on part dès le début d'un système d'extensions galoisiennes $L_\nu/\Q$ (finies où non) de
groupes de Galois $G_\nu$, et qu'étant données des $D_\nu \subset G_\nu$ on s'intéresse au nombre $\pi(D,x)$ de nombres premiers $p<x$ tels que $\Frob_{p,G_\nu} \in D_\nu$ pour tout $\nu$.

Dans cet article, nous développons la méthode de Serre dans deux directions essentiellement indépendantes. Premièrement, nous 
montrons comment la meilleure version du théorème de Chebotarev que nous utilisons permet de donner des résultats plus forts si l'on connaît
des majorations des complexités de Littlewood qui interviennent qui sont meilleures que celles données par la borne de Cauchy-Schwarz. 
Deuxièmement, nous utilisons la  méthode du grand crible, combinée avec certaines propriétés élémentaires de la norme de Littlewood, pour obtenir des bornes qui sont meilleures même quand on n'a pas de meilleure majoration des complexités qui interviennent que celles de Cauchy-Schwarz. 
Nous renvoyons à \S\ref{sectdenszero} pour les énoncés, qui sont un peu techniques, et pour les méthodes utilisées, et à \S\ref{syscompLT}
pour des applications à la conjecture de Lang-Trotter pour des systèmes compatibles, en particulier pour des variétés abéliennes.

Contentons-nous ici de donner les résultats que nous obtenons pour une courbe elliptique $E/\Q$, non CM, ayant bonne réduction
hors $M \geq 1$. On note $a_p = |E(\F_p)|-1-p$ pour $p \nmid  M$ et pour $a$ un entier;
on pose $\pi(a,x) = |\{p < x, p \nmid M, a_p=a\}|$. 

Quand $a=0$, $\pi(0,x)$ est le nombre de nombres premiers $p<x$ en lesquels $E$ a bonne réduction supersingulière. Serre a montré sous (GRH) que $\pi(0,x)=\bigo{x^{3/4}}$. Murty-Murty-Saradha, avec leur version sous Artin et (GRH) ont reconsidéré le problème mais n'ont pu obtenir qu'une autre preuve du même résultat (sous (GRH) seule, car une astuce permet dans ce cas de se ramener à une situation où Artin est connue), résultat que Elkies a réussi à obtenir inconditionnellement. Nous obtenons dans cet article, sous (GRH) et Artin: $$\pi(0,x) = \bigo{ x^{2/3+\epsilon}}\text{ pour tout $\epsilon >0$}.$$
Quand $a \neq 0$, Serre avait obtenu sous (GRH), $\pi(a,x)=\bigo{x^{7/8} (\log x)^{-1/2}}$, résultat que Murty-Murty-Saradha ont amélioré en
$\pi(a,x) = \bigo{x^{4/5} (\log x)^{-1/5}}$. Reconsidérant le problème en lui appliquant une méthode du crible, sous (GRH) et Artin, 
Zywina a récemment obtenu une nouvelle preuve de la même formule. Nous prouvons ici:
$$\pi(a,x)=\bigo{x^{5/7+\epsilon}} \text{ pour tout $\epsilon>0$.}$$

\par \bigskip 
{\bf Notations: } Dans tout l'article, les nombres $c_1,c_2,c_3,$ etc. sont des constantes positives. Ce sont des constantes absolues, sauf mentions explicites
du contraire qui resteront assez rares et préciseront alors de quoi elles dépendent. On emploie la notation de Landau tel que celui-ci et ensuite Bourbaki,
l'employait: si $f(x)$ et $g(x) >0$ sont deux fonctions d'une variable $x$ définies sur une partie $P$ non majorée de $\R$, $f(x)=O(g(x))$ s'il existe une constante
$A$ et une constante $C>0$ telles que pour tout $x \in P$, $x \geq A \Longrightarrow | f(x) | < C g(x)$. Quand $f(x)$ et $g(x)$ dépendent de paramètres
autres que $x$, les constantes $A$ et $C$ peuvent en dépendre aussi.   
\par \bigskip
{\bf Remerciements: } Je tiens à remercier en tout premier lieu la communauté du forum mathoverflow. J'y ai posé durant la rédaction de cet article 
de nombreuses questions, dont les réponses (ainsi parfois que les réponses à des questions posées par d'autres) m'ont donné des références, des idées, parfois mêmes des preuves, que j'ai utilisées (dans ce dernier cas, les preuves sont attribuées nommément ci-dessous à leurs auteurs sur mathoverflow).
Je remercie également l'université de Yale, qui m'a accueilli pendant la rédaction de cet article, et où j'ai baigné dans une ambiance beaucoup plus "théorie analytique des nombres" que celle, plus orientée "théorie algébrique des nombres", dont je bénéficie habituellement à Boston. Ce travail a été  
motivé en premier lieu par  (et sera appliqué ultérieurement à) des questions concernant les formes modulaires modulo $2$ et la fonction de partition: c'est grâce à Jean-Louis Nicolas et Jean-Pierre Serre 
que je me suis intéressé au sujet. Je les en remercie ici. Je remercie également  Emmanuel Kowalski pour d'intéressantes conversations, et pour avoir attiré mon attention sur son livre \cite{kowalski} ainsi que tous ceux qui ont manifesté un intérêt pour ce travail durant sa réalisation.

\section{La norme de Littlewood $\lambda(f)$ et la complexité de Littlewood $\lambda(D)$}

\subsection{Égalité avec la norme d'algèbre $\NA{f}$}
\label{egaliteNA}

\subsubsection{Rappel: la norme d'algèbre pour un groupe localement compact, d'après Eymard (\cite{eymard})}

Dans ce~\S, nous ne supposons pas que $G$ est fini, mais que $G$ est un groupe localement compact, muni d'une mesure de Haar $\mu_G$ invariante à gauche,
ce qui permet de parler des espaces $L^1(G)$ et $L^2(G)$, munis de leur normes naturelles que nous noterons $\NN{f}_1$ et $\NN{f}_2$ et, en ce qui concerne  $L^2(G)$, du produit hermitien noté $\herm{\, , \, }$, ainsi que de la représentation régulière gauche $L$ de $G$ sur $L^2(G)$: $L(x)(f)(y)=f(xy)$.

Pour toute fonction $f \in L^1(G)$, et toute représentation unitaire $r$ de $G$ sur un espace de Hilbert $V$, on définit l'opérateur continu $r(f)$ de $V$
par la formule usuelle: $$r(f)(v) = \int_G f(x) r(x)(v) d\mu_G(x).$$
On notera $f^\ast$ la fonction définie par $f^\ast(x)=f(x^{-1})$ si bien que $r(f^\ast)$ est l'adjoint 
$r(f)^\ast$ de $r(f)$. En particulier, on dispose pour $f \in L^1(G)$ d'un opérateur continu $L(f)$ sur $L^2(G)$, qui n'est autre que la convolution à gauche par $f$. On notera $\NN{L(f)}$ la norme d'opérateur de $L(f)$. 

\begin{definition}  Soit $f \in L^1(G) \cap L^2(G)$. On définit la quantité $\NA{f} \in \R^+ \cup \{\infty\}$, qu'on appellera la {\it norme d'algèbre} de $f$,
par $$\NA{f}= \sup_{g \in L^2(G), \NN{L(g)} \leq 1} \herm{f,g}.$$
\end{definition}

Notons que la borne supérieure définissant $\NA{f}$ peut être infinie: c'est le cas par exemple si $G=\R/\Z$ et $f$ est la fonction caractéristique de l'intervalle $[0,1/2]$.
Par ailleurs, $\NA{f}$ est indépendant du choix de la mesure de Haar choisie (contrairement \`a $\NN{f}_1$,
$\NN{f}_2$, $L(f)$, $\herm{f,g}$, $\NN{L(f)}$ qui sont homogènes de degré $1$, $1/2$, $1$, $1$ et $1$ en la mesure de Haar choisie).
 
\subsubsection{Comparaison avec la norme de Littlewood dans le cas d'une fonction centrale sur un groupe fini}

Dans ce paragraphe, $G$ est un groupe fini muni d'une mesure de Haar $\mu_G$ que nous choisissons de normaliser en demandant que $\mu_G(G)=1$. 
C'est la normalisation qui revient à voir le groupe compact et discret $G$ comme un groupe compact plutôt que discret, point de vue qui sera le seul naturel quand dans nos applications $G$ sera un groupe de Galois. (Cependant, certaines formules concernant la norme de Littlewood ci-dessous 
seraient légèrement plus simples si l'on adoptait le point de vue $G$ discret, avec une mesure de Haar où chaque élément est de mesure $1$).

Avec cette mesure de Haar, si $f,g:G \rightarrow \Cb$ sont des fonctions, on a $\NN{f}_1 = \frac{1}{|G|} \sum_{x \in G} |f(x)|$, 
$\NN{f}_2^2 =   \frac{1}{|G|} \sum_{x \in G} |f(x)|^2$, $\herm{f,g} = \frac{1}{|G|} \sum_{x \in G} \overline{f(x)} g(x)$,
$\pi(f)=\frac{1}{|G|} \sum_x f(x)\pi(x)$ pour $\pi \in \Gd$.

Rappelons qu'on a défini dans l'introduction, pour $f:G \rightarrow \Cb$ une fonction centrale, la transformée de Fourier $\hat{f} : \Gd \rightarrow \Cb$ par la formule $\hat{f}(\pi) = \frac{1}{|G|} \sum_{x \in G} f(x) \tr \pi(x^{-1})$ si bien qu'on a également $\hat{f}(\pi) = \frac{1}{|G|} \sum_{x \in G} f(x) \overline {\tr \pi(x)} = \herm{\chi_\pi,f}$ où $\chi_\pi$ est le caractère de $\pi$, ou encore $\hat{f}(\pi)=\tr \pi(f^\ast)$ où $f^\ast$ est la fonction définie par $f^\ast(x)=f(x^{-1})$. La fonction $\hat{f}$ est aussi définie par la formule $f = \sum_{\pi \in \Gd} \hat{f}(\pi) \chi_\pi = f = \sum_{\pi \in \Gd} \hat{f}(\pi) \chi_\pi.$ Rappelons finalement qu'on a défini la norme de Littlewood
de $f$ par la formule $$\lambda_G(f) = \sum_{\pi \in \Gd} |\hat{f}(\pi)| \dim\pi.$$

\begin{theoreme} Soit $G$ un groupe fini, et $f: G \rightarrow \Cb$ une fonction centrale. Alors
$$\NA{f} = \lambda_G(f).$$
\end{theoreme}

\begin{pf} Nous utilisons le lemme 5.2 de \cite{sanders} qui fait le gros du travail. Ce lemme affirme (pour tout $f: G \rightarrow \Cb$, non nécessairement centrale)
que $\NA{f} = \displaystyle{\sum_{i=1}^{|G|^2}} s_i$, où $s_1,s_2,\dots,s_{|G|^2}$ sont les valeurs singulières de l'opérateur $L(f)$, c'est-à-dire les racines carrées des valeurs propres (répétées selon leur multiplicité comme zéros du polynôme caractéristique) de $L(f)^\ast L(f)$, où $L(f)^\ast=L(f)$ est l'adjoint de $L(f)$.

Si $f$ est centrale, et $\pi$ une représentation irréductible de $G$, l'opérateur $\pi(f)$ commute à $\pi(G)$ donc est la multiplication par un scalaire d'après le lemme de Schur. Ce scalaire est évidemment $\tr \pi(f) / \dim \pi = \herm{\chi_\pi,f^\ast} / \dim \pi$. De même $\pi(f^\ast)$ est l'opérateur scalaire 
$\herm{\chi_\pi,f} / \dim \pi$.

Comme la représentation régulière est la somme directe, pour $\pi \in \Gd$, de $\dim \pi$ copies de la représentation $\pi$, l'opérateur $L(f) L(f)^\ast=L(f) L(f^\ast)$ est la somme directe pour $\pi \in \hat{G}$ de $\dim \pi$ copies de l'opérateur scalaire $\herm{\chi_\pi,f^\ast} \herm{\chi_\pi,f} /(\dim \pi)^2 = 
|\herm{\chi_\pi,f}|^2 / (\dim \pi)^2$
sur l'espace de $\pi$. Les valeurs propres de $L(f) L(f)^\ast$ sont donc les valeurs $|\herm{\chi_\pi,f}|^2 / (\dim \pi)^2$ pour $\pi \in \Gd$ chacune avec la multiplicité $(\dim \pi)^2$ (un facteur $\dim \pi$ parce qu'un opérateur scalaire sur un espace de dimension $\pi$ a $\dim \pi$ fois la même valeur propre, et un facteur $\dim \pi$
puisque la représentation $\pi$ apparaît $\dim \pi$ fois). Les valeurs spectrales de $L_f$ sont les racines carrées des précédentes, à savoir les $|\herm{\chi_\pi,f}|/\dim\pi$ avec multiplicité $(\dim \pi)^2$. Il vient donc
$$\NA{f} = \sum_{\pi \in \Gd} \dim \pi \ |\herm{\chi_\pi,f}|,$$
ce qui est exactement la définition de $\lambda_G(f)$.
\end{pf}

\subsubsection{Compléments}

On peut se demander si l'on peut étendre la définition~\ref{defNL} de la norme de Littlewood au cas de groupes localement compacts plus généraux, et démontrer dans ce contexte l'égalité avec la norme d'algèbre. On peut en tout cas le faire pour une classe relativement restreinte de groupes, à savoir celle des groupes $G$ qui satisfont:
\begin{center} (*) \ \ Toute représentation unitaire irréductible de $G$ est de dimension finie, et $G$ est séparable \end{center}
Un groupe satisfaisant la première partie de cette condition est parfois appelé un {\it groupe de Moore}.
Un groupe compact et un groupe localement compact abélien sont de Moore. Plus généralement, un {\it $Z$-groupe}, i.e. un groupe $G$ tel que le quotient $G/Z(G)$ est compact, est de Moore, et de même un {\it $Z$-groupe virtuel}, i.e. un groupe qui contient un sous-groupe fermé d'indice fini qui est un $Z$-groupe.
Moore a montré (\cite{moore}) que tout groupe de Moore est limite projective de groupes de Lie qui sont des $Z$-groupes virtuels. 
L'intérêt de cette modeste généralisation est d'inclure dans un même formalisme le cas des groupes finis qui nous intéresse ici, et le cas classique depuis Littlewood
des groupes localement compact abéliens (séparables), en premier lieu le cas de $\Z$.

Rappelons (cf. \cite{moore}) qu'un groupe qui satisfait (*) est nécessairement unimodulaire, et de type $1$. Son dual est $\widehat{G} = \coprod_{n \geq 1} \widehat{G}_n$, où 
$\widehat{G}_n$ est l'ensemble des représentations unitaires irréductibles de $G$ de dimension $n$ à équivalence près, muni de sa topologie naturelle (dite topologie de Fell).
On sait alors que la tribu de Mackey sur $\Gd$ est juste sa tribu borélienne, et qu'on dispose d'une mesure canonique sur $\Gd$, la {\it mesure de Plancherel
$\mu_P$}, caractérisé par la formule de Plancherel (cf. \cite[Theorem 18.8.2]{dixmier}):
\begin{eqnarray} \label{plancherel}  \forall f \in L^1(G) \cap L^2(G),\  \int_G |f|^2 d\mu_G = \int_{\Gd} \tr  (\pi(f) \pi(f)^\ast) d\mu_P(\pi).\end{eqnarray}

\begin{definition} Soit $G$ un groupe localement compact satisfaisant (*) et soit $f$ une fonction centrale sur $G$, appartenant à $L^1(G) \cap L^2(G)$. 
On pose $\lambda_G(f) = \int_{\Gd} |\tr \pi(f)| d\mu_P(\pi) \in \R_+ \cup \{\infty\}$ et on appelle cette quantité la {\it norme de Littlewood} de $f$.
\end{definition}
 
Remarquons que si $G$ est fini, $\mu_P(\{\pi\})=\dim \pi$ (cela se voit aisément en appliquant (\ref{plancherel}) à $f=\chi_\pi$), si bien que la définition de $\lambda_G(f)$ redonne dans ce cas celle donnée dans l'introduction.

\begin{prop} Soit $G$ un groupe localement compact satisfaisant (*) et soit $f$ une fonction centrale sur $G$, appartenant à $L^1(G) \cap L^2(G)$. 
Alors on a l'égalité dans $\R_+ \cup \{\infty\}$:
 $$\NA{f}=\lambda_G(f).$$
\end{prop}
\begin{pf} Laissée au lecteur. \end{pf}

\subsection{Propriétés élémentaires de la norme et de la complexité de Littlewood}

Dans toute cette partie, $G$ est un groupe fini, muni de sa mesure de Haar $\mu_G$ de masse totale $1$, et nous développons les sorites concernant la norme de Littlewood d'une fonction centrale sur $G$. Un certain nombre
de ces résultats gardent un sens et restent vrais dans le cadre plus général de la norme d'algèbre d'une fonction quelconque sur $G$. Quand c'est le cas, 
nous donnons le plus souvent deux preuves: une, générale, valable pour la norme d'algèbre d'une fonction quelconque (c'est parfois une simple référence à \cite{sanders}) et une autre, plus concrète et utilisant la théorie des représentations de $G$, valable pour la norme de Littlewood d'une
 fonction centrale (seul cas dont nous aurons besoin  dans cet article). 
 
 \subsubsection{Propriétés algébriques}

\begin{prop} $\NA{\,}$ est une norme sur l'espace des fonctions de $G$ dans $\Cb$. En particulier, $\lambda_G$ est une norme sur l'espace des fonctions centrales de $G$ dans $\Cb$.
\end{prop}
C'est évident.
\begin{cor} \label{lambdaunion} Si $D_1$ et $D_2$ sont deux sous-ensembles de $G$ invariants par conjugaison, $$\lambda(D_1 \cup D_2) \leq \lambda(D_1)+\lambda(D_2) + \lambda(D_1 \cap D_2).$$ Si $D_1$ et $D_2$ sont disjoints, $$|\lambda(D_1)-\lambda(D_2)| \leq \lambda(D_1 \cup D_2) \leq \lambda(D_1)+\lambda(D_2).$$
\end{cor}
\begin{pf}
En effet, on a sous ces hypothèses $\un_{D_1+D_2}=\un_{D_1} + \un_{D_2} - \un_{D_1 \cap D_2}$, d'où la première assertion. En particulier,
si $D_1$ et $D_2$ sont disjoints, $\un_{D_1+D_2}=\un_{D_1} + \un_{D_2}$, d'où la seconde.
\end{pf}

\begin{remarque} \label{suppspec} Notons qu'on n'a pas en général égalité $\lambda(D_1 \cup D_2) = \lambda(D_1) + \lambda(D_2)$ si $D_1$ et $D_2$ sont disjoint.
On a en revanche évidemment cette égalité quand $D_1$ et $D_2$ sont disjoints et de {\it support spectral} disjoints.
\end{remarque}

\begin{cor} \label{lambdacomplement}  Si $D$ est un sous-ensemble de $G$ invariant par conjugaison, $\lambda(D)-1 \leq \lambda(G-D) \leq \lambda(D)+1$. \end{cor}
\begin{pf}
On a $\lambda(G)=1$ puisque $\widehat{\un_G}(\pi)=0$ sauf si $\pi$ est la représentation triviale, auquel cas $\widehat{\un_G}(\pi)=1$. On applique alors le corollaire précédent à $D_1=D$, $D_2 = G-D$.
\end{pf} 

\begin{prop} \label{invaut} La norme $\NA{\ }$ est invariante par automorphismes, anti-automorphismes et par translations de $G$.
En particulier, si $\sigma: G \rightarrow G$ est soit un automorphisme de groupe, soit un anti-automorphisme, soit de la forme $x \mapsto zx$ pour $z$ fixé dans le centre $Z(G)$ de $G$, et
si $f: G \rightarrow \Cb$ est centrale, alors $\lambda(f)=\lambda(f \circ \sigma)$, et si $D \subset G$  est invariant par conjugaison,
 $\lambda(\sigma(D)) = \lambda(D)$.
\end{prop}
\begin{pf} C'est évident sur la définition de $\NA{f}$. 
Voici une preuve directe du cas particulier concernant  la norme de Littlewood d'une fonction centrale:
 si $\psi$ est un automorphisme, ou un antiautomorphisme de $G$ c'est trivial, et si $\psi$ est la multiplication par $z \in Z(G)$, il suffit de remarquer
que pour tout $\pi \in \Gd$, $\pi(z)$ est un scalaire par le lemme de Schur, et une racine de l'unité, donc de module 1, et que
$\widehat{f\circ\sigma}(\pi) = \overline{\pi(z)} \widehat{f}(\rho)$.
\end{pf}

\begin{prop}On a $\NA{fg} \leq \NA{f} \NA{g}$, autrement dit $\NA{\ }$ est une norme d'algèbre. En particulier,
pour $f,g$ centrale, $\lambda_G(fg) \leq \lambda_G(f) \lambda_G(g)$. \end{prop}
\begin{pf} Ceci est prouvée par Sanders: cf. \cite{sanders}, Prop. 5.4.\end{pf}

\begin{cor} \label{cor1D} Pour tout sous-ensemble non vide $D$ de $G$, $\NA{\un_D} \geq 1$. En particulier, pour $D$ invariant par conjugaison, $\lambda_G(D) \geq 1$. \end{cor}
\begin{pf} En effet $\un_D^2=\un_D$, d'où $\NA{\un_D} \leq \NA{\un_D}^2$ et $ \NA{\un_D} \neq 0$ si $D$ est non vide. \end{pf}

\subsubsection{Minoration}

\begin{prop} \label{propmin} On a $\NA{f} \geq \NN{f}_\infty $. En particulier, si $f$ est centrale, $\lambda(f) \geq \NN{f}_\infty $.
\end{prop} 
\begin{pf} 
C'est le lemme 5.3 de \cite{sanders}, mais nous donnons une preuve plus directe dans le cas général.
Pour $x \in G$, $L(\un_x)$ est $\frac{1}{|G|}$ fois la translation à gauche par $x$ sur $L^2(G)$ donc la norme d'opérateur $\NN{L(\un_x)}$ est $\frac{1}{|G|}$.
Pour $f$ une fonction sur $G$, écrivons $f = \sum_{x \in G} f(x) \un_{x}$, d'où $\NN{L(f)} \leq \sum_{x \in G} \NN{f(x) L(\un_x)} = \frac{1}{|G|} \sum_{x \in G} |f(x)| =  \NN{f}_1$. On a donc $\NA{f}= \sup_{g, \NN{L(g)} \leq 1} |\herm{f,g}| \geq \sup_{g, \NN{g}_1 \leq 1} |\herm{f,g}| = \NN{f}_{\infty}$.
\end{pf}

On retrouve en particulier le corollaire~\ref{cor1D}: $\NA{\un_D} \geq 1$ si $D$ est non vide. Les cas d'égalité ont été déterminés par Sanders:
\begin{theoreme} \label{thmin} Soit $D$ un sous-ensemble non vide de $G$. On a l'égalité $\NA{\un_D} \geq 1$ si et seulement
si $D$ est une classe à gauche\footnote{ou, ce qui revient bien sûr au même, une classe à droite pour un (autre) sous-groupe de $G$.} pour un sous-groupe $H$ de $G$. En particulier, si $D$ est invariant par conjugaison, on a $\lambda(D)=1$ si et seulement si $D$ est de la forme $aH$ avec $H$ un sous-groupe distingué de $G$, et $a$ un élément de $G$ dont l'image canonique dans $G/H$ est centrale. \end{theoreme}
\begin{pf} La première partie est prouvée dans les sections $6$ et $7$ de \cite{sanders}. La seconde en résulte si l'on remarque qu'un classe à gauche $aH$
est invariante par conjugaison si et seulement si $H$ est distingué dans $G$ et l'image de $a$ dans $G/H$ est centrale, ce qui est facile: si $aH$ est invariant par conjugaison, $H=\{xy^{-1}, x \in aH, y \in aH\}$ l'est aussi, et l'image de $aH$, i.e. l'image de $a$, dans $G/H$ l'est également; la réciproque est triviale.

Comme la preuve de Sanders est longue et difficile, nous donnons une autre preuve, plus simple, de la seconde assertion, i.e. du cas où $D$ est stable par conjugaison. Nous re-prouvons du même coup dans ce cas l'inégalité $\lambda(D) \geq 1$.

On a d'après l'égalité de Parseval $$\frac{|D|}{|G|} = \NN{\un_D}_2^2 = \sum_{\pi\in \Gd} |\widehat{\un_D}(\pi)|^2.$$
D'autre part, pour $\pi \in \Gd$, on a $| \widehat{\un_D}(\pi)| = \frac{1}{|G|} \left| \sum_{d \in D} \tr(\pi(d))\right| \leq \frac{|D|}{|G|} \dim\pi$, avec égalité si et seulement si $\pi(d)$ est un scalaire indépendant de $d$ pour tout $d \in D$. On obtient donc
\begin{eqnarray*} \frac{|D|}{|G|} \leq  \sum_{\pi\in \Gd} |\widehat{\un_D}(\pi)| \frac{|D|}{|G|} \dim\pi =  \frac{|D|}{|G|} \lambda(D) \end{eqnarray*}
où encore $1 \leq \lambda(D)$, avec égalité si et seulement si, $D$ satisfait la condition suivante:
\begin{center} ($*$) \ \ pour tout $\pi \in \Gd$, on a $\begin{cases} 
\text{ soit} & \widehat{\un_D}(\pi)=0 \\ \text{ soit} & \pi(d) \text{ est un scalaire indépendant de $d$ pour tout $d \in D$ }\end{cases}$ \end{center}
Cette condition ($*$) est satisfaite si $D=aH$ avec $H$ distingué et $a$ central dans $G/H$, car pour les $\pi$ tels que $\pi(H) \neq 1$, on a 
$\widehat{\un_D}(\pi) = 0$, et pour les autres, $\pi(d) = \pi(a)$ est indépendant de $d$ et est un scalaire par le Lemme de Schur, puisque $\pi$ se factorise par une représentation irréductible de $G/H$ et l'image de $a$ dans $G/H$ est centrale. Réciproquement, si $D$ satisfait ($*$), 
définissons $H$ comme l'intersection des sous-groupes $\ker \pi$ pour $\pi$ tels que $\widehat{\un_D}(\pi) \neq 0$; 
c'est un sous-groupe normal de $G$, et pour toute représentation $\pi$ de $G/H$, et tout $d$ dans l'image de $D$ dans $G/H$, $\pi(d)$ est un scalaire indépendant de $d$. Le lemme suivant, appliqué à $G/H$, montre alors que l'image de $D$ dans $G/H$ est réduite à un seul élément $a$, qui plus est dans le centre de $G/H$, ce qui termine la preuve.
\end{pf}
\begin{lemme} \label{lemmescalaire} Soit $G$ un groupe fini. Les seuls éléments $g$ de $G$ qui sont tels que $\pi(g)$ est scalaire pour toute représentation irréductible de
$G$ sont les éléments du centre de $G$. De plus, pour de tels éléments, $g$ est uniquement déterminé par l'application  $\Gd \rightarrow \Cb^\ast,\ \pi \mapsto \pi(g)$.
\end{lemme}
\begin{pf}
On a un isomorphisme d'algèbre $\prod_\pi \pi : \Cb[G] \rightarrow \prod_{\pi} \End(V_\pi)$, où $V_\pi$ est l'espace de la représentation $\pi$, et l'hypothèse implique que l'image de $g$ par cette application est centrale, donc que $g$ est dans le centre de $\Cb[G]$, ce qui implique qu'il est dans le centre de $G$. La deuxième assertion provient de ce que la somme des représentations irréductibles $\pi$ de $G$ est fidèle.
\end{pf}

\subsubsection{Majoration}

 \begin{prop}[Majoration par Cauchy-Schwarz] \label{majcs}
 Pour $f$ une fonction complexe sur $G$, on a $\NA{f} \leq  \sqrt{|G|} \NN{f}_2$, avec égalité si et seulement si le support de $f$ est vide ou réduit à un point.
 En particulier, si $f$ est centrale, $\lambda(f) \leq  \sqrt{|G|} \NN{f}_2$ avec égalité si et seulement si le support de $f$ est vide ou réduit à un point (du centre de $G$).  Si $D$ est un sous-ensemble de $G$ stable par conjugaison, $\lambda(D) \leq \sqrt{|D|}$ avec égalité si et seulement si $D$ est vide ou un singleton
 (inclus dans le centre de $G$). \end{prop}

\begin{pf} En effet, si $\NN{L_f}_{HS}$ dénote la norme de Hilbert-Schmidt de l'opérateur $L_f$, on a
$\NN{L_f} \geq \frac{1}{\sqrt{|G|}} \NN{L_f}_{HS}$ d'après les lemmes 3.3 et 3.4 de \cite{sanders},
et $\NN{L_f}_{HS}=\NN{f}_2 $ d'après \cite[Th. 4.2]{sanders}. D'où par dualité, $\NA{f} \leq \sqrt{|G|} \NN{f}_2$. Nous laissons l'analyse du cas
d'égalité au lecteur.
\par \medskip
Voici une preuve directe dans le cas où $f$ est centrale, qui justifie le nom de "majoration par Cauchy-Schwarz". 
Comme $f = \sum \widehat{f}(\pi) \chi_\pi$, et que les $\chi_\pi$ forment une base orthonormée de $L^2(G)$,
on a $$\NN{f}_2^2 = \sum_{\pi} |\widehat{f}(\pi)|^2.$$ 
Par ailleurs et par Cauchy-Schwarz $$\lambda(f) = \sum_\pi |\hat{f}(\pi)| \dim \pi \leq \sqrt{ \sum_\pi | \hat{f}(\pi)|^2 } \sqrt{\sum_\pi (\dim \pi)^2}.$$
La première racine carrée est égale à $\NN{f}_2$, la seconde à $\sqrt{|G|}$ et l'inégalité est prouvée. 
Pour avoir égalité, il faut et il suffit
que le vecteur $(\widehat{f}(\pi))_{\pi \in \Gd}$ soit proportionnel au vecteur $(\dim \pi)_{\pi \in \Gd}$, ce qui est équivalent à l'assertion $\forall x \in \supp f, \ \forall \pi \in \Gd, \pi(x)$ est un scalaire indépendant de $x$. Le lemme~\ref{lemmescalaire} montre que ceci est encore équivalent à l'assertion que $\supp f$ est vide ou un singleton contenu dans le centre.  

Enfin, l'assertion pour $D$ découle de la précédente puisque $\NN{\un_D}_2=\frac{\sqrt{|D|}}{\sqrt{|G|}}$.
\end{pf}
\begin{remarque} La quantité $\sqrt{|G|} \NN{f}_2$ s'écrirait simplement $\NN{f}_2$ si l'on avait choisi la mesure de Haar sur $G$ qui donne une masse $1$ à chaque élément. Cette mesure est la mesure naturelle si l'on considère le groupe fini $G$ comme un groupe discret, tandis que celle que nous avons choisie
revient à considérer $G$ comme un groupe compact, ce qui est plus naturel  dans les applications où $G$ est un groupe de Galois.
\end{remarque}
La majoration de Cauchy-Schwarz $$\lambda(\un_D) \leq \sqrt{|D|}$$ sera l'étalon auquel nous mesurerons les autres majorations de $\lambda(D)$ que nous 
obtiendrons. 
\begin{cor}[Majoration triviale] \label{majtriv} On a  pour tout $f:G \rightarrow \Cb$, $\NA{f} \leq |G| \NN{f}_1$. En particulier $\lambda(f) \leq |G|\NN{f}_1$ si $f$ est centrale,
et $\lambda(D) \leq |D|$. Les cas d'égalité sont les mêmes que dans la proposition.
\end{cor}

\subsubsection{Comportement par produit}

Soit $G_1,G_2$ deux groupes finis, $f_1 : G_1 \rightarrow \Cb$ et $f_2: G_2 \rightarrow \Cb$ deux fonctions. On définit la fonction
$f_1 \otimes f_2 : G_1 \otimes G_2 \rightarrow \Cb$ par la formule $(f_1 \otimes f_2)(g_1,g_2)=f_1(g_1) f_2(g_2)$.\begin{prop} \label{compprod}
On a $$\NAA{f_1 \otimes f_2}{G_1 \times G_2}=\NAA{f_1}{G_1}\NAA{f_2}{G_2}.$$
Si $f_1$ et $f_2$ sont invariantes par conjugaison, $f_1 \otimes f_2$ l'est également, et l'on a 
$$\lambda_{G_1\times G_2}(f_1 \otimes f_2)=\lambda_{G_1}(f_1) \lambda_{G_2}(f_2).$$
En particulier, si $D_1 \subset G_1$ et $D_2 \subset G_2$ sont deux sous-ensembles  stables par conjugaison, on a
$\lambda_{G_1\times G_2}(D_1 \times D_2)=\lambda_{G_1}(D_1) \lambda_{G_2}(D_2).$
\end{prop}
\begin{pf} C'est évident sur la définition de la norme d'algèbre puisque $L^2(G_1 \times G_2) = L^2(G_1) \otimes L^2(G_2)$. 
C'est d'ailleurs tout aussi évident pour une fonction centrale $f$, avec la norme de Littlewood puisque l'application $(\pi_1,\pi_2) \mapsto \pi_1 \otimes \pi_2$
identifie $\widehat{G_1} \times \widehat{G_2}$ avec $\widehat{G_1 \times G_2}$, et que la transformé de Fourier de $f_1 \otimes f_2$ est juste $\widehat{f_1} \otimes \widehat{f_2}$ modulo cette identification.
\end{pf}

\subsubsection{Comportement par passage au quotient}

\begin{prop} \label{compquot} Soit $G$ un groupe fini, $U$ un sous-groupe distingué de $G$, $s$ la surjection canonique $G \rightarrow G/U$,
et $f_U$ une fonction centrale sur $G/U$. Alors $\lambda_G(f \circ s)=\lambda_{G/U}(f)$. 
\end{prop}
\begin{pf}
L'application $\pi \mapsto \pi \circ s$ identifie $\widehat{G/U}$ à un sous-ensemble de $\Gd$, et l'on a $\widehat{f \circ s} = \widehat{f}$ sur
 $\widehat{G/U}$, $\widehat{f \circ s}=0$ sur le complement de $\widehat{G/U}$ dans $\widehat{G}$. La proposition s'en suit.
\end{pf}

\begin{cor} \label{compquotcor} Si $D_U$ est une partie de $G/U$ stable par conjugaison, et $D$ son image inverse dans $G$, alors $\lambda_G(D)=\lambda_{G/U}(D_U)$.
\end{cor}

\begin{remarque} On retrouve le fait que si $D=a U$ avec $U$ normal et $a$ central dans $G/U$, $\lambda_G(aU)=1$, puisque le corollaire~\ref{compquotcor} et la proposition~\ref{invaut} donnent
$\lambda_G(aU) = \lambda_{G/U}(\{a\})=\lambda_{G/U}(\{1\})=1$.
\end{remarque}

\subsubsection{Comportement par restriction à un sous-groupe}

\begin{prop} \label{compres} Soit $G$ un groupe fini, $f : G \rightarrow \Cb$ une fonction centrale, et $H$ un sous-groupe de $G$. 
On a alors $\lambda_H(f_{|H}) \leq \lambda_G(f)$. En particulier, si $D$ est un sous-ensemble de $G$ stable par conjugaison,
$\lambda_H(D \cap H) \leq \lambda_G(D)$.
\end{prop}
\begin{pf} Écrivons $f=\sum_{\pi\in \Gd} c_\pi \tr \pi$. Alors $f_{|H} = \sum_{\pi \in \Gd} c_{\pi} \tr \pi_{|H} = \sum_{\rho \in \Hd} \sum_{\pi \in \Gd} c_\pi m(\rho,\pi) \tr \rho$, où pour $\pi$ une représentation irréductible de $G$, $\rho$ une représentation irréductible de $H$, $m(\rho,\pi)$ est la multiplicité de $\rho$ dans $\pi_{|H}$. On 
a donc \begin{eqnarray*} \lambda_H(f_{|H}) &=& \sum_{\rho \in \Hd} |\sum_{\pi \in \Gd} c_{\pi} m(\rho,\pi)| \dim \rho \\
& \leq & \sum_{\pi \in \Gd} |c_{\pi}| \sum_{\rho \in \Hd} m(\rho,\pi) \dim \rho \\ & =&  \sum_{\pi \in \Gd} |c_\pi| \dim \pi = \lambda_G(f)\end{eqnarray*}
 \end{pf}

\begin{remarque} On obtient une nouvelle preuve de ce que $\lambda_G(f_1 f_2) \leq \lambda_G(f_1) \lambda_G(f_2)$ pour deux fonctions centrales $f_1$ et $f_2$ sur $G$, en remarquant que  $f_1 f_2$ est la restriction de $f_1 \otimes f_2$ au sous-groupe diagonal $G$ de $G \times G$.
\end{remarque}

\subsubsection{Comportement par induction}

Rappelons que pour $G$ un groupe fini, $H$ un sous-groupe de $G$, et $f$ une fonction centrale sur $H$, on définit une fonction centrale sur $G$,
$\Ind_H^G f $,  par la formule:
\begin{eqnarray} \label{indf} (\Ind_H^G f)(x) = \frac{1}{|H|} \sum_{y \in G, \ yxy^{-1} \in H} f(yxy^{-1}) \end{eqnarray}
\begin{prop} \label{compind} Soit $G$ un groupe fini, $H$ un sous-groupe de $G$, et $f$ une fonction centrale sur $H$. 
Alors $$\lambda_G(\Ind_H^G f) \leq  \frac{|G|}{|H|} \lambda_H(f).$$ On a égalité si et seulement si la condition suivante est satisfaite:
 pour tout $\pi \in \Gd$, la fonction $\hat{f}$ a un argument complexe constant sur l'ensemble des sous-représentations irréductibles de la restriction de
 $\pi$ à $H$.
  \end{prop}
Précisons ce que nous voulons dire par "argument complexe constant". Une application $h: X \rightarrow \Cb$ a {\it un argument complexe constant}
s'il existe $\theta \in [0,1)$, tel que $\forall x \in X$, $h(x) = |h(x)| e^{2 i \pi \theta}$.  Autrement dit, la restriction  de $x \mapsto h(x)/|h(x)|$ est constante sur le sous-ensemble de $X$ des $x$ tels que $h(x) \neq 0$.

\begin{pf} Comme dans la preuve de la proposition~\ref{compres}, pour $\pi \in \Gd$, $\rho \in \Hd$  on écrit $m(\rho,\pi)$ pour la multiplicité de $\rho$ dans 
$\pi_{|H}$, ou ce qui revient au même, pour la multiplicité de $\pi$ dans $\Ind_H^G \rho$. Pour $\rho \in \Hd$ fixé, on a donc
$\sum_{\pi \in \Gd} m(\rho,\pi) \dim \pi = \dim \Ind_H^G \rho = \frac{|G|}{|H|} \dim \rho$.

Écrivons $f = \sum _{\rho \in \Hd} \hat{f}(\rho) \tr \rho$, d'où $$\Ind_H^G f =  \sum_{\rho \in \Hd} \hat{f}(\rho) \tr \Ind_H^G \rho = \sum_{\rho,\pi} \hat{f}(\rho) m(\rho,\pi) \dim \pi,$$ 
et donc \begin{eqnarray*} \lambda_G(\Ind_H^G f)  &=& \sum_{\pi \in \Gd}  \left| \sum_{\rho \in \Hd}  \hat{f}(\rho) m(\rho,\pi) \right | \dim \pi \\
& \leq & \sum_{\pi \in \Gd}  \sum_{\rho \in \Hd}   |\hat{f}(\rho)| m(\rho,\pi)\dim \pi   \\
&=& \sum_{\rho \in \Hd} |\hat{f}(\rho)| \sum_{\pi \in\Gd} m(\rho,\pi) \dim \pi\\
&=& \sum_{\rho \in \Hd}  |\hat{f}(\rho)| \frac{|G|}{|H|} \dim \rho\\
&=& \frac{|G|}{|H|} \lambda_H(f)
\end{eqnarray*}

Pour qu'on ait égalité, il faut et il suffit que pour tout $\pi \in \Gd$, on ait 
$$  \left| \sum_{\rho \in \Hd}  \hat{f}(\rho) m(\rho,\pi) \right |  = \sum_{\rho \in \Hd}   |\hat{f}(\rho)| m(\rho,\pi)$$
ce qui est équivalent à ce que l'argument complexe de $\hat{f}$ soit constant sur l'ensemble des $\rho$ tels que $m(\rho,\pi) \neq 0$, ce qui est la condition de l'énoncé.
\end{pf}

 \subsubsection{La méthode de Serre}
 
 Le résultat suivant est en quelques sorte la traduction dans le langage de la complexité de Littlewood d'une technique employée par Serre (cf. \cite[\S2.7]{serre}) pour obtenir des majorations améliorées du terme d'erreur dans le théorème de Chebotarev.
 
 \begin{theoreme} \label{methserre}
 Soit $G$ un groupe fini, $D$ un sous-ensemble de $G$ stable par conjugaison. Soit $H$ un sous-groupe de $G$, et $U$ un sous-groupe normal de $H$,
 sur lesquels on fait les hypothèses suivantes, pour lesquelles on note $C(d)$ la classe de conjugaison dans $G$ d'un élément $d$ de $D$:
 \begin{itemize}
\item[(a)] $|C(d)|$ est indépendant de $d \in D$
\item[(b)] $|C(d) \cap H|$ est indépendant de $d \in D$.
\item[(c)] $U(D \cap H)=D \cap H$.
\end{itemize}
Notons $s$ la projection canonique de $H$ sur $H/U$.
Alors 
$$\lambda_G(D) \leq \frac{|C(d)|}{|C(d) \cap H|} \lambda_{H/U}(s(D \cap H)),$$
où $d$ est un élément quelconque de $D$.
De plus, on a égalité si et seulement si la condition suivante est satisfaite:
\begin{itemize}
\item[(d)] Pour tout $\pi \in \Gd$, la fonction $\widehat{\un_{s(D \cap H)}}$ a un argument complexe constant sur l'ensemble des
représentations irréductibles $\rho$ de $H/U$  telles que $\rho \circ s$ apparaît dans $\pi_{|H}$.
\end{itemize}
 \end{theoreme}
 \begin{pf} 
Par la formule~(\ref{indf}), $\Ind_H^G (\un_{D \cap H})(x) = \frac{1}{|H|} |\{ y \in G, y x y^{-1} \in D \cap H\}|.$
L'ensemble $\{ y \in G, y x y^{-1} \in D \cap H\}$ est vide si $x \not \in D$. Pour $x=d \in D$, le nombre d'éléments
de la forme $y d y^{-1}$ qui sont dans $D \cap H$ est $|C(d) \cap H|$ et chacun de ces éléments est obtenu
pour $|Z_G(d)|$ valeurs de $y$, où $Z_G(d)$ est le centralisateur de $d$ dans $G$. On a $|Z_G(d)|=|G|/|C(d)|$, d'où
il vient 
$\Ind_H^G (\un_{D \cap H})(x) =  \frac{|C(d) \cap H| |G|}{|C(d)||H|}$ si $x=d \in D$, $0$ si $x \not \in D$. Comme la valeur  $\frac{|C(d) \cap H| |G|}{|C(d)||H|}$
est constante par les hypothèses (a) et (b), on a
$$\un_D = \frac{|C(d)||H|}{|C(d) \cap H||G|} \Ind_{H}^G \un_{D \cap H}.$$
Par l'hypothèse (c) et la proposition~\ref{compquot}, on a $\lambda_H(D \cap H) = \lambda_{H/U} (s(D \cap H))$, et 
le support de $\widehat{\un_{D \cap H}}$ ne contient que des représentations irréductibles de $H$ qui sont triviales sur $U$.
La proposition~\ref{compind} donne donc:
$$\lambda_G(D) = \frac{|C(d)|}{|C(d) \cap H|} \lambda_{H/U}(s(D \cap H))$$ avec égalité
quand l'hypothèse (d) est satisfaite. Ceci termine la preuve.
 \end{pf}

\subsection{Exemples de calculs de complexités de Littlewood dans le cas abélien}

Les calculs ci-dessous ont deux buts: préparer certaines applications arithmétiques qui seront développées dans le reste de cet article ou dans des articles ultérieurs, et permettre de dégager une idée intuitive de ce que mesure la complexité de Littlewood d'un sous-ensemble $D$ d'un groupe abélien. Je serais tenté d'exprimer cette idée ainsi: plus l'ensemble $D$ est facile à décrire en terme de la structure de groupe de $G$, plus sa complexité de Littlewood $\lambda(D)$ est petite.
 
\subsubsection{Progressions arithmétiques dans $\Z/n\Z$}

\begin{prop} Il existe une constante absolue $\ccl{propint}$ telle que pour toute progression arithmétique $D=\{a + i r \ \mid \  i=0,\dots,x\}$
avec $a,r \in \Z/n\Z$, $x \geq 0$ un entier, on ait
on ait $\lambda_{\Z/n\Z}(D) \leq \ccr{propint} \log |D|$.
\end{prop}
C'est un calcul facile et classique, laissé au lecteur.

\subsubsection{Intervalles dans $(\Z/n\Z)^\ast$}

Soit $x,n$ des entiers tels que $1 \leq x < n$ 
 et soit $D \subset (\Z/n\Z)^\ast$ l'ensemble des résidus modulo $n$ des entiers positifs plus petits que $x$ et  premiers à $n$.
Par Cauchy-Schwarz, $\lambda_{(\Z/n\Z)^\ast}(D) \leq \sqrt{x}$ et cette borne n'est pas loin d'être optimale, du moins quand $x$ n'est pas trop petit par rapport à $n$.
\begin{prop} Sous (GRH), pour tout $\epsilon >0$ il existe une constante $\ccl{cint}>0$ telle que
$$\lambda_{(\Z/n\Z)^\ast}(D) \geq \ccr{cint} \sqrt{x} n^{-\epsilon}.$$
\end{prop}
\begin{pf} Soit $\chi$ un caractère de $(\Z/n\Z)^\ast$. On a, pour $\chi \neq 1$, l'estimation sous (GRH) bien connue
$\widehat{\un_D}(\chi) = \frac{1}{\phi(n)} \sum_{a=1}^x \chi(a) \leq \frac{1}{\phi(n)} \ccr{cint} \sqrt{x} n^\epsilon$ où $\ccr{cint}$ est une constante absolue.
D'où \begin{eqnarray*} x &=&
 \sum_{\chi \in (\Z/n\Z)^\ast} |\widehat{\un_D}(\chi)|^2 \\ & \geq & \sum_{\chi \in 
 (\Z/n\Z)^\ast,\chi \neq 1} |\widehat{\un_D}(\chi)| \frac{1}{\phi(n)} \ccr{cint} \sqrt{x} n^{\epsilon} \\ 
 &=& \ccr{cint} \sqrt{x} n^{\epsilon} (\lambda(D) - |\widehat{\un_D}(1)| )\end{eqnarray*}
 ce qui prouve le résultat puisque $ |\widehat{\un_D}(1)| = |D|/\phi(n) <1$
\end{pf}

\subsubsection{Sous-ensemble des générateurs d'un groupe cyclique}

Pour $n \geq 1$ un entier, on note $\phi(n)$ la fonction d'Euler et $\omega(n)$ le nombre de diviseurs premiers (comptés sans multiplicité) de $n$. 
\begin{prop}  \label{lambdagen} Soit $G=\Z/n\Z$ et soit $D$ l'ensemble des générateurs de $G$. Alors $\lambda(D)=2^{\omega(n)} \frac{\phi(n)}{n}$.
\end{prop}
\begin{pf}
Par une application de la formule d'inclusion-exclusion, on a $\un_D = \sum_{d \mid n} \mu(d) \un_{d \Z/n \Z}$, d'où $\widehat{\un_D} = 
\sum_{d \mid n} \mu(d) \widehat{\un_{d \Z/n \Z}}$. Pour $a$ un diviseur de $n$, soit $\chi$ un caractère de $\Z/n\Z$ tel que $\ker \chi =a \Z/n\Z$.
On voit directement sur la définition $$\widehat{\un_{d \Z/n\Z}}(\chi) = \begin{cases}  1/d & \text{ si $a \mid d$} \\
0 & \text{sinon} \end{cases}$$
d'où
$$\widehat{\un_D}(\chi) =  \sum_{d,\ a \mid d \mid n} \mu(d)/d.$$
Soit $P$ l'ensemble des facteurs premiers de $n$.
Si $a$ a un facteur carré, tous les $d$ apparaissant dans la somme ci-dessus en ont aussi un et $\mu(d)=0$, et donc $\widehat{\un_D}(\chi)  = 0$.
Si $a$ n'a pas de facteurs carré, on peut écrire $a = \prod_{p \in S} p$ pour un unique sous-ensemble $S$ de $P$, et on voit qu'alors 
$\widehat{\un_D}(\chi) = (-1)^{|S|} \prod_{p \in P-S} (p-1)  / \prod_{p \in P} p$.

Pour chaque diviseur $a$ de $n$, il y a exactement $\phi(a)$ caractères $\chi$ tel que $\ker \chi = a \Z/n\Z$, puisque ces derniers s'identifient aux caractères fidèles de $\Z/n\Z$. Si $a = \prod_{p \in S} p$, $\phi(a) = \prod_{p \in S} (p-1)$. On obtient donc
$$\lambda(D) = \sum_{S \subset P} \prod_{p \in P} (p-1)/p = 2^{|P|} \phi(n)/n,$$
d'où le résultat puisque $|P| = \omega(n)$.
\end{pf}

\subsubsection{Sous-tores de $\Gl_n$}
\label{lambdatore}

Dans ce \S, nous nous contentons de poser un problème que nous ne savons pas résoudre, sauf dans deux cas particuliers. Sa solution générale
serait pourtant très utile dans les applications à la conjecture de Lang-Trotter, cf. \S\ref{syscompLT}.

Soit $N$ un entier $\geq 1$, et soit $T$ un sous-tore de dimension $r$ du schéma en groupe $\Gl_n$ sur $\Spec \Z[1/N]$. Soit $\Lambda$ l'ensemble des nombres premiers $\ell$ ne divisant pas $N$
tels que $T_{\F_\ell} = T \otimes_{\Spec Z[1/N]} \Spec \F_\ell$ soit un tore {\bf déployé}. Comme $\Lambda$ est l'ensemble des $\ell$ qui sont totalement décomposé dans l'extension finie de $\Q$ fixée par le noyau de l'application naturelle $\GQ \rightarrow \Gl(X^\ast(T))$, $\Lambda$ est un ensemble frobénien de densité positive, en particulier infini. 

Soit $a$ un élément de $\Z$. Pour tout nombre premier $\ell$, on note $T(\F_\ell)^{\ta,\dr}$
l'ensemble des $g \in T(\F_\ell)$ qui vus comme éléments de $\Gl_n(\F_\ell)$, sont diagonalisables réguliers (i.e. à valeurs propres distinctes)
et de trace $a \pmod{\ell}$.

\begin{question} \label{questiontore} Peut-on donner une estimation asymptotique de $\lambda_{T(\F_\ell)}(T(\F_\ell)^{\ta,\dr})$ quand $\ell$ tend vers l'infini tout en restant dans $\Lambda$ ? Dans quel cas cette estimation est-elle d'un ordre de grandeur meilleur que l'estimation de Cauchy-Schwarz $\lambda(T(\F_\ell)^{\ta,\dr}) \leq \sqrt{ | T (\F_\ell)^{\ta,\dr} | } $ ? 
\end{question}
Notons que si $T(\F_\ell)^{\ta,\dr}$ est non vide, alors $\ell \in \Lambda$ d'où la restriction à ces $\ell$ dans l'énoncé du problème.

\begin{theoreme} \label{torediagonal} Soit $T$ le sous-tore diagonal de $\Gl_n$ sur $\spec \Z$.
Alors quand $\ell$ tend vers l'infini, $\lambda( T(\F_\ell)^{\ta,\dr}) \sim \ell^{(n-1)/2}$ si $a \neq 0$, et $\lambda(T(\F_\ell)^{\tz,\dr}) 
\sim \ell^{(n-2)/2}$.
\end{theoreme} 
\begin{pf}
Soit $D_a$ l'ensemble des matrices de trace $a$ dans $T(\F_\ell)$ et $D_a^\nr$  l'ensemble des éléments non réguliers de $D_a$. 
Comme, toue les éléments de $D_a$ sont diagonalisables, on a $D_a =  T(\F_\ell)^{\ta,\dr} \coprod D_a^\nr$. 
On a $|D_a^{\nr}| = O(\ell^{n-2})$ quand $\ell$ tend vers l'infini, d'où $\lambda(D_a^\nr) = O(\ell^{(n-2)/2})$ par la majoration de Cauchy-Schwarz (cf. prop.~\ref{majcs}); quand $a=0$ on a même $\lambda(D_0^\nr) = O( \ell^{(n-3)/2})$ puisque dans ce cas, si $Z$ est le centre de $\Gl_n(\F_\ell)$,
on voit que $Z D_0^{\nr}=D_0^{\nr}$ et donc, par la prop.~\ref{compquot}, $\lambda_{T(\F_\ell)} = \lambda_{T(\F_\ell)/Z}(D_0^{\nr}/Z)$ auquel on peut ensuite appliquer la majoration de Cauchy-Schwarz. On voit donc que $\lambda(D_a^\nr)$
est dans tous les cas négligeable par rapport à l'équivalent de $\lambda( T(\F_\ell)^{\ta,\dr})$ qu'on veut prouver, et donc qu'il suffit de prouver cet
équivalent pour $\lambda(D_a)$ au lieu de $\lambda( T(\F_\ell)^{\ta,\dr})$

Soit $\chi=(\chi_1,\dots,\chi_n)$ un caractère $T(\F_\ell) \rightarrow \Cb^\ast$. On a par définition
$$\widehat{\un_{D_a}}(\chi) = \frac{1}{(\ell-1)^n} \sum_{\substack{x_1,\dots,x_n \in \F_q^\ast\\ x_1+\dots+x_n=a }} \chi_1(x_1) \dots \chi_n(x_n), $$
et on reconnait dans la somme un somme de Jacobi généralisée $J_a(\chi_1,\dots,\chi_n)$, comme définie dans  \cite{irelandrosen}. 
On voit par un changement le changement de variable $x_i \mapsto a x_i$ que si $a \neq 0$, $|J_a(\chi_1,\dots,\chi_n)|=|J_1(\chi_1,\dots,\chi_n)|$.
Si $a \neq 0$ dans $\F_q$, et si les $\chi_i$ ne sont pas tous triviaux, on a d'après le théorème 4 de \cite[Chapter 8]{irelandrosen}, 
$|J_1(\chi_1,\dots,\chi_n)|=\ell^{(n-1)/2}$ sauf si $\chi_1\dots\chi_n=1$ auquel cas $  |J_1(\chi_1,\dots,\chi_n)|=\ell^{n/2-1}$. Comme il y a $(\ell-1)^n$ caractères $\chi$, on voit donc que la contribution des caractères $\chi$ tels que $\chi_1 \dots \chi_n=1$ est négligeable et que 
$\lambda_G(D_a) \sim \ell^{(n-1)/2}$ si $a \neq 0$ dans $\F_\ell$.
Bien sûr, si $a \neq 0$ dans $\Z$, on aura $a \neq 0$ dans $\F_\ell$ pour tout $\ell>a$, d'où le résultat dans ce cas.

Supposons au contraire $a=0$ dans $\Z$, dans dans $\F_\ell$ pour tout $\ell$. Si $\chi_1\dots \chi_n \neq 1$, on a $J_0(\chi_1,\dots,\chi_n) = 0$ d'après \cite[8.5.1]{irelandrosen}, tandis que si $\chi_1 \dots \chi_{n} = 1$, et au moins l'un des $\chi_i$, disons $\chi_n$, est non trivial,
on a $|J_0(\chi_1,\dots,\chi_n)| = (\ell-1) |J_1(\chi_1,\dots,\chi_{n-1})| = (\ell-1) \ell^{(n-3)}/2$. Comme on a $(\ell-1)^{n-1}$ caractères $\chi$ tels que $\chi_1 \dots \chi_n=1$, et que le cas particulier où tous les $\chi_i$ sont triviaux est négligeable, on obtient
$\lambda_{T(\F_\ell)}(D_0) \sim \ell^{(n-2)/2}$.
\end{pf}

\begin{remarque} On a $|T(\F_\ell)^{\ta,\dr}| \sim \ell^{n-1}$ donc la majoration de Cauchy-Schwarz donne 
$\lambda(T(\F_\ell)^{\ta,\dr}) =O(\ell^{(n-1)/2})$. On voit donc que cette majoration  donne le bon ordre de grandeur pour $a \neq 0$. Dans le cas $a=0$, elle ne le donne pas tout à fait, mais si l'on remarque que $Z T(\F_\ell)^{\tz,\dr}=T(\F_\ell)^{\tz,\dr}$ (où $Z$ est le centre de $\Gl_n(\F_\ell)$), on obtient bien en utilisant la proposition~\ref{compquot} avant la majoration de Cauchy-Schwarz que $\lambda_{T(\F_\ell)}(D_0) =\bigo{ \ell^{(n-2)/2}}$. Autrement dit, la majoration de Cauchy-Schwarz, tenant compte de l'amélioration triviale ci-dessus dans le cas $a=0$, donne un résultat optimal dans la situation du théorème~\ref{torediagonal}. 
\end{remarque}

\begin{theoreme} \label{toresymplectique} 
Soit $T$ le tore "symplectique" des matrices diagonales $\diag(x,y,y^{-1},x^{-1})$ de $\Gl_4$ sur $\spec \Z$.
Alors il existe une constante absolue $\ccl{ctsymp}$ telle que $1 \leq \lambda(T(\F_\ell)^{\tz,\nr} )< \ccr{ctsymp}$. 
Le même résultat est vrai si $T$ est remplacé par le tore des "similitudes symplectiques",
i.e. des matrices diagonales de $\Gl_4$ de la forme $\diag(zx,zy,zy^{-1},zx^{-1})$.
\end{theoreme} 
\begin{pf}
Le cas des "similitudes symplectiques" se ramène au cas "symplectique" par la prop.~\ref{compquot}. Traitons donc ce cas.
Par la proposition~\ref{propmin}, on a $\lambda(T(\F_\ell)^{\tz,\nr}) \geq 1$; il suffit donc de prouver la majoration.

Soit $D_0$ l'ensemble des matrices de $T(\F_\ell)$ dont la trace est nulle. On a $D_0 =  T(\F_\ell)^{\tz,\nr} \coprod D_0^{\nr}$,
où $D_0^{\nr}$ est l'ensemble des matrices de $D_0$ non régulières. On voit facilement que $|D_0^{\nr}| \leq \cc $ indépendamment de $\ell$.
Il suffit donc de prouver le résultat pour $\lambda(D_0)$.
 
Or il se trouve que $D_0 = D^+_0+D^-_0$ avec $D^\pm_0=\{(x,y,y^{-1},x^{-1}) \in D_0, y= - x^{\pm 1}  \}$, et 
$D^\pm_0$ est $(1,-1) H^\pm$ où $H^\pm$ est le sous-groupe de $T(\F_\ell)$ défini par $x=y^{\pm 1}$ (je dois cette observation à Felipe Voloch).

 On a donc $\lambda(D_0) \leq  \lambda(D^+_0) + \lambda(D^-_0) + \lambda(D^+_0 \cap D^-_0)$, et
 $\lambda(D^\pm_0) =1$ par le théorème~\ref{thmin}, tandis que $D^+_0 \cap D^-_0$ est de cardinal $2$. On a donc $\lambda(D_0) \leq 2+\sqrt{2}$,
 ce qui termine la preuve.
\end{pf}
 \begin{remarque} Dans ce théorème, la majoration qu'on  obtient est meilleure que celle donnée par Cauchy-Schwarz même en tenant compte de l'amélioration qui consiste à passer au quotient par $Z$. Par exemple, si $T$ est le tore des similitudes symplectiques, l'ordre de grandeur de $|T(\F_\ell)^{\tz,\nr}|$ est $\ell^2$
 et celui de $|T(\F_\ell)^{\tz,\nr}/Z|$ est $\ell$, si bien que la majoration de Cauchy-Schwarz pour ce quotient, combinée avec la proposition~\ref{compquot}
 donne $\lambda(T(\F_\ell)^{\tz,\nr}) =O(\ell^{1/2})$, moins bon que le $O(1)$ du théorème. Il serait intéressant de savoir si un phénomène semblable se produit pour $a \neq 0$, ou pour les tores symplectiques en dimension $2g > 4$. 
 \end{remarque}

\subsection{Exemples de calculs de complexités de Littlewood dans le cas non-abélien}

\subsubsection{Sous-ensemble des permutations transitives de $S_n$}

\begin{lemme} Soit $\sigma$ un $n$-cycle de $S_n$, $\lambda$ un diagramme de Young à $n$ cases, et $\rho_\lambda$ la représentation
de $S_n$ attachée à $\lambda$. Alors $ \tr \rho_\lambda(\sigma) = 0$, sauf si $\lambda$ est un diagramme à $l$ lignes, le première
ligne ayant $n-l+1$ cases et les $l-1$ autres une case chacune. Dans ce cas, $ \tr \rho_\lambda(\sigma) = (-1)^{l-1}$, et $\tr \rho_\lambda(1)={n-1 \choose l-1}$.
\end{lemme}
\begin{pf} C'est une simple application de la règle de Murnaghan-Nakayama (cf. e.g. \cite{james}). Si  $\sigma \in S_n$ est une permutation,
donc la décomposition en produit de cycles à support disjoints contient des cycles de longueur $\mu_1, \mu_2, \dots, \mu_r$, avec $\mu_1 + \dots + \mu_r =n$ et $\mu_1 \geq \mu_2 \geq \dots \geq \mu_r$, cette règle nous demande de considérer les {\it remplissages du diagramme $\lambda$ par $\mu$}, i.e. par
 les entiers de $1$ à $r$, respectivement répétés $\mu_1, \dots, \mu_r$ fois, qui satisfont les propriétés suivantes:
\begin{itemize}
\item chaque ligne et chaque colonne du diagramme remplie est croissante;
\item pour $i=1,\dots,r$, l'ensemble des cases remplies par $i$ est connexe
\item pour $i=1,\dots,r$ l'ensemble des cases remplies par $i$ ne contient pas de carré $2 \times 2$.
\end{itemize}
Dans le cas où $\sigma$ est un $n$-cycle, on a $r=1$, $\mu_1=n$, et il s'agit de remplir toutes les cases du diagramme $\lambda$ avec des $1$, ce qui évidemment n'est possible que d'une seule manière. Ce remplissage satisfait automatiquement les deux premières conditions requises; il satisfait la troisième si et seulement si le diagramme $\lambda$ lui-même ne contient pas de carré $2 \times 2$, i.e. si et seulement si la première ligne de $\lambda$ à $n-l+1$ cases et  les $l-1$ autres lignes une seule case chacune. 

La règle de Murnaghan-Nakayama attribue à chaque remplissage un poids $\pm 1$, et stipule que $\tr \rho_\lambda(\sigma)$ est la somme de ses poids
pour tous les remplissages de $\lambda$ par $\mu$. Dans le cas du diagramme à $l$ lignes considéré ci-dessus, ce poids est $(-1)^{l-1}$.

On applique la même règle pour $\sigma=1$, qui correspond à $r=n$ et $\mu_1=\dots=\mu_n=1$. Il s'agit donc de remplir $\lambda$ avec tous les nombres de $1$ à $n$, de telle manière que la première ligne, et la première colonne soient croissantes (les deux autres conditions étant automatique). Il est clair que pour un tel remplissage la première case sera toujours marquée $1$, et que le remplissage est complètement déterminé par le choix des $l-1$ éléments de $2$ à $n$ qui finissent de remplir la première colonne. De plus un tel remplissage a un poids $1$. On a donc $\tr \rho_\lambda(1) = {n-1 \choose l-1}$.
 
\end{pf}

\begin{prop} \label{phisym1} Soit $n$ un entier et $D \subset S_n$ l'ensemble des $n$-cycles (i.e. des permutations transitives). 
Alors $\lambda_{S_n}(D) = \frac{2^{n-1}}{n} $.\end{prop}
\begin{pf} On calcule $\widehat{\un_D}(\rho_{\lambda}) = \frac{1}{n!} \sum_{\sigma \in D} \tr \rho_\lambda(\sigma)$.
D'après le lemme précédent, les seuls diagrammes de Young $\lambda$ pour lesquels on obtient un résultat non nul sont ceux dont toutes les lignes
sauf la première n'ont qu'une case, et l'on a alors  $\widehat{\un_D}(\rho_{\lambda})=\frac{(-1)^{l-1}}{n}$, où $l$ est le nombre de lignes de $\lambda$.
On a donc    $\lambda_{S_n}(D) = \sum_{l=1}^n \frac{1}{n} {n-1 \choose l-1} = 2^{n-1}{n}$.
\end{pf}

\begin{remarque} Notons que bien qu'exponentiel en $n$, ce résultat est bien meilleur que la borne de Cauchy-Schwartz $\lambda(D) \leq \sqrt{(n-1)!}$ qui elle est super-exponentielle.
\end{remarque}

\subsubsection{Groupes munis d'un système de Tits}

Dans ce numéro, nous supposons que
\begin{itemize} \item[(i)]
$G$ est un groupe fini muni d'un système de Tits $(B,N)$ (cf. \cite[ch. IV, \S2]{bourbaki}).
\end{itemize}
On pose comme d'habitude $T=B \cap N$, et $W = N/T$: c'est le groupe de Weyl du système de Tits. Pour $w$ in $W$, $t \in T$, on écrit $t^w$ pour $ntn^{-1}$
si $n$ est un représentant de $w$ dans $N$. Nous supposerons de plus donné 
\begin{itemize} \item[(ii)]
un sous-groupe normal $U$ de $B$ qui contient le groupe dérivé de $B$ et tel que $B=TU$ et $U \cap T =\{1\}$. 
\end{itemize}
Ceci implique en particulier que $T$ est abélien. Notons $s$ la surjection  $s : B \rightarrow B/U \simeq T$, où l'isomorphisme $B/U \simeq T$ est l'inverse de l'isomorphisme $T \simeq B/U$ induit par l'injection canonique de $T$ dans $B$. On a $\ker s = U$ et $s_{|T} = \Id_T$.
Finalement, nous supposerons que
\begin{itemize} \item[(iii)] si $n \in N$, $nUn^{-1} \cap B \subset U$.
\end{itemize} 

\begin{exemple} \label{exembleredbn} Soit $G$ un groupe réductif connexe sur $\F_q$. Alors, $G$ est quasi-déployé (théorème de Lang) et si $T$ est le centralisateur
d'un tore maximal déployé, $B$ un Borel contenant $T$, et $U$ le radical unipotent de $B$, alors $G$ muni de $(B,N)$ et de $U$ 
est un système de Tits satisfaisant les conditions requises. Voir \cite{carter}.
\end{exemple}

\begin{lemme} \label{bbu} Soit $d$ un élément de $T$, dont le centralisateur dans $G$ est $T$. Soit $u$ un élément de $U$. Alors $ud$ et $du$ sont conjugués à  $d$ dans $B$.
\end{lemme} 
\begin{pf} (Je remercie Will Sawin pour cet argument). Considérons l'application $f: B \rightarrow B$, $b \mapsto  b^{-1} d b d^{-1}$. Les fibres de cette application sont des classes à gauche sous le centralisateur de $d$, c'est-à-dire $T$, donc l'image de $B$  est de cardinal $|B|/|T|$. Comme l'image de $f$ est contenue dans $U$ par l'hypothèse (ii), et que $|U|=|B|/|T|$, $f$ est surjectif sur $U$. Il existe donc $b \in B$ tel que $f(b)=u$, c'est-à-dire tel que $b^{-1} d b = ud$. Le cas de $du$ se traite de même.
\end{pf}
\begin{lemme} \label{bsb} Si $b \in B$, $b$ est conjugué à $s(b)$ dans $B$. Si $n \in N$, $w$ l'image de $n$ dans $W$, $b$ in $B$ et si 
$nbn^{-1} \in B$, alors $nbn^{-1}$ est conjugué à $s(b)^w$ dans $B$.
\end{lemme}
\begin{pf} Si l'on pose $u=b^{-1} s(b)$, on a $s(u) = s(b^{-1}) s(s(b)) = s(b^{-1}) s(b)=1$ donc $u \in U$, et donc $s(b)=bu$ est conjugué à $b$ dans $B$ par le lemme précédent. Pour la deuxième assertion, écrivant encore $b = s(b)u^{-1}$, on a $nbn^{-1} = s(b)^w nun^{-1}$. Par hypothèse, $nbn^{-1} \in B$, donc $nun^{-1} \in B$, et par l'hypothèse (iii), $nun^{-1} \in U$, si bien que $s(nbn^{-1})= s(b)^w$ le résultat en découle.
\end{pf}

 \begin{theoreme}\label{tits} Soit  $D$ un sous-ensemble de $G$ stable par conjugaison. Supposons que pour tout $d \in D$, le centralisateur de $d$ dans $G$ soit conjugué
 à $T$. Alors on a $\lambda_G(D) = \frac{|G|}{|W| |B|} \lambda_T(D \cap T)$.
 \end{theoreme} 
 \begin{pf} On va appliquer le cas d'égalité du 
 théorème~\ref{methserre} avec $G=G$, $H=B$, $U=U$.
 Vérifions-en les quatre hypothèses:
 \begin{itemize}
 \item[(a)] Soit $d \in D$. Par hypothèse, $d$ est conjugué à un élément $d'$ de $T$, et comme $C(d)=C(d')$,
  on est ramené au cas $d \in T \cap D$.
 Comme le centralisateur de $d$ est $T$, $$|C(d)|=|G|/|T|,$$ qui est bien indépendant de $d$.
 \item[(b)] Soit $d \in D$, que nous pouvons à nouveau supposer dans $T$. On a \begin{eqnarray} \label{CdB} C(d) \cap B = \coprod_{w \in W} C_B(d^w),\end{eqnarray} où  $C_B$ est la classe de conjugaison de $d^w$ dans $B$. En effet, soit $d' \in C(d)$
et  choisissons $g \in G$ tel que $d'=g d g^{-1}$. 
Le choix de l'élément $g$ n'est pas unique, mais la classe $gT$ l'est. Soit $B w B$, pour $w \in W$, la cellule, évidemment bien déterminée, de la décomposition de Bruhat à laquelle $g$ appartient. Écrivons $g=b_1 n b_2$ avec $b_1,b_2 \in B$, et $n$ un représentant de $w$ dans $N$. On a alors
$d'=gdg^{-1}=b_1 n b_2 d b_2^{-1} n^{-1} b_1^{-1}$ et l'on voit que $d'$ est dans $B$ si et seulement si $n b_2 d b_2^{-1} n^{-1}$ l'est, auquel cas cet élément, et donc aussi $d'$, est conjugué dans $B$ à $d^w$ par le lemme~\ref{bsb}, ce qui prouve (\ref{CdB}).
Comme le normalisateur dans $B$ d'un élément régulier de $T$ est toujours $T$, $|C_B(d^w)| = |B|/|T|$ et il en résulte que $$|C(d) \cap B| = |W| |B| / |T|,$$
qui est bien indépendant de $B$.
\item[(c)] Il faut voir que $U(D \cap B) = D \cap B$, ce qui suit immédiatement du lemme~\ref{bbu} puisque $D$ est invariant pas conjugaison.
\item[(d)] Soit $\pi \in \Gd$. Il faut montrer que sur les caractères de $T$ qui apparaissent dans $\pi_{|B}$, la fonction $\widehat{\un_{s(D \cap B)}}$ a
un argument constant. Or ces caractères forment une classe de conjugaison de $\Td$ sous $W$, par le critère de Mackey.
Il suffit donc de voir que la fonction $\widehat{\un_{s(D \cap B)}}$ sur $T$ est invariante par $W$, donc de voir que le sous-ensemble $s(D \cap B)$ de $T$
est invariant par $W$, ce qui est encore une conséquence du lemme~\ref{bsb}.
\end{itemize}
Le théorème~\ref{methserre} donne donc
$$\lambda_G(D) = \lambda_T(s(D \cap B))$$
et comme $s(D \cap B)=D \cap T$ par le lemme~\ref{bsb}, le théorème est prouvé.
 \end{pf}

Nous appliquons maintenant le théorème aux sous-groupes réductifs connexes de $\Gl_n$:
soit $N\geq 1$ un entier, $G$ un sous-schéma en groupes réductif connexe\footnote{Dans la terminologie de \cite{sga3}, tous les schémas en groupes réductifs sont connexes, mais nous n'adopterons pas cette terminologie.} de $\Gl_n$ sur $\Z[1/N]$, $T$ un tore maximal de $G$. On écrit $G_\Q$ (resp. $G_{\F_\ell}$)
pour les fibres de $G$ au dessus du point générique (resp. du point spécial de caractéristique $\ell$) de $\spec \Z[1/N]$, $G_{\bar \Q}$ pour 
$G_\Q \times_{\spec \Q} \spec \bar \Q$ et de même pour $T_\Q$, $T_{\F_\ell}$, $T_{\bar Q}$. Soit $W$ le groupe de Weyl de $T_{\bar \Q}$ dans $G_{\bar \Q}$.
On note $d$ la dimension commune des $G_\Q$, $G_{\F_\ell}$ pour $\ell \nmid N$, $r$ (le rang réductif) celle des $T_\Q$, $T_{\F_\ell}$.
Il suit facilement de \cite[Théorème 2.5, exposé 19]{sga3} qu'il existe un ensemble frobénien $\Lambda$ de densité $>0$ de $\ell$ ne divisant pas $N$ tel que, pour $\ell \in \Lambda$, $T_{\F_\ell}$ soit un tore maximal déployé de $G_{\F_\ell}$ et le groupe de Weyl $W_\ell$ de $T_{\F_\ell}$ dans $G_{F_\ell}$ s'identifie canoniquement avec le groupe de Weyl $W$.
 
\begin{cor} \label{corGT} Gardons les notations du paragraphe précédent.
 Fixons un entier $a \in \Z$.  Pour $\ell$ un nombre premier, on note $G(\F_\ell)^{\ta,\dr}$ (resp. $T(\F_\ell)^{\ta,\dr}$) l'ensemble des éléments de
 $G(\F_\ell)$ (resp. $T(\F_\ell)$) qui, dans $\Gl_n(\F_\ell)$ sont de trace $a \pmod{\ell}$ et diagonalisables réguliers.
 Alors, quand $\ell \in \Lambda$, $\ell \rightarrow \infty$, on a
 $$\lambda_{G(\F_\ell)} (G(\F_\ell)^{\ta,\dr}) \sim \frac{1}{|W|} \ell^{(d-r)/2} \lambda_{T(\F_\ell)} (T(\F_\ell)^{\ta,\dr}) $$
 \end{cor}
 \begin{pf} Pour $\ell \in \Lambda$, $G(\F_\ell)$ possède un système de Tits $(B,N)$, avec $B=B(\F_\ell)$ les points d'un Borel contenant $T(\F_\ell)$, $N$ le normalisateur de $T$, $T=T(\F_\ell)$, et de groupe de Weyl $W_\ell=W$, satisfaisant 
 les hypothèses du théorème précédent. Comme il est clair que $G(\F_\ell)^{\ta,\dr})  \cap T(\F_\ell) = T(\F_\ell)^{\ta,\dr})$, et que 
 $|B(\F_\ell)| \sim \ell^{(n+r)/2}$, $T(\F_\ell) = (\ell-1)^r \sim \ell^r$ ce théorème donne l'équivalence voulue.
 \end{pf}
 
 Donnons deux cas particuliers:
 \begin{cor} \label{corGln} Soit $G=\Gl_n$ sur $\spec \Z$, $T$ le tore diagonal, $a \in \Z$. Alors quand $\ell$ tend vers l'infini, 
 $$\lambda_{G(\F_\ell)} (G(\F_\ell)^{\ta,\dr}) \sim \begin{cases} \frac{1}{n!} \ell^{(n^2-1)/2} & \text{ si } a \neq 0 \\
 \frac{1}{n!} \ell^{(n^2-2)/2}  & \text{ si } a = 0\end{cases} $$
 \end{cor}
 \begin{pf} Cela résulte du corollaire précédent et du théorème~\ref{torediagonal}.\end{pf}

\begin{cor} \label{corGsp4} Soit $G=\GSP_4 \subset \Gl_4$, $T$ le tore maximal des matrices diagonales de $G$. Alors
$$\lambda_{G(\F_\ell)} (G(\F_\ell)^{\tz,\dr}) = O (\ell^4)$$ 
\end{cor}
\begin{pf} La dimension de $\GSP_4$ est $d=11$, celle de son tore maximal est $r=3$, et on a $ \lambda_{T(\F_\ell)} (T(\F_\ell)^{\tz,\dr}) \leq \cc$ d'après le théorème~\ref{toresymplectique}.
\end{pf}

\subsubsection{Ensemble des matrices de trace donnée dans $\Gl_n$}

Fixons un entier $n \geq 1$. Pour $\ell$ un nombre premier, soit $\F_\ell$ le corps fini de cardinal $\ell$, $a \in \Z$, $G = \Gl_n(\F_\ell)$, 
$D_a$ l'ensemble des matrices de $\Gl_n(\F_\ell)$ dont la trace est $a \pmod{\ell}$. On se donne pour objectif de calculer l'ordre de grandeur, quand $\ell$
tend vers l'infini, de $\lambda_G(D_a)$. Nous allons voir que cet ordre de grandeur est le même que celui de $\lambda(\Gl_n(\F_\ell)^{\ta,\dr})$ que nous avons calculé ci-dessous. Autrement dit, pour $\Gl_n$, il revient au même, en terme d'ordres de grandeur du moins, de travailler avec toutes les matrices de trace $a$,
ou bien seulement celles qui sont régulières et diagonalisables.

Commençons par majorer $\lambda_G(D_a)$. On a $|D_a| \leq \ell^{n^2-1}$ puisque $\ell^{n^2-1}$ 
est le nombre de matrices de trace $a$ dans $M_n(\F_\ell)$, d'où $\lambda_G(D_a) \leq \ell^{\frac{n^2-1}{2}}$ par la majoration de 
Cauchy-Schwarz. Si $a \equiv 0 \pmod{q}$, soit $Z$ le centre de $\Gl_n(\F_\ell)$. On a $ZD_0 \subset D_0$, d'où $\lambda_G(D_0) = \lambda_{G/Z}(D_0/Z) \leq \sqrt{|D_0|/|Z|} \leq \sqrt{ \ell^{n^2-1}/(\ell-1)} \leq 2 \ell^{\frac{n^2-2}{2}}$ dès que $\ell \geq 2$.  Ces deux majorations sont implicitement appliquées dans \cite{effective}, au moins pour $n=2$. Nous allons voir qu'elles donnent le bon ordre de grandeur pour $\lambda_G(D_a)$ et $\lambda_G(D_0)$.

\begin{theoreme} \label{thmlambdagln} Pour $n$ et $a$ fixé, on a pour $\ell>a$: 
\begin{eqnarray*}\lambda_G(D_a) & \asymp &  \ell^{\frac{n^2-1}{2}} \text{ si }a \neq 0 \\
\lambda_G(D_0) & \asymp & \ell^{\frac{n^2-2}{2}} 
 \end{eqnarray*}
\end{theoreme}

%\begin{eqnarray*} \frac{1}{n!+1} q^{\frac{n^2-1}{2}} & \leq & \lambda_G(D_a) \leq q^{\frac{n^2-1}{2}} \text{ si } a \neq 0 \\
% \frac{1}{n!+1} q^{\frac{n^2-2}{2}} \leq \lambda(D_a) & \leq & 2 q^{\frac{n^2-2}{2}} \text{ si }a \neq 0 
%\end{eqnarray*}
%
\begin{pf}
 Pour minorer asymptotiquement $\lambda_G(D_a)$ écrivons $D_a = D_a^\reg \coprod D_a^\nr$, où $D_a^\reg$ est l'ensemble des éléments réguliers de $D_a$ (i.e. ceux dont les
 valeurs propres dans $\bar \F_q$ sont deux-à-deux distinctes) et $D_a^\nr$ l'ensemble de ceux qui ne le sont pas. On a
 \begin{eqnarray*} \lambda_G(D_a^\nr) &=& O(\ell^{\frac{n^2-2}{2}})  \text{ si }a \neq 0 \\
 \lambda_G(D_a^\nr) &=& O(\ell^{\frac{n^2-3}{2}}) \text{ si } a = 0.
\end{eqnarray*}
 En effet, $|D_a^\nr| = O(\ell^{n^2-2})$ d'où le résultat pour $a \neq 0$ par la majoration de Cauchy-Schwartz. Et quand $a=0$, $D_a^\nr$ est encore invariant par $Z$, et $|D_a^\nr|/|Z| = O(\ell^{n^2-3})$ et le résultat en découle par  le corollaire~\ref{compquotcor} et la majoration de Cauchy-Schwartz. 
Nous sommes donc ramenés à démontrer $\lambda_G(D_a^\reg) \gg \ell^{\frac{n^2-1}{2}}$ pour $a \neq 0$, et $\lambda_G(D_0^\reg) \gg \ell^{\frac{n^2-2}{2}}$.

Pour cela, rappelons trois résultats sur les séries principales. Soit $T$ le tore diagonal de $\Gl_n(\F_\ell)$, $B$ le Borel des matrices triangulaires supérieures, $W = S_n$ le groupe de Weyl. Pour $\chi : T \rightarrow \Cb^\ast$ un caractère, on notera $I(\chi)$ la représentation $\Ind_B^G \chi$. 
\begin{itemize}
\item[(a)] la représentation $I(\chi)$ est irréductible si et seulement si le caractère $\chi$ est {\it régulier}, i.e. $\chi^w \neq \chi$ pour tout $w \in W$.
\item[(b)] Si $g \in G^\reg$, $\chi$ un caractère quelconque de $G$, et si $\tr I(\chi)(g) \neq 0$, alors $g$ est diagonalisable (i.e. conjugué à un élément de $T$).
\item[(c)] Soit $g \in G^\reg$ est diagonalisable et régulier, $\pi \in \Gd$; si $\tr \pi(g) \neq 0$ alors $\pi \subset I(\chi)$ pour un caractère $\chi$ de $T$.
\end{itemize}
Les deux premiers résultats sont bien connus et élémentaires: le premier se déduit par exemple du critère d'irréductibilité de  Mackey (cf. \cite[prop. 23]{serrerep}), et le deuxième encore plus facilement de la formule du caractère d'une induite ({\it loc. cit.}, prop. 20).  Le troisième semble plus difficile: il 
résulte de la théorie de Deligne-Lusztig (cf. \cite[formula (7.6.2)]{DL}).

Écrivons alors \begin{eqnarray*} \lambda_G(D_a^\reg) &=& \sum_{\pi \in \Gd} |\widehat{\un_{D_a^\reg}}(\pi)| \dim \pi  \\
&\geq& \sum_{\chi} \sum_{\pi \subset I(\chi)} |\widehat{\un_{D_a^\reg}}(\pi)| \dim \pi \\
& = & \sum_{\chi\text{ régulier }} |\widehat{\un_{D_a^\reg}}(I(\chi))| \dim I(\chi) +  \sum_{\chi\text{ non régulier }} |\widehat{\un_{D_a^\reg}}(I(\chi))| \dim I(\chi) 
\end{eqnarray*}
On peut borner le second terme par Cauchy-Schwarz
$$   \sum_{\chi\text{ non régulier }} |\widehat{\un_{D_a^\reg}}(I(\chi))| \dim I(\chi) \leq  \sqrt{\frac{|D_a^\reg|}{|G|} \sum_{\chi\text{ non régulier}} (\dim I(\chi)^2)}.$$ On a $|G| \sim \ell^{n^2}$, $|D_a^\reg| \sim \ell^{n^2-1}$, $\dim I(\chi)^2 = |G/B|^2 \sim \ell^{n(n-1)}$, et le nombre de $\chi$ non-réguliers
est $O(\ell^{n-1})$, d'où l'on voit que la contribution des $\chi$ non réguliers est $O( \sqrt{ \ell^{-1} \ell^{n(n-1)} \ell^{n-1}}) = O(\ell^{(n^2-2)/2})$. 
Quand $a = 0$, on peut travailler dans $G/Z$ et montrer de le même façon que cette contribution est $ O(\ell^{(n^2-3)/2})$. Dans tous les cas, elle est négligeable au vu du résultat à prouver, et on est donc ramené à prouver que 
\begin{eqnarray*} \sum_{\chi\text{ régulier }} |\widehat{\un_{D_a^\reg}}(I(\chi))| \dim I(\chi) &\gg& \ell^{\frac{n^2-1}{2}} \text{ si } a \neq 0, \\ &\gg&  \ell^{\frac{n^2-1}{2}} \text{ sinon }\end{eqnarray*}
Cette somme est égale à  $ \sum_{\chi\text{ régulier }} |\widehat{\un_{ \Gl_n(\F_\ell)^{\ta,\dr} }}(I(\chi))| \dim I(\chi) $ par le point (b) rappelé ci-dessus.
Par le même calcul utilisant Cauchy-Schwartz que plus haut, on voit  que   $$ \sum_{\chi\text{ non-régulier }} |\widehat{\un_{\Gl_n(\F_\ell)^{\ta,\dr}}}(I(\chi))| \dim I(\chi) $$ est négligeable,
et on est donc ramené à prouver la minoration du théorème pour $$ \sum_{\chi} |\widehat{\un_{\Gl_n(\F_\ell)^{\ta,\dr}}}(I(\chi))| \dim I(\chi).$$ Comme cette somme
est égale à   $\lambda_G(\Gl_n(\F_\ell)^{\ta,\dr})$ par le point (c), il faut prouver que 
\begin{eqnarray*}  \lambda_G(\Gl_n(\F_\ell)^{\ta,\dr}), &\gg& \ell^{\frac{n^2-1}{2}} \text{ si } a \neq 0, \\ &\gg&  \ell^{\frac{n^2-1}{2}} \text{ sinon },\end{eqnarray*}
ce qu'on a déjà fait plus haut, au corollaire~\ref{corGln}.

\end{pf}

\section{L'invariant $\varphi_G(D)$}

Cette section est destinée à l'étude de l'invariant $\varphi(D)=\varphi_G(D)$ d'un sous-ensemble $D$ non-vide invariant par conjugaison d'un groupe fini $G$. Reppelons-en la définition:

$$\varphi_G(D) = \inf_f \frac{\lambda_G(f)}{\mu_G(f)},$$ quand $f$ parcourt l'ensemble des fonctions sur $G$ à valeurs réelles,
qui satisfont les deux conditions suivantes
\begin{itemize} \item[(a)] si $f(g)>0$, alors $g \in D$.
\item[(b)] $\mu_G(f) >  0$.
\end{itemize}

Cet invariant s'avère beaucoup plus subtil que la norme de Littlewood $\lambda(D)$, et pour l'instant il y a bien peu de cas où nous pouvons le calculer exactement, et à peine plus où nous savons en donner une estimation intéressante.   

\subsection{Propriétés élémentaires}

\begin{prop} \label{extf} La borne inférieure définissant $\varphi(D)$ est atteinte. Plue précisément, il existe une fonction $f_D$ centrale sur $G$, et une seule à multiplication par un réel strictement positif près, tel que:
\begin{itemize} 
 \item[(a)] si $f_D(g)>0$, alors $g \in D$;
\item[(b)] $\mu_G(f_D) > 0$;
\item[(c)] $\varphi_G(D)=\frac{\lambda_G(f_D)}{\mu_G(f_D)}$.
\end{itemize} 
\end{prop}
\begin{pf} Soit $E$ l'ensemble des fonctions $f : G \rightarrow \R$ centrales, satisfaisant (a), (b) et la condition $\sum_g |f(g)| = 1$. 
Comme toute fonction satisfaisant (a) et (b)
est, à multiplication par un réel strictement positif unique près, dans $E$, on voit qu'il suffit de montrer que la fonctionnelle $f \mapsto \frac{\lambda_G(f)}{\mu_G(f)}$
admet un unique minimum sur $E$. 
Or $E$ est convexe, et cette fonctionnelle est strictement convexe sur $E$, ce qui prouve que la borne inférieure, si elle est atteinte, l'est par une unique fonction $f \in E$.

Si $\epsilon$ est un réel positif, soit $E_\epsilon$ le sous-ensemble de $E$ défini par la
condition $\mu_G(f) \geq \epsilon$. L'ensemble $E_\epsilon$ est compact, étant fermé dans la sphère unité pour la norme $L^1$ de l'espace des fonctions de $G$ dans $\Cb$. Pour $f \in E-E_\epsilon$ on a $\frac{\lambda_G(f)}{\mu_G(f)}> \frac{\lambda_G(f)}{\epsilon}  \geq \frac {\NN{f}_\infty}{\epsilon} \geq \frac{1}{|G|\epsilon}$ (en utilisant la proposition~\ref{propmin} puis le fait que $\NN{f}_\infty \geq \NN{f}_1 / |G|$). On voit donc que 
la borne inférieure de $f \mapsto \frac{\lambda_G(f)}{\mu_G(f)}$ sur $E$ est égal à la borne intérieure sur $E_\epsilon$ pour $\epsilon$ suffisamment petit, qui est atteinte par compacité.
\end{pf}

\begin{cor} On a $\varphi_G(D)>0$.
\end{cor}
\begin{pf} En effet, la fonction $f_D$ est non nulle par (a), donc $\lambda_G(f_D) >0$ puisque $\lambda_G$ est une norme, et donc $\varphi_G(D)>0$ par (c).
\end{pf}

\begin{cor} Soit $\sigma : G \rightarrow G$ une bijection qui est soit  un automorphsime de groupe, soit un anti-automorphisme, 
soit de la forme $g \mapsto zg$ où $z$ est un élément du centre. Alors $\varphi(\sigma(D))=\varphi(D)$. Si de plus $\sigma(D)=D$, alors la fonction $f_D$ de la proposition précédente satisfait $f_D \circ \sigma = f_D$.
\end{cor}  
\begin{pf} D'après la proposition~\ref{invaut}, pour $\sigma$ comme dans l'énoncé on a $\lambda_G(f)=\lambda_G(f \circ \sigma)$ pour toute fonction centrale $f$ sur $G$, et comme évidemment $\mu_G(f) = \mu_G(f \circ \sigma)$, la première assertion en découle aisément.
 Pour la seconde, remarquons que si $\sigma(D)=D$, l'espace des fonctions $f$ satisfaisant (a) et (b) est invariant par $\sigma$,  ainsi comme on vient de le voir que la fonctionnelle $f \mapsto \lambda_G(f)/\mu_G(f)$,
donc $f_D \circ \psi=\nu f_D$ pour un réel $\nu$ par la proposition précédente, et évidemment $\nu=1$.
\end{pf}

\begin{cor} Soit $G$ un groupe abélien, $H$ un sous-groupe, $A \subset G/H$. Alors $\varphi_G(AH)=\varphi_{G/H}(A)$.
\end{cor}
\begin{pf}
On applique le corollaire précédent aux bijections $\psi_z: x \mapsto zx$ de $G$ pour $z \in H$.
Comme ces bijections stabilisent $AH$, on voit que la fonction $f_{AH}$ qui réalise le minimum est invariante par translation par $H$, et provient donc d'une
fonction $\tilde f$ sur $G/H$ vérifiant les conditions (a) et (b) relative au sous-ensemble $A$ de $G/H$. On a $\lambda_G(f_{AH})=\lambda_{G/H}(\tilde f)$ par la proposition~\ref{compquot}, d'où $\varphi_{G/H} (A) \leq \varphi_G(AH)$. L'autre inégalité est claire.
\end{pf} 
En particulier, pour l'anti-automorphisme $x \mapsto x^{-1}$, on obtient:
\begin{cor} Soit $G$ un groupe fini, $D$ un sous-ensemble invariant par conjugaison et par l'anti-automorphisme $g \mapsto g^{-1}$. 
Alors $\varphi_G(D)$ est la borne inférieure de $\lambda_G(f)/\mu_G(f)$ quand $f$ parcourt l'ensemble des fonctions centrales  à valeurs réelles telles que
\begin{itemize} 
\item[(a)] si $f(g)>0$, alors $g \in D$;
\item[(b)] $\mu_G(f) > 1$;
\item[(c)] Pour tout $g$ dans $G$, $f(g)=f(g^{-1}).$
\end{itemize}
\end{cor}
Plaçons-nous sous les hypothèses du corollaire précédent. Pour $f$ satisfaisant (c), la fonction $\widehat{f}$ est à valeurs réelles, si bien que $\lambda_G(f)$ est {\bf combinaison linéaire} à coefficients entiers de {\bf valeurs absolues de forme linéaires} $f \mapsto \widehat{f}(\rho)$ sur l'espace des fonctions $f$ centrales satisfaisant (c).
Autrement dit, calculer le minimum $\varphi_G(D)$ et trouver la fonction $f_D$ qui le réalise est dans ce cas {\bf un problème de programmation linéaire}.

\begin{prop} Soit $D \subset G$ un sous-ensemble non vide invariant par conjugaison. On a  $$\varphi(D) \leq \frac{\lambda(D)|G|}{|D|}$$ 
avec égalité si et seulement si $f_D = \un_D$ à multiplication par un réel près.
\end{prop}
C'est clair.

\begin{cor} \label{corphics} Soit $D \subset G$ un sous-ensemble non vide invariant par conjugaison. On a
\begin{eqnarray} \varphi(D) &\leq& \frac{|G|}{\sqrt{|D|}} 
\end{eqnarray}
\end{cor}
\begin{pf} Cela résulte de la proposition précédente vue la majoration de Cauchy-Schwarz (proposition~\ref{majcs}.)\end{pf}

\begin{prop} On a $$\varphi(D) \geq \frac{|G|}{|D|}$$ avec égalité si et seulement si $\lambda(D)=1$, i.e. si et seulement si $D$ est de la forme $aH$ avec $H$ un sous-groupe distingué dans $G$ et $a$ un élément de $G$ dont l'image dans $G/H$ est centrale. Lorsque que ces conditions sont satisfaites, $f_D=\un_D$ à un scalaire près.
\end{prop}
\begin{pf} 
Soit $f$ satisfaisant (a) et (b). 
On a $\frac{\lambda_G(f)}{\mu_G(f)} \geq \frac{\lambda_G(f)|G|}{\sum_{g \in G, f(g)>0} f(g)} \geq \frac{\NN{f}_\infty |G|}{\sum_{g \in G, f(g)>0} f(g)}$
d'après la proposition~\ref{propmin}. Comme d'après (b), la somme du dénominateur a au plus $|D|$ termes, et que chacun d'eux est positif et plus petit que $\NN{f}_\infty$, on voit que $$\frac{\lambda_G(f)}{\mu_G(f)} \geq \frac{|G|}{|D|}.$$ L'égalité n'est possible que s'il n'y a pas de $g$ avec $f(g)<0$ et que tous les
$g$ tels que $f(g)>0$ sont tels que $f(g)=\NN{f}_\infty$, autrement dit que si $f=\un_D$ à un scalaire près. Dans ce cas, on a $\lambda_G(D)=1$ et donc $D$ est bien de la forme annoncée d'après le théorème~\ref{thmin}.
\end{pf}  

\subsection{Exemples de calculs et d'estimations de $\varphi(D)$}

\subsubsection{Générateurs de certains groupes cycliques}

\begin{prop} \label{phigen} Soit $n$ un entier tel que $\sum_{p \mid n} 1/p < 1$. Soit $D$ l'ensemble des générateurs de $\Z/n\Z$. 
Alors $\varphi_{\Z/n\Z}(D) \leq \frac{\omega(n)+1}{1 - \sum_{p \mid n} 1/p}$ 
\end{prop}
\begin{pf} Soit $f = \un_{\Z/n\Z} -  \sum_{p \mid n} \un_{p \Z/n\Z}$. Il est clair que $f(a) \leq 0$ dès que $a$ n'est pas un générateur de $\Z/n\Z$. Par ailleurs on a $\mu_G(f) = 1 - \sum_{p \mid n} 1/p$ qui est strictement positif sous nos hypothèses. On a donc $\varphi(f) \leq \lambda(f)/\mu(f) \leq  (1 + 
\sum_{p \mid n} 1) / (1 - \sum_{p \mid n} 1/p) = (1+\omega(n))/(1 - \sum_{p \mid n} 1/p)$.
\end{pf}
On peut comparer ce résultat à celui qu'on obtient à partir du calcul de $\lambda(D)$ (prop.~\ref{lambdagen}):
 $\varphi_{\Z/n\Z}(D) \leq \frac{|G|}{|D|} \lambda(D) = 2^{\omega(n)}$ valable lui sans restriction sur $n$. On voit donc facilement qu'il y a une infinité de $n$ tels que 
 $ \varphi_{\Z/n\Z}(D)  < \frac{|G|}{|D|} \lambda(D) \lambda(d)$, i.e. telle que la fonction optimale $f_D$ n'est pas $\un_D$.
  Il serait très intéressant de connaître l'ordre de grandeur, ou même la valeur exacte de  $\varphi_{\Z/n\Z}(D)$. Comme nous l'avons vu, c'est un problème de programmation linéaire.

\subsubsection{Certains sous-ensemble des permutations de $S_n$}

Dans ce numéro, on note $f$ la fonction qui à une permutation $g$ du groupe symétrique $S_n$ associe le nombre de points que $g$ laisse fixes dans $\{1,\dots,n\}$.

\begin{lemme} Soit $G$ un sous-groupe du groupe symétrique $S_n$.
Alors \begin{eqnarray*} \lambda_G(f) &=& n \\ 
\lambda_G(f^2) &=& n^2 \\
\mu_G(f) & \geq & 1 \text{ avec égalité ssi $G$ agit transitivement sur $\{1,\dots,n\}$} \\
\mu_G(f) & \geq & 2 \text{ avec égalité ssi $G$ agit doublement transitivement sur $\{1,\dots,n\}$} 
\end{eqnarray*}
\end{lemme}
\begin{pf} La fonction $f$ est la trace de la représentation naturelle de $G$ sur $V=\Cb^n$ par permutation des coordonnées. Écrivons
cette représentation $\oplus_\rho \rho^{m_\rho}$, on a $\lambda_G(f) = \sum_{\rho} |m_\rho| \dim \rho = \sum_{\rho} m_\rho \dim \rho = \dim \Cb^n =n$.
On raisonne de même avec $f^2$ qui est la trace de la représentation $V \otimes V$.
Les résultats sur $\mu_G$ sont bien connus (cf. \cite[exercice 2.3]{serrerep}).
\end{pf}
On a  plus généralement avec la même preuve $\lambda_G(f^k) = n^k$ pour tout $k \geq 0$. En revanche on se gardera bien de croire
 que $\mu_G(f^k)=k$ si $G$ agit $k$-transitivement sur $\{1,\dots,n\}$.

\begin{prop} \label{propphisym} Soit $G$ un sous-groupe du symétrique $S_n$. On a
\begin{eqnarray*} \varphi_G( \{g \in G \, | \, f(g) \geq 1\}) & \leq & n\\ 
\varphi_G( \{g \in G \, | \, f(g) \geq 2\}) & \leq & n^2+1\\ 
\varphi_G( \{g \in G \, | \, f(g) = 0 \}) & \leq & 2n^2+2n   \text{ si $G$ agit transitivement sur $\{1,\dots,n\}$}
\end{eqnarray*}
\end{prop} 
\begin{pf} 
La fonction $f$ est par définition nulle sur le complémentaire de $\{g \in G \, | \, f(g) \geq 1\}$ et l'on a $\mu_g(f) \geq 1 >0$, donc par définition $\varphi(\{g \in G \, | \, f(g) \geq 1\}) \leq \lambda(f)/\mu(f) \leq n/1 = n$.

Considérons la fonction $f'=f^2-1$. On voit aisément que $f'$ est négative ou nulle sur le complémentaire de $\{g \in G \, | \, f(g) \geq 2\}$, et 
$\mu_G(f') = \mu_G(f^2)-1 \geq 1$. D'où $\varphi_G( \{g \in G \, | \, f(g) \geq 2\}) \leq \lambda(f')/\mu_G(f') \leq n^2+1$.

Considérons la fonction $f'' =  f^2 - (n+1) f  + n$.  Si $G$ agit transitivement sur $\{1,\dots,n\}$,
on a $\mu_G(f'') \geq 2 - (n+1) + n = 1$. Par ailleurs, on voit aisément que $f''$ est négatif ou nul sur le complémentaire de $\{g \in G \, | \, f(g) = 0 \}$, c'est-à-dire dès que $f(g) \geq 1$. Donc $\varphi(\{g \in G \, | \, f(g) = 0 \}) \leq \lambda(f'')/\mu_G(f'') = n^2 + (n+1)n + n = 2n^2+2n$.
\end{pf}
\begin{remarque} Dans les trois exemples ci-dessus, on a donné une borne polynomiale en $n$ pour $\varphi_G(D)$, alors que la majoration de Cauchy-Schwarz
$\varphi_G(D) \leq |G|/\sqrt{D}$ donnerait une borne super-exponentielle.
\end{remarque}

\section{Les théorèmes principaux}

\subsection{Le théorème de Chebotarev effectif}

\label{densitefrob}

Comme dans l'introduction, soit $L$ un corps de nombres galoisien sur $\Q$, $G=\Gal(L/\Q)$, et $M$ le produit des nombres premiers qui 
sont ramifiés dans $L$. Soit $f$ une fonction centrale à valeurs complexes sur $G$, et $\pi(f,x) = \sum_{p < x} f(\Frob_p)$.
Rappelons l'énoncé de la version du théorème de Chebotarev que nous utilisons dans cet article.
\renewcommand{\thetheoremerep}{\ref{chef}}
\begin{theoremerep}[Chebotarev effectif] Supposons vraie (GRH) et la conjecture d'Artin pour les
 fonctions $L$ d'Artin associées aux représentations irréductibles de $\Gal(L/\Q)$ qui appartiennent au support spectral de $f$. On a, pour $x \geq 3$,
\begin{eqnarray*}  \left| \pi(f,x) - \mu(f) \Li(x)  \right| < \ccr{ct} \, x^{1/2}\,  \lambda(f)\, (\log x+\log M+\log |G| ) \end{eqnarray*}
\end{theoremerep}
Montrons d'abord comment déduire ce théorème des résultats de \cite{ik} et \cite{effective}.
Nous noterons  $\Lambda_f(n)$ la fonction de Von Humbolt associée à $f$, définie par:
$$\Lambda_f(n) = \begin{cases} f(\Frob_p^k) \log p &\text{si }  n=p^k, p \nmid M, \\
0 & \text{sinon.} \end{cases}$$
On pose
\begin{eqnarray} 
\label{defpsifphi}\psi(f,x) &=& \sum_{n \leq x} \Lambda_f(n) 
\end{eqnarray}

On écrit $f = \sum_{\pi \in \Gd} c_\pi \chi_\pi$, où $c_\pi = \widehat{f}(\pi)$. On a donc par définition $\lambda_G(f) = \sum_\pi |c_\pi| \dim \pi$,
et l'on voit immédiatement que pour prouver le théorème~\ref{chef}, il suffit de le faire quand $f = \chi_\pi$ pour $\pi \in \Gd$.
Notons que dans ce cas, $\mu_G(\chi_\pi)=0$ si $\pi$ n'est pas la représentation triviale, que $\mu_G(\chi_\pi)=1$ si $\pi$ est la représentation triviale,
et que $\lambda(\chi_\pi)=\dim \pi$. En supposant que $L(\pi,s)$ satisfasse l'hypothèse de Riemann et la conjecture d'Artin (i.e. soit holomorphe sauf peut-être en $s=1$), le théorème des nombres
premiers généralisé (cf. \cite[Theorem 5.15]{ik}) donne
\begin{eqnarray} \label{515} | \psi(\chi_\pi,x) - \mu_G(\chi_\pi) x | < \ccl{ctp} \, x^{1/2} (\log x) \log(x^{\dim \pi} q(\pi)),\end{eqnarray}
 où $q(\pi)$ est le conducteur d'Artin de $\pi$. On utilise alors la majoration (cf. \cite{effective} ou \cite[Prop 7.4]{eff})
\begin{eqnarray} \label{majorationdisc} \log q(\pi) \leq 2  (\dim \pi)   \, (\log M + \log |G|) \end{eqnarray}
pour obtenir
\begin{eqnarray} \label{chefpsi} | \psi(\chi_\pi,x) - \mu_G(\chi_\pi) x | < \Cr{ct} \, x^{1/2} (\log x) (\dim \pi) (\log x + \log |G|+\log |M|) ,\end{eqnarray}
avec $\ccr{ct}=2\ccr{ctp}$.
Un argument standard d'intégration par partie \cite{davenport} permet d'en déduire:
$$ | \pi(\chi_\pi,x) - \mu_G(\chi_\pi) \Li(x) | <  \ccr{ct} x^{1/2} (\log x)( \dim \pi) (\log x + \log |G|+\log |M|) ,$$
ce qui est la formule voulue et termine la preuve du théorème~\ref{chef}.
\begin{remarque} \label{remMM} Dans le théorème précédent, $M$ est le produit des nombres premiers ramifiés dans $L$. En fait, le théorème reste vrai
si on remplace $M$ par n'importe quel nombre $.0$ divisible par tous les nombres premiers ramifiés dans $L$. C'est clair, car ce changement de
$M$ ne fait pas diminuer le terme d'erreur, et que le changement induit dans $\pi(D,x)$ (qui ne compte que les $p$ ne divisant pas $M$) est au plus $\log M$,
ce qui est absorbé sans effort par le terme d'erreur.
\end{remarque}

\begin{remarque} Si l'on suppose que toutes les fonctions $L$ d'Artin satisfont GRH, mais pas nécessairement la conjecture d'Artin,
il ne semble pas possible de prouver, dans l'état actuel des connaissances, une formule aussi précise que (\ref{515}).
Certes, il résulte de GRH et du fait que $\zeta_L(s)$ n'a qu'un pôle en $s=1$, simple qui plus est, que les fonctions $L$
d'Artin ont tous leur zéros et tous leur pôles non triviaux sur la droite critique $\Re s =1/2$, ce qui permet d'écrire une {\it formule explicite}
pour $\psi(\chi_\pi,x)$ avec un terme d'erreur qui est une somme de termes $x^z$ correspondant aux zéros et pôles non-triviaux $z$ 
de $L(\pi,s)$, donc de termes
de module $x^{1/2}$. C'est encourageant, mais pour conclure il faudrait alors connaître une majoration similaire à \cite[(5.33)]{ik} pour le nombre total de zéros et de pôles de $L(\pi,s)$ 
sur un segment $-T \leq \Im s < T$ de la droite critique, et c'est cette majoration qui fait défaut, car la méthode habituelle, 
qui consiste à intégrer $L'/L$ sur un contour adéquat,  ne donnerait qu'une estimation de la {\it différence} du nombre de zéros et de pôles, au lieu de leur {\it somme}. 
\end{remarque}
%\begin{remarque} Notons une lègère variante que nous aurons à utiliser une fois ci-dessous. Soit $R$ un ensemble fini de représentations de dimension finie
%de $G$ {\it non nécessairement irréductibles}, dont les caractères sont . Une fonction centrale $f:G \rightarrow \Cb$ est dite $R$-{\it admissible} si elle peut s'écrire
%\begin{eqnarray} \label{Radm} f = \sum_{\rho \in R} c_\rho \tr \rho \end{eqnarray} avec les $c_\rho \in \Cb$, et on 
%pose alors $\lambda_R(f) = \sum_{\rho \in R} |c_\rho| \dim \rho$ de $f$
%Alors, en admettant (GRH) et la conjecture d'Artin pour les représentations $\rho \in R$ (i.e. l'holomorphie de $L(\rho,s)$ sauf peut-être en $s=1$), on a pour tout $f$ $R$-admissible
%\begin{eqnarray} \label{chefe} \left| \pi(f,x) - \mu(f) \Li(x)  \right| < \ccr{ct} \, x^{1/2}\,  \lambda_R(f)\, (\log x+\log M+\log |G| ) \end{eqnarray}
%La preuve est exactement la même que celle du théorème~\ref{chefe}: on se ramène par linéarité au cas $f=\tr \rho$, et ce cas se traite exactement comme celui
%où $f=\tr \pi$, $\pi \in \Gd$. L'intérêt de cette variante est qu'il est parfois plus facile de vérifier la conjecture d'Artin pour les représentations non irréductibles $\rho \in R$ que pour leurs composantes irréductibles. 
%
%Ainsi, si $f$ satisfait 

Si l'on combine le théorème~\ref{chef} avec l'estimation de $\lambda(f)$ de la proposition~\ref{majcs} (Cauchy-Schwarz), on obtient
\begin{cor}[Chebotarev effectif, "version de Murty-Murty-Saradha"] \label{cormms}
 Supposons vraie (GRH) et la conjecture d'Artin pour les
 fonctions $L$ d'Artin associées aux représentations irréductibles de $\Gal(L/\Q)$ qui appartiennent au support spectral de $f$. On a, pour $x \geq 3$,
\begin{eqnarray} \label{chefemms} \left| \pi(f,x) - \mu(f) \Li(x)  \right| < \ccr{ct} x^{1/2} \NN{f}_2 \sqrt{|G|}  (\log x +\log M +\log |G| ) \end{eqnarray}
et si $D$ est un sous-ensemble de $G$ invariant par conjugaison.
\begin{eqnarray} \label{chefemmsd} \left| \pi(D,x) - \frac{|D|}{|G|} \Li(x)  \right| < \ccr{ct} x^{1/2}  \sqrt{|D|}  (\log x +\log M +\log |G| ) \end{eqnarray}
\end{cor}
En utilisant la majoration triviale, on obtient une conclusion moins forte, mais qui a l'avantage d'être valable sans la conjecture d'Artin, et de ne supposer GRH
que pour les fonctions Zeta des corps de nombres:
\begin{theoreme}[Chebotarev effectif, "version de Lagarias-Odlyzko-Serre"]
 Supposons vraie (GRH) pour la fonction du corps $L$. On a, pour $x \geq 3$,
\begin{eqnarray} \label{chefelosf} \left| \pi(f,x) - \mu(f) \Li(x)  \right| < \ccl{ctlos} x^{1/2} \NN{f}_1 |G|  (\log x +\log M +\log |G| ) \end{eqnarray}
et si $D$ est un sous-ensemble de $G$ invariant par conjugaison.
\begin{eqnarray} \label{chefelosd} \left| \pi(D,x) - \frac{|D|}{|G|} \Li(x)  \right| < \ccr{ctlos} x^{1/2}  |D| (\log x +\log M +\log |G| ) \end{eqnarray}
\end{theoreme}
Pour (\ref{chefelosd}), cf. \cite[Remarques 3), page 134]{serre}, et (\ref{chefelosf}) en découle par additivité. 

Notons qu'une version très légèrement plus précise  de ces deux théorèmes (cf. \cite[Théorème 4]{serre} and \cite[Proposition 3.2]{effective}) est en fait énoncée dans les articles cités: le facteur logarithmique $(\log x +\log M +\log |G|)$ du terme d'erreur est remplacé par $(\log x + \frac{\log d_L}{|G|})$, où $d_L$ est la valeur absolue du discriminant de l'extension
$L/\Q$. D'après Hecke et Serre (\cite{serre}), on a $\frac{1}{2} \log M \leq \frac{\log d_L}{|G|} \leq \log M + \log |G|$. Ces estimations montrent
que  nos énoncés "perdent" par rapport à leurs variantes plus précises au pire un facteur $\log |G|$, ce qui n'a pas beaucoup d'importance dans la suite. 
De plus, dans la plupart des applications que nous avons en vue $\frac{\log d_L}{|G|}$ est du même ordre de grandeur que $\log M + \log |G|$ (par exemple pour
les corps cyclotomiques, ou bien les corps engendrés par les coordonnées des points de torsion d'une courbe elliptique), donc on ne perd même pas ce facteur $\log |G|$. 

\subsection{Plus petit nombre premier d'un ensemble frobénien}

Rappelons le théorème que nous devons prouver:
\renewcommand{\thetheoremerep}{\ref{thmpluspetit}}
\begin{theoremerep}  Soit $L/\Q$ une extension finie galoisienne. $G=\Gal(L/\Q)$, 
$D \subset G$ non-vide et invariant par conjugaison, et $M=\prod_{p \text{ ramifié dans } L} p$. Supposons vraie (GRH) et la conjecture d'Artin pour les fonctions $L$ d'Artin associées aux représentations irréductibles de $G$. Alors 
le plus petit nombre premier
$p$ tel que $\Frob_p \in D$ vérifie:
$$p < \ccr{cpp} \varphi_G(D)^2 (\log M + \log |G|)^2.$$
\end{theoremerep} 
\begin{pf} 
Par la proposition~\ref{extf}, il existe une fonction centrale $f$ sur $G$, satisfaisant:
\begin{itemize} 
 \item[(a)] si $f(g)>0$, alors $g \in D$;
\item[(b)] $\sum_{g \in G} f(g) > 0$;
\item[(c)] $\varphi_G(D)=\frac{\lambda_G(f)}{\sum_{g \in G} f(g)}$.
\end{itemize} 

Définissons, pour $x \in \R^+$, 
\begin{eqnarray*} \psi_1(f,x) &=& \sum_{n=1}^x \Lambda_f(n) (x-n) \\
\theta_1(f,x) &=& \sum_{p < x} \Lambda_f(p) (x-p) 
\end{eqnarray*}

Si $\pi$ une représentation irréductible de $G$, de caractère $\chi_\pi$, on a par le lemme~\ref{lemmepsi1} ci-dessous, 
$$|\psi_1(\chi_\pi,x) - \mu_G(\chi_\pi) x^2 | < \ccr{cpsi1}\, x^{3/2} (\log q(\pi) + \dim \pi),$$
et donc  par \cite[Prop. 2.5]{effective}, on en déduit 
$$| \psi_1(\chi_\pi,x) - \mu_G(\chi_\pi)x^2| < \ccr{cpsi1}\, x^{3/2} (\dim \pi)  (\log M+\log|G|),$$
d'où, par linéarité et définition de $\lambda_G(f)$
\begin{eqnarray} \label{estpsi1}  | \psi_1(f,x) - \mu_G(f) x^2|   < \ccr{cpsi1}\, x^{3/2} \lambda_G(f)  (\log M + \log |G|  ) \end{eqnarray}
On a \begin{eqnarray*} |\theta_1(f,x) - \psi_1(f,x) |  &=&  \sum_{p,k \geq 2, p^k <x} \Lambda_f(p^k) (x-p^k) \\
&\leq & x  \sum_{k \geq 2} \psi(f,x^{1/k})
\end{eqnarray*}
On a par le théorème~\ref{chef} et sa preuve (cf.~(\ref{chefpsi}))
$$\psi(f,x^{1/2}) < \ccl{cpsixh} (x^{1/2} + x^{1/4} \log^2 x \, \lambda_G(f) (\log M+\log|G|)),$$
et pour $k \geq 3$,
$$\psi(f,x^{1/k}) < \ccr{cpsixh} (x^{1/3}+x^{1/6}  \log^2 x \, \lambda_G(f) (\log M+\log|G|)).$$
Comme par ailleurs $\psi(f,x^{1/k})=0$ si $x^{1/k} < 2$ i.e. si $k > \log x / \log 2$, on obtient:
\begin{eqnarray*} |\theta_1(f,x) - \psi_1(f,x) |  & < & \frac{\ccr{cpsixh}}{\log 2} \left( x^{3/2}+x^{5/4} \log^2 x \lambda_G(f) (\log M+\log |G| ) \right. \\ & &+ \left. x^{4/3} \log x + x^{7/6} \log^2 x   \lambda_G(f) (\log M+\log |G|)\right)  \\ &<& \cc (x^{3/2} \lambda_G(f)  (\log M)  \end{eqnarray*}

On obtient donc pour $\theta_1$ le même estimation que pour $\psi_1$, avec une constante différente :
\begin{eqnarray} \label{theta_1}  |\theta_1(f,x) - \mu_G(f) x^2 | < \ccl{cpps} x^{3/2} \lambda_G(f)  (\log M + \log|G|) \end{eqnarray}
Puisque $\mu_G(f) >0$ par (b), cette estimation montre
$$\theta_1(f,x) > 0 \text{  dès que } x >  \Cr{cpps} ^2 \frac{\lambda_G(f)^2}{\mu_G(f)^2}  (\log M + \log |G|)^2,$$
donc dès que $ x >  \Cr{cpps} ^2 \varphi_G(f)^2  (\log M)^2$.
Par définition de $\theta_1$, si $\theta_1(x)>0$ il existe un nombre premier $p<x$ tel que $f(\Frob_p) > 0$, ce qui par (a) entraîne que $\Frob_p \in D$.
Le résultat suit, avec $\ccr{cpp}=\ccr{cpps}^2$.
\end{pf}

\begin{lemme} \label{lemmepsi1}
Soit $\pi$ une représentation irréductible de $G$ et supposons que $L(\pi,s)$ satisfasse (GRH) et la conjecture d'Artin. On a :
$$|\psi_1(\chi_\pi,x) - \mu_G(\chi_\pi) x^2 | < \ccl{cpsi1} x^{3/2} (\log q(\pi) + \dim \pi),$$
\end{lemme}
\begin{pf} 
Partons de l'égalité élémentaire
$$- \frac{L'}{L}(\pi,s) = \sum_{n \geq 1} \Lambda_{\chi_\pi}(n) n^{-s}.$$
Multipliant les deux membres par $ \frac{x^{s+1}}{s(s+1)} $, où $x \geq 3$ est un paramètre réel, et intégrant sur la droite $\Re s =2$, on obtient:
$$\psi_1(\chi_\pi,x) = \frac{1}{2 i \pi} \int_{2-i \infty}^{2 + i \infty} \frac{L'}{L}(\chi_\pi,s) \frac{x^{s+1}}{s(s+1)} ds.$$
En déplaçant la droite $\Re s = 2$ infiniment loin vers la gauche, on écrit $\psi_1(\chi_\pi,x)$ comme la somme des résidus aux pôles de $\frac{L'}{L}(\chi_\pi,s) \frac{x^{s+1}}{s(s+1)}$; le pôle simple de $L(\pi,s)$ en $s=1$ quand $\pi$ est triviale, ou son absence quand $\pi$ est non-triviale, donne un terme $\mu_G(\chi_\pi) x^2$, et tous les autres pôles triviaux de $\frac{L'}{L}(\chi_\pi,s) \frac{x^{s+1}}{s(s+1)}$ (le plus important étant celui en $s=0$) donnent au total une erreur $\leq \cc (\log q(\pi)) x\log x$ ({\it compare} \cite[chapter 19]{davenport}). Donc
$$\psi_1(\chi_\pi,x) - \mu_G(\chi_\pi) x^2 - \sum_{z} \epsilon(z) \frac{x^{z+1}}{z(z+1)} < \cc (\log q(\pi))  x \log x,$$
où la somme porte sur les zéros et les pôles $z$ non triviaux de de $L(\pi,s)$, comptés plusieurs fois selon leur multiplicité, et où le signe $\epsilon(z)$ est 1 pour un zéro, $-1$ pour un pôle. Puisque nous supposons vraie la conjecture d'Artin,  les pôles non triviaux n'existent pas, et puisque nous supposons que $L(\pi,s)$ satisfait (GRH), tous les zéros non triviaux  satisfont $\Re z = 1/2$, donc $|x^{z+1}| = x^{3/2}.$ Par ailleurs, le nombre de zéros $z$ tels que $T \leq |\Im z| < T+1$
est $\leq \cc (\log q(\pi) + \dim \pi \log(T+3))$ d'après \cite[Prop. 5.7(1)]{ik} (on utilise encore une fois la conjecture d'Artin pour s'assurer que  $L(\pi,s)$ satisfait aux conditions imposées aux fonctions $L$ dans tout le chapitre 5 de \cite{ik}.)
On a  donc  $$\sum_{z} \frac{x^{z+1}}{z(z+1)} < \cc x^{3/2} \sum_{T=1}^\infty \frac{\log q(\pi) + \dim \pi (\log (T+3))}{T^2} < \cc x^{3/2} (\log q(\pi) + \dim \pi),$$ et le résultat suit.
\end{pf}
\begin{remarque} Le théorème précédent reste vrai si l'on remplace $M$ par n'importe quel entier positif divisible par tous les nombres premiers ramifiés dans $L$.
Il suffit dans la preuve de voir que la formule \ref{theta_1} reste vraie avec la nouvelle définition de $M$ (qui change aussi la définition de $\theta_1$). C'est clair,
car le changement de $M$ ne fait que grandir le terme d'erreur, et change $\theta_1(x)$ au plus par $x \log M$, ce qui peut être absorbé dans le terme d'erreur.
\end{remarque}

\begin{remarque} \label{remarqueppm}
En supposant seulement (GRH) pour la  fonction Zeta de $L$, Lagarias et Odlyzko obtiennent le résultat
\begin{eqnarray} \label{ppplo}  p < \ccl{cpplo}(\log d_K)^2 \leq \ccr{cpplo} |G|^2 (\log |G| + \log M)^2 \end{eqnarray}
moins précis sauf quand $|D|=1$.
Sous les mêmes hypothèses et conjectures que le théorème~\ref{thmpluspetit}, Murty et Murty énoncent (sans donner la preuve):
\begin{eqnarray} \label{pppmm} p < \ccl{cppmm} \frac{ |G|^2}{|D|} (\log |G| + \log M)^2.\end{eqnarray}
Comme nous l'avons expliqué l'introduction, notre résultat est plus fort d'une part car en prenant $f = \un_D$, i.e. en utilisant l'estimation
$\varphi(D) \leq \lambda(D) |G| / |D|$, on obtient un meilleur résultat que celui de Murty et Murty si $\lambda(D) < \sqrt{|D|}$, et d'autre part parce qu'il est souvent possible d'utiliser une fonction $f$ satisfaisant les conditions (a) et (b) ci-dessus qui donne un meilleur résultat que $\un_D$, autrement dit
parce que souvent $\varphi(D) \leq \lambda(D) |G| / |D|$.

Il est important d'observer cependant que cette astuce d'utiliser une fonction $f$ autre que $\un_D$ afin de majorer le plus petit nombre premier $p$ tel que $\Frob_p \in D$ ne donne quelque chose de nouveau que parce qu'on a utilise la forme avec la complexité de Littlewood $\lambda(f)$
du théorème~\ref{chef}, et non la forme de Murty-Murty-Saradha (\ref{chefemms}) où $\lambda(f)$ est remplacé par $\NN{f}_2 \sqrt{|G|}$.
En effet, appliquer le principe de démonstration du théorème~\ref{thmpluspetit} avec (\ref{chefemms}) amènerait à chercher à minimiser
l'expression $\frac{\NN{f}_2}{\mu_G(f)}$ sur les fonctions $f$ satisfaisant les conditions (a) et (b), et on voit aisément que le minimum de cette fonctionnelle est toujours atteint pour $f=\un_D$. (Si $f$ est une fonction satisfaisant (a) et (b), soit $f'$ définie par $f'=\max(f,0)$. Il est clair que $f'$  satisfait aussi (a)  et (b) et  que l'on a $\mu_G(f') \geq \mu_G(f)$ et $\NN{f'}_2 \leq \NN{f}$, d'où $\frac{\NN{f'}_2}{\mu_G(f')} \leq \frac{\NN{f}_2}{\mu_G(f)}$. De plus comme $f'$ est positive et à support dans $D$, $\mu_G(f') = \frac{1}{|G|} \sum_{x \in D} f'(x) \leq \frac{1}{|G|}\sqrt{ \sum_{x \in D} |f'(x)|^2 \sum_{x \in D} 1} = \NN{f'}_2 \sqrt{|G||D|}$ avec égalité quand $f'$ est proportionnel à $\un_D$.)

Cette observation illustre la thèse centrale de cet article: que la version de Chebotarev avec complexité de Littlewood $\lambda(f)$ est à la fois meilleure et plus souple que celle où l'on remplace  $\lambda(f)$ par la norme $L^2$ (ou {\it a fortiori} $L^1$) de $f$.
\end{remarque}

\subsection{Application à certains ensembles de densité $0$}
\label{sectdenszero}

\subsubsection{Position du problème}

Dans cette partie on se donne un entier $M \geq 1$, un ensemble $\Lambda$, et pour chaque $\nu \in \Lambda$, une extension finie galoisienne
 $L_\lambda$ de $\Q$, non ramifiée hors des nombres premiers divisant $M$,  de groupe de Galois $G_\nu$. On se donne également une famille $D=(D_\nu)_{\nu \in \Lambda}$ de sous-ensembles $D_\nu \subset G_\nu$ stables par conjugaison. On pose 
 \begin{eqnarray} \label{tildeDdef} \tilde D = \{p  \,\mid\, p \nmid M \text{ et } \forall \nu \in \Lambda, \Frob_{p,G_\nu} \in D_\nu \  \} = \cap_{\nu \in \Lambda} \tilde D_\nu. \end{eqnarray} et \begin{eqnarray} \pi(D,x) = | \tilde D \cap [1,x] | =  \{p \leq x  \,\mid\, p \nmid M \text{ et } \forall \nu \in \Lambda, \Frob_{p,G_\nu} \in D_\nu  \}  \end{eqnarray}
 
Si $\inf_\nu \frac{|D_\nu|}{|G_\nu|} = 0$ (ce qui sera le cas dans les applications), l'ensemble $\tilde D$ a une densité naturelle égale à zéro, et notre but
est d'estimer sa "rareté" en donnant une majoration  (valable sous GRH et la conjecture d'Artin pour les fonctions $L$ d'Artin associées aux représentations de $G_\nu$, $\nu \in \Lambda$) de la forme $\Pi_D(x) =\bigo{x^u (\log x)^v} $, où $0 \leq u<1$ et $v$ sont des réels dépendant de la taille des $G_\nu$, $D_\nu$ et $\lambda_{D_\nu}(G_\nu)$.

Nous donnons ci-dessous deux théorèmes donnant de telles estimations. L'idée de la preuve du premier est très simple, et nous l'avons reprise de \cite{serre}; c'est elle aussi qui est utilisée par Murty-Murty-Saradha dans \cite{effective}.
Pour majorer $\pi_D(x)$, on observe que $\pi(D,x) \leq \pi(D_\nu,x)$ pour tout $\nu$ et on choisit le $\nu=\nu(x)$ qui donne le meilleur résultat. Le théorème (cf. théorème 
\ref{denszero})  que nous obtenons est donc est donc meilleur que celui énoncé dans \cite{serre} (resp. celui implicite dans \cite{effective}) quand 
et dans la mesure où nous disposons d'une estimation des $\lambda(D_\nu)$ meilleure que l'estimation triviale (resp. de Cauchy-Schwarz).

Il est clair que le méthode décrite ci-dessus peut être améliorée si l'on arrive à utiliser l'information $\tilde D \subset \tilde D_\nu$ pour tous les $\nu$ à la fois, ou du moins pour beaucoup de $\nu$, et non seulement pour un seul. C'est difficile, car ces informations sont en partie redondantes. Heureusement, il existe une théorie bien documentée pour traiter ce genre de question: la méthode du crible. Plus précisément, comme les $|G_\nu - D_\nu|/|G_\nu|$ ne tendent pas
vers $0$ dans nos hypothèses, c'est le {\it grand crible} qui s'applique ici. Ce n'est pas la première fois que le grand crible est appliqué à des questions de ce genre: cf. \cite{kowalski}  et surtout  le travail \cite{zywina} de Zywina. Cependant, le grand crible que nous utilisons plus bas, et qui donnera, sous des hypothèses un peu plus fortes mais souvent satisfaites en pratique, un résultat (théorème~\ref{denszerocrible}) meilleur que celui du théorème~\ref{denszero}, n'utilise pas comme Zywina le {\it crible de conjugaison} (cf. \cite[3.1]{kowalski}) qui par nature ne peut donner qu'une estimation qui dépend de $|D_\nu|$ et $|G_\nu|$ mais pas de $\lambda_{G_\nu}(D_\nu)$. Nous renvoyons à la preuve du théorème~\ref{denszerocrible} et à la remarque qui la suit pour de plus amples détails sur notre méthode et comment elle permet d'obtenir des résultats plus forts que ceux de Zywina.

 \subsubsection{Borne obtenue par application directe de Chebotarev}

\begin{theoreme} \label{denszero}
Soit $\alpha$ et $\beta$ deux réels tels que $0<\alpha,\beta \leq 1$. On fait les hypothèses suivantes:
\begin{itemize}
\item[(a)] Il existe un nombre $L>0$ tels que tout intervalle de $\R$ de longueur $L$ contienne au moins un $\log |G_\nu|$
\item[(b)] Il existe un nombre $Q>0$ tel que $|D_\nu|/|G_\nu| \leq Q/|G_\nu|^\alpha$.
\item[(c)] Il existe un nombre $R>0$ tel que $\lambda_{G_\nu}(D_\nu) \leq R |G_\nu|^\beta$.
\end{itemize}
Supposons que toutes les fonctions $L$ d'Artin des représentations de  $G_\nu=\Gal(L_\nu/\Q)$ satisfassent (GRH) et la conjecture d'Artin. Alors 
on a \begin{eqnarray} \label{formdenszero} \pi_D(x) =\bigo{ x^{\frac{\alpha+2 \beta}{2 \alpha+2 \beta}} (\log x)^{\frac{\alpha-\beta}{\alpha+\beta}} }. \end{eqnarray}
\end{theoreme}
\begin{pf}  Fixons $x$ assez grand.  Choisissons $\nu =\nu(x)$ tel que $$e^{-L} x^{1/(2\alpha+2\beta)} (\log x)^{-2/(\alpha+\beta)}
\leq |G_\nu| \leq  x^{1/(2\alpha+2\beta)} (\log x)^{-2/(\alpha+\beta)},$$ ce qui est possible par (a). Appliquons le théorème de Chebotarev:
$$ \pi(D_\nu,x) < \frac{|D_\nu|}{|G_\nu|} \Li(x)  +  \Cr{ct} (x^{1/2} \lambda(D_\nu) (\log x +\log M+\log |G_\nu| ).$$
Majorons le terme principal:
\begin{eqnarray*}
\frac{|D_\nu|}{|G_\nu|} \Li(x) & \leq & Q x (\log x)^{-1} |G_\nu|^{-\alpha} \\ &\leq & Q e^{\alpha L} x^{1-\alpha/(2\alpha+2 \beta)} (\log x)^{-1 + 2\alpha/(\alpha+\beta)} \\
&  = &  \bigo{ x^{(\alpha+2 \beta)/(2 \alpha+2 \beta)} (\log x)^{\frac{\alpha-\beta}{\alpha+\beta}}}
\end{eqnarray*}

Majorons le terme d'erreur: comme $\log |G_\nu| = O( \log x)$ et $\log M$ est constant, $\log x+\log M +\log |G_\nu| = O( \log x)$ et
\begin{eqnarray*}
x^{1/2} \lambda(D_\nu) (\log x+\log M +\log |G_\nu|) & = & \bigo{ x^{\frac{1}{2}} |G_\nu|^\beta (\log x) } \\
&=& \bigo{ x^{\frac{1}{2}} (\log x) x^{\beta/(2\alpha+2\beta)} \log(x)^{- 2 \beta/(\alpha+\beta)}}\\
&=& \bigo{ x^{(\alpha+2 \beta)/(2 \alpha+2 \beta)} (\log x)^{\frac{\alpha-\beta}{\alpha+\beta}}}.
\end{eqnarray*}
D'où $\pi(D,x) \leq \pi(D_\nu,x) = \bigo{x^{(\alpha+2 \beta)/(2 \alpha+2 \beta)} (\log x)^{\frac{\alpha-\beta}{\alpha+\beta}}}$ et le
 théorème est prouvé.
\end{pf}

\begin{remarques} 
\begin{itemize}
\item[(i)] En utilisant la majoration triviale $\lambda(D_\nu) \leq |D_\nu|$, on voit que (b) implique (c) avec $\beta = 1-\alpha$.
On voit donc que sous les hypothèses (a) et (b), le théorème entraîne  
$$\pi(D,x) < \cc x^{1-\frac{\alpha}{2}} (\log x)^{2 \alpha-1}.$$
 Ce résultat est celui qu'obtient Serre (\cite[Théorème 11(ii)]{serre}), sous GRH mais sans supposer Artin.
 
\item[(ii)] Si au lieu d'utiliser la majoration triviale, on utilise la majoration de Cauchy-Schwartz $\lambda(D_\nu) \leq |D_\nu|^{1/2}$,
 on voit que (b) implique (c) avec $\beta = (1-\alpha)/2$, et donc que le théorème donne, sous les hypothèses (a) et (b):
\begin{eqnarray} \label{denszeromms} \pi(D,x) < \cc x^{\frac{1}{1+\alpha}} (\log x)^{\frac{3 \alpha-1}{\alpha+1}}.\end{eqnarray}
 Ce résultat est celui qu'on obtient avec la version de Murty-Murty-Saradha du théorème de Chebotarev (en utilisant (GRH) et Artin). Bien qu'il ne soit pas explicitement énoncé en général dans~\cite{effective}, il l'est dans le cas particulier où les $G_\nu$ sont des sous-groupes de $\Gl_2(\F_\ell)$ attachées à une forme modulaire propre pour les opérateurs de Hecke, cas dans lequel on a $\alpha=1/4$: voir la première partie du premier théorème page 3 de \cite{effective}.
 Une astuce que nous utiliserons également permet d'ailleurs aux auteurs dans ce cas précis de se passer de la conjecture d'Artin.
 
 \item[(iii)] L'exposant de $x$ dans le théorème est toujours $\geq 1/2$ avec égalité si $\beta=0$. \end{itemize}
\end{remarques}

\subsubsection{Borne obtenue par une application du grand crible}

Nous énonçons maintenant la majoration de $\pi(D,x)$ utilisant le grand crible. Comme les hypothèses que nous utilisons ressemblent à celles du 
théorème~\ref{denszero}, nous les numérotons d'une manière qui tente de mettre en lumière cette ressemblance.

\begin{theoreme} \label{denszerocrible}
Soit $d \geq 1$ un entier, et soit $\alpha$ et $\beta$ deux réels tels que $0<\alpha,\beta \leq 1$. On suppose que $\Lambda$ est l'ensemble de tous
les nombres premiers ne divisant pas un entier fixé $N \geq 1$ (on notera $\ell$ au lieu de $\nu$ un élément générique de $\Lambda$) et on fait les hypothèses suivantes:
\begin{itemize}
\item[(b')] Il existe des nombres $P>0$ et $P'>0$ tel que pour tout $\ell \in \Lambda$, $\ell^{-d \alpha} - P \ell^{-d\alpha -1} \leq  \frac{|D_\ell|}{|G_\ell|}  \leq  \ell^{- d \alpha} + P' \ell^{-d \alpha-1}$.
\item[(c')] Il existe un nombre $R>0$ tel que pour tout $\ell \in \Lambda$,  $\lambda_{G_\ell}(D_\ell) \leq R |G_\ell|^\beta$.
\item[(d')] Pour tout sous-ensemble fini $m$ de $\Lambda$, l'application naturelle $\GQ \rightarrow \prod_{\ell \in m} G_\ell$ est surjective.
\end{itemize}
Supposons que toutes les fonctions $L$ d'Artin des représentations de  $G_\ell$, $\ell \in \Lambda$ satisfassent (GRH) et la conjecture d'Artin.
Alors on a, pour tout $\epsilon>0$:
$$\pi(D,x) =\bigo{  x^{\frac{d\alpha+ 2d \beta + 1}{2d \alpha+2d\beta+2}+\epsilon }     }$$
\end{theoreme}

\begin{pf} On utilise la méthode du grand crible, pour laquelle notre référence  est \cite[\S2.1 sqq]{kowalski}, dont on reprend les notations et la terminologie.

\par \bigskip
{\bf Mise en place du crible:}
Tout d'abord, la {\it donnée de grand crible} $(Y,\Lambda,(Y_\ell)_{\ell \in \Lambda},(\pi_\ell)_{\ell \in \Lambda})$ que nous utiliserons est la suivante:
\begin{itemize}
\item L'ensemble $Y$ est l'ensemble des nombre premiers $p$ ne divisant pas $M$.
\item L'ensemble $\Lambda$ est comme dans l'énoncé du théorème: $\Lambda=\{ \ell \, \mid \, \ell \nmid N\}$.
\item Pour $\ell \in \Lambda$, on définit un ensemble $Y_{\ell} = \{\oui,\non\}$, et une application surjective $\pi_\ell : Y \rightarrow Y_\ell$,
telle que $\pi_\ell(p)$ soit la réponse à la question: "est-ce que l'élément $\Frob_p$ de $G_\ell$ appartient à $D_\ell$"?
\end{itemize}
Avant d'aller plus loin, munissons $Y_\ell$ de la mesure de probabilité $\mu_\ell$ tel que $\{\oui\}$ soit de mesure $\mu_\ell(\{\oui\})=\frac{|D_\ell|}{|G_\ell|}$.  Ainsi, l'application naturelle 
$\psi_\ell: \ G_\ell \rightarrow Y_\ell$ qui envoie $D_\ell$ sur $\oui$ et son complémentaire sur $\non$ réalise l'espace mesuré $Y_\ell$ comme un quotient de l'espace 
$G_\ell$ muni de sa mesure invariante de probabilité.

Définissions ensuite l'{\it ensemble à cribler} ("siftable set" en anglais) $X$:
$$ X = \{p \in Y, p \leq x\}.$$ 
Munissons cette ensemble de la mesure de comptage (notée $\mu$ dans \cite{kowalski} mais pour laquelle nous n'utiliserons pas de notation), et de son injection canonique dans $Y$ (notée $F$ dans \cite{kowalski}, sous-entendue ici). Le diagramme commutatif suivant a pour but d'aider à visualiser les
différents objets introduits:
$$\xymatrix{ X \ar@{^{(}->}[rr] & & Y \ar[rr]^{p \mapsto \Frob_{p,G_\ell}} \ar[rrd]^{\pi_\ell} & & G_\ell \ar[d]^{\psi_\ell} \\ & & & & Y_\ell = \{\oui,\non\}}$$

Soit $\Lc^\ast$ un sous-ensemble fini de $\Lambda$ (le {\it support premier du crible}) et $\Lc$ un ensemble de parties de $\Lc^\ast$ (le {\it 
support du crible}), que nous choisirons tous deux plus tard. Finalement, soit $\Omega=(\Omega_\ell)_{\ell \in \Lc^\ast}$ la {\it famille d'ensemble criblants}, avec $\Omega_\ell \subset Y_\ell$ for $\ell \in \Lc^\ast$, qu'on définit simplement comme
$\Omega_\ell = \{\non\}$ pour tout $\ell \in \Lc^\ast$.

On peut alors définir l'ensemble criblé (cf. \cite[definition 2.1]{kowalski} -- cet ensemble dépend des données ci-dessus à l'exception de $\Lc$ et des mesures choisies sur les $Y_\ell$):
\begin{eqnarray*} S(X,\Omega,\Lc^\ast) &=& \{p \in X, \pi_\ell(p) \not \in \Omega_\ell \text{ pour tout } \ell \in \Lc^\ast \} \\
&=& \{p \leq x, p \nmid M, \Frob_p \in D_\ell  \text{ pour tout } \ell \in \Lc^\ast \} \\ & \supset & \tilde D \cap [1,x] \end{eqnarray*}
On a donc 
\begin{eqnarray} \pi(D,x) = |\tilde D \cap [1,x]| \leq  | S(X,\Omega,\Lc^\ast) | \end{eqnarray}

\par \bigskip

{\bf Majoration de la constante de crible $\Delta$}

Aux données ci-dessus on peut attacher une constante de crible $\Delta=\Delta(X,\Lc)$ définie \cite[prop. 2.4]{kowalski}, que nous nous donnons pour but de majorer. Mais pour ce faire, il nous faut d'abord introduire, suivant \cite{kowalski} quelques objets construits à partir des  données précédentes, et des notations pour les nommer.

Pour $\ell \in \Lc^\ast$, soit $\phi_\ell$ la fonction à valeurs réelles sur $Y_\ell=\{\oui,\non\}$ qui à $\oui$ associe $\sqrt{\frac{|G_\ell|-|D_{\ell}|}{|D_\ell|}}$
et à $\non$ associe $- \sqrt{\frac{|D_\ell|}{|G_\ell|-|D_{\ell}|}}$. Pour la structure hermitienne sur $L^2(Y_\ell,\mu_\ell,\Cb)$, cette fonction $\phi_\ell$ est orthogonale à la fonction constante $1$, et de norme $1$. C'est même, à une phase sans importance près,  
l'unique telle fonction dans $L^2(Y_\ell,\Cb)$. On pose $\Bc_{\ell}^\ast = \{\phi_\ell\}$. 

De même, pour $m \in \Lc$, on définit $Y_m = \prod_{\ell \in m} Y_\ell$ (produit comme espace mesuré), $\pi_m: Y \rightarrow Y_m$ comme le produit des $\pi_\ell$, $\ell \in m$, la fonction  $\phi_m \in L^2(Y_m,\Cb)$
par $\phi_m = \otimes_{\ell \in M} \phi_\ell$, et finalement $\Bc_m^\ast = \{\phi_m\}$. On a donc par définition, pour tout $p \in Y$,
$$\phi_m(\pi_m(p)) = \prod_{\ell\in m} \phi_\ell(\pi_\ell(p)).$$ 
On définit $L_m$ comme l'extension de $\Q$ composée des $L_\ell$ pour $\ell \in M$, $G_m = \Gal(L_m/\Q)$, si bien que $G_m = \prod_{\ell \in m} G_\ell$ par l'hypothèse (d'). On a des application naturelles surjective $G_m \rightarrow G_\ell \rightarrow Y_\ell$ pour $\ell \in m$, et donc, en prenant leur produit une application $\psi_m: G_m \rightarrow Y_m$. On a donc pour tout $m$ le diagramme commutatif suivant:
$$\xymatrix{ X \ar[rr] & & Y \ar[rr]^{p \mapsto \Frob_{p,G_m}} \ar[rrd]^{\pi_m} & & G_m \ar[d]^{\psi_m} \\ & & & & Y_m \ar[d]^{\phi_m} \\ & & & &  \R}$$
Ces diagrammes sont compatibles en un sens évident pour $m \subset m'$ et redonne le diagramme précédent pour $m=\{\ell\}$.

 La proposition 2.9 de \cite{kowalski} donne:
\begin{eqnarray} \label{prop29} \Delta \leq \max_{m \in \Lc} \sum_{n \in \Lc} |W(\phi_m,\phi_{n})| \end{eqnarray}
où
$$W(\phi_m,\phi_{n}) = \sum_{p \in X} \phi_m(\pi_m(p)) \overline{\phi_n(\pi_n(p))}$$
Comme $\phi_n$ est à valeurs réelles, on peut supprimer la conjugaison complexe; par la commutativité et la compatibilité des diagramme ci-dessus pour $m$, $n$ et $m\cup n$, il vient
$$ W(\phi_m,\phi_{n}) = \sum_{p \leq x, p \nmid M} (\phi_m \circ \psi_m)(\phi_n \circ \psi_n)(\Frob_{p,G_{m\cup n}})$$
La fonction à valeurs réelles $(\phi_m \circ \psi_m)(\phi_n \circ \psi_n)$ sur le groupe $G_{m \cup n}=G_s \times G_i$, où  $s$ est la différence symétrique de $m$ et $n$, et $i$ leur intersection (i.e. $s=(m \cup n)-(m \cap n),\ \ i=m \cap n$) peut s'écrire
$$(\phi_m \circ \psi_m)(\phi_n \circ \psi_n) = (\phi_s \circ \psi_s) \otimes (\phi_i \circ \psi_i)^2.$$
On voit en particulier que la moyenne des valeurs de cette fonction, $\mu_{G_{m \cup n}}((\phi_m \circ \psi_m)(\phi_n \circ \psi_n) )$
est nulle sauf si $s=\emptyset$, i.e. sauf si $m=n$, auquel cas cette moyenne est $1$. 

Appliquons le théorème de Chebotarev effectif, pour $f=(\phi_m \circ \psi_m)(\phi_n \circ \psi_n)$. On obtient
\begin{eqnarray} \label{chebphi} 
\left| W(\phi_m,\phi_n) - \delta_{m,n} Li(x) \right| <  \ccr{ct} \lambda_{G_{m \cup n}}((\phi_m \circ \psi_m)(\phi_n \circ \psi_n)) x^{1/2} (\log x + \log |G_{m \cup n}| + \log M).
\end{eqnarray} 

Il nous faut maintenant estimer la norme de Littlewood apparaissant dans la formule ci-dessus. On a d'après la proposition~\ref{compprod}
\begin{eqnarray*} \lambda_{G_{m \cup n}}((\phi_m \circ \psi_m)(\phi_n \circ \psi_n)) &=&
\prod_{\ell \in s} \lambda_{G_\ell} (\phi_\ell \circ \psi_\ell)  \times \prod_{\ell \in i} \lambda_{G_\ell} (\phi_\ell^2 \circ \psi_\ell)
\end{eqnarray*}
Or, comme $\phi_\ell \circ \psi_\ell =  \sqrt{\frac{|G_\ell|-|D_{\ell}|}{|D_\ell|}} \un_{D_\ell} +  \sqrt{\frac{|D_\ell|}{|G_\ell|-|D_{\ell}|}} \un_{G_\ell-D_\ell}$, 
on a puisque $\lambda$ est une norme:
\begin{eqnarray*} \lambda_{G_\ell} (\phi_\ell \circ \psi_\ell) &\leq &  \sqrt{\frac{|G_\ell|-|D_{\ell}|}{|D_\ell|}} \lambda(D_\ell)  + \sqrt{\frac{|D_\ell|}{|G_\ell|-|D_{\ell}|}} \lambda(G_\ell-D_\ell) \\ 
\lambda_{G_\ell} (\phi_\ell^2 \circ \psi_\ell) &\leq & \frac{|G_\ell|-|D_{\ell}|}{|D_\ell|} \lambda(D_\ell)  + \frac{|D_\ell|}{|G_\ell|-|D_{\ell}|} \lambda(G_\ell-D_\ell)
\end{eqnarray*}
Or d'après le corollaire~\ref{lambdacomplement} 
$\lambda(G_{\ell}-D_\ell)\leq \lambda(D_\ell)+1$.
On obtient donc facilement, si $\frac{|G_\ell|}{|D_\ell|} \geq 2$ (ce qui par (b') ne peut être en défaut que pour un nombre fini de $\ell$, qu'on peut
alors enlever de l'ensemble $\Lambda$ en changeant l'entier $N$): 
$$  \lambda_{G_\ell} (\phi_\ell^2 \circ \psi_\ell) < \frac{|G_\ell|}{|D_\ell|} \lambda(D_\ell)$$
et 
$$  \lambda_{G_\ell} (\phi_\ell \circ \psi_\ell) < \frac{|G_\ell|}{|D_\ell|} \lambda(D_\ell).$$

En appliquant ces estimations dans (\ref{chebphi}), on obtient
$$ \left| W(\phi_m,\phi_n) - \delta_{m,n} Li(x) \right| < \ccr{ct} \left( \prod_{\ell \in m \cup n} \frac{|G_\ell|}{|D_\ell|} \lambda(D_\ell) \right) x^{1/2} (\log x + \log |G_{m \cup n}|+\log M).$$ D'après (b') on a $\frac{|G_\ell|}{|D_\ell|} < 2 \ell^{d \alpha}$ pour $\ell$ suffisamment grand, donc pour tout $\ell \in \Lambda$
quitte à changer l'entier $N$. On a donc, en utilisant aussi (c'):
$$ \left| W(\phi_m,\phi_n) - \delta_{m,n} Li(x) \right| < \ccr{ct} \left( \prod_{\ell \in m \cup n} (2R) \ell^{d\alpha+d\beta} \right) x^{1/2}
 (\log x + \log |G_{m \cup n}|) $$

Il est temps de restreindre le choix du support premier de crible $\Lc^\ast$, et du support de crible $\Lc$, qui jusqu'à présent était quelconque.
Soit $Q \geq 1$ un paramètre réel, que nous  choisirons ultérieurement. Soit $\Lc^\ast_Q$ l'ensemble des $\ell \in \Lambda$ tel que $\ell \leq Q$, et $\Lc_Q$ l'ensemble des $m \in \Lc^\ast$ tels que $ \prod_{\ell  \in m} \ell  \leq Q$. Pour ce choix de $\Lc^\ast$, $\Lc$, on a donc
\begin{eqnarray*} \left| W(\phi_m,\phi_n) - \delta_{1,d} Li(x) \right|  &<& \ccr{ct} \left(\prod_{\ell \in m \cup n} 2R \right) Q^{d(\alpha+\beta)} x^{1/2} (\log x + d \log Q + \log M) \end{eqnarray*}
Le nombre de facteurs dans le produit  $\prod_{\ell \in m \cup n}$ est le nombre de facteurs premiers d'un entier sans facteur carré plus petit que $Q$. Ce nombre
est donc $\leq c \log Q / \log \log Q$ pour $Q$ suffisamment grand, d'après le théorème des nombres premiers (où $c>0$ est une constante). On a donc 
$\prod_{\ell \in m \cup n} 2R \leq ((2R)^c)^{\log Q / \log \log Q} = O(Q^\epsilon)$ pour tout $\epsilon >0$ (la constante implicite dépendent de $\epsilon$).
On obtient donc
\begin{eqnarray*} \left| W(\phi_m,\phi_n) - \delta_{1,d} Li(x) \right|  &<& \ccl{cRc} Q^{d(\alpha+\beta+\epsilon)} x^{1/2} \log x \end{eqnarray*}
pour $Q$ suffisamment grand, où la constante $\ccr{cRc}$ dépend de $R$, $\epsilon$, $c$ mais pas de $Q$ ni de $x$,
et donc, par (\ref{prop29}) et puisque $|\Lc_Q|<Q$,
\begin{eqnarray} \Delta(X,\Lc_Q) < Li(x) + \ccr{cRc} Q^{d(\alpha+\beta)+1+\epsilon} x^{1/2} \log x  \end{eqnarray}
C'est la majoration de $\Delta$ que nous cherchions.

\par \bigskip
{\bf Minoration de la constante de crible $H$}:

La constante $H$ est définie (cf. \cite[(2.4)]{kowalski}) par
$$H = \sum_{m \in \Lc} \prod_{\ell \in m} \frac { \mu_\ell(\{\non\}) }{\mu_\ell(\{\oui\})} = \sum_{m \in L_c} \prod_{\ell \in m} 
\frac{1 - |D_\ell|/|G_\ell|}{|D_\ell|/|G_\ell|} $$

Posons $f(m)
= \frac{1}{m^{d\alpha}} \prod_{\ell \mid m, \ell \nmid N} 
\frac{1 - |D_\ell|/|G_\ell|}{|D_\ell|/|G_\ell|}$ pour $m$ sans facteur carré.
En particulier, pour $\ell$ un nombre premier ne divisant pas $N$,
$f(\ell)=\frac{1 - |D_\ell|/ |G_\ell|}{\ell^{d \alpha}  |D_\ell|/|G_\ell| }$, et $f(\ell)=1$ si $\ell \mid N$. 
Par l'hypothèse (b'), on a donc $| f(\ell) - 1| = O( \ell^{-1})$ quand $\ell \rightarrow \infty$.  
D'après le théorème de Lau et Wu (\cite[Theorem 1]{lauwu} or \cite[Theorem G.1]{kowalski}), on a
$$\left| \sum^\flat_{m \leq Q} f(m) - c Q \right| <  \cc Q (\log Q)^{-1} \log\log Q,$$
où $c = \prod_{\ell} (1-\ell^{-1})(1+f(\ell)\ell^{-1})$ est une constante strictement positive et $\sum^\flat$ désigne la somme restreinte aux entiers sans facteurs carrés.
On a alors
\begin{eqnarray} 
H &=&  \sum^\flat_{m \leq Q} m^{d \alpha} f(m) \\
&\geq &   (Q/2)^{d \alpha} \sum^\flat_{Q/2 \leq m \leq Q}  f(m)\\
& \geq &  (Q/2)^{d \alpha} c Q/2 + O( Q (\log Q)^{-1} \log\log Q),
\end{eqnarray} 
ce qui implique 
$$H^{-1}  < \ccl{chinv} Q^{-d \alpha -1 }$$
quand $Q$ est assez grand, la constante $\ccr{chinv}$ ne dépendant que de $P$, $P'$ et $N$.

\par \bigskip
{\bf Fin de la preuve:} 
D'après le théorème du grand crible (\cite[prop. 2.3]{kowalski}),
on a
$\Pi(D,x) <  \Delta H^{-1}$, ce qui donne, utilisant notre majoration de $\Delta$ et notre minoration de $H$:
$$\Pi(d,x) < \ccr{chinv} \frac{\Li(x)}{Q^{1 + \alpha d}} + \ccr{chinv}\ccr{cRc} Q^{d\beta+\epsilon} x^{1/2} \log x,$$
pour $x \geq 3$ et $Q$ suffisament grand. Prenons alors 
 $Q=Q(x) = x^{\frac{1}{2d\beta+2d\alpha+2}}$, ce qui est possible si l'on suppose que $x$ est suffisamment grand.
 On obtient alors $$\pi_D(x) < \left(\frac{\ccr{chinv}}{\log x} + \ccr{chinv}\ccr{cRc} \log x \right)  x^{\frac{d\alpha+ 2d \beta + 1+\epsilon}{2d \alpha+2d\beta+2} }$$
 pour $x$ suffisamment grand, et le théorème.
 \end{pf}

\begin{remarque}
\begin{itemize}
\item[(i)] L'hypothèse (d') est évidemment une hypothèse d'indépendance: elle revient à dire que les extensions $L_\ell$ de $\Q$ pour $\ell \in \Lambda$ sont
linéairement disjointe. 
\item[(ii)] Les hypothèses (b') et (c') du théorème~\ref{denszerocrible} impliquent évidemment les hypothèses (b), et (c)  du théorème~\ref{denszero} avec les mêmes valeurs de $\alpha$ et $\beta$.
Sous ces hypothèses légèrement plus restrictives et l'hypothèse d'indépendance (d'), la conclusion du  théorème~\ref{denszerocrible} est plus forte que celle du théorème~\ref{denszero}
car on a $\frac{d\alpha+ 2d \beta + 1}{2d \alpha+2d\beta+2} < \frac{d\alpha+ 2d \beta}{2d \alpha+2d\beta}$ si $\beta>0$.
\item[(iii)] Comme pour le théorème~\ref{denszero}, au lieu d'utiliser une estimation maison (c') de $\lambda(D_\ell)$, on peut sous les hypothèses (a') et (b') seules, déduire que (c') est vrai avec $\beta=1-\alpha$ (si l'on utilise la majoration triviale $\lambda(D_\ell) \leq |D_\ell|$) ou même avec $\beta=(1-\alpha)/2$
(si l'on utilise la majoration de Cauchy-Schwarz). On obtient alors, pour tout $\epsilon>0$,
\begin{eqnarray} \label{denszerocribletriv} \pi_D(x) &=& \bigo{ x^{ \frac{2d - d \alpha + 1}{2d +2} +\epsilon } } \text{ avec la majoration triviale,}\\
\label{denszerocriblecs} \pi_D(x) &=& \bigo{ x^{\frac{d+1}{d +d\alpha+2}+\epsilon }  }\text{ avec la majoration de Cauchy-Schwarz}.\end{eqnarray}
\item[(iv)] Dans \cite{zywina}, une version différente du grand crible est utilisée pour le même problème dans un cas particulier. Au lieu de
prendre pour $Y_\ell$ un ensemble à deux éléments comme nous l'avons fait, avec $\Omega_\ell$ l'un de ces éléments.
Zywina prend pour $Y_\ell$ l'ensemble des classes de conjugaisons du groupe $G_\ell^\sharp$ et pour $\Omega_\ell$ l'image du complémentaire de $D_\ell$ dans $G_\ell^\sharp$.
Ainsi, notre choix $(Y_\ell,\{ \non \})$ est un quotient du choix $(G_\ell^\sharp,\Omega_\ell)$ par la relation d'équivalence évidente (celle dont les classes d'équivalences sont $D_\ell$ et $G_\ell-D_\ell$), et ceci est également vraie pour les mesures que Zywina et nous mettons sur
$G_\ell^\sharp$ et $Y_\ell$. Un aspect important du choix de Zywina est que l'estimation donné par le grand crible ne dépend que de la taille  $|D_\ell|$ des $D_\ell$, non de leur position dans $G_\ell$. Autrement dit, toute information qu'on pourrait avoir sur $\lambda(D_\ell)$ est perdue. Ceci provient du fait que
la définition de la constante $\Delta$ ne fait pas intervenir les $\Omega_\ell=D_\ell$, tandis que celle de la constante $H$ ne fait intervenir que la taille des $\Omega_\ell$ ou plus précisément, leur mesure dans $Y_\ell$. Cet aspect, qui souvent est un des avantages du crible (cf. \cite[\S2.5]{kowalski}),
est ici évidemment un inconvénient.

Ainsi, dans le cas particulier considéré par Zywina, on a $d=4$, $\alpha=1/4$, la borne que Zywina obtient est $\Pi_D(x)=O(x^{4/5})$ (négligeant les termes logarithmiques) est bien celle qu'on obtient à partir de (\ref{denszerocribletriv}), i.e. avec le crible et la majoration triviale de $\lambda(D_\ell)$. Cette borne
se trouve être la même que la majoration de \cite{effective} (cf. (\ref{denszeromms})). En d'autres termes, ce que le crible fait gagner à l'estimation de Zywina par rapport à celle de Murty-Murty-Saradha est reperdue par l'emploi implicite d'une moins bonne majoration de $\lambda(D_\ell)$. 

\item[(v)] Si l'on est prêt à admettre (GRH)  mais pas le conjecture d'Artin, on peut quand même démontrer une version
plus faible du théorème précédent: la même conclusion reste vraie à condition de renforcer l'hypothèse (c') sur $\beta$ en demandant que
\begin{itemize}
\item[(c'')] il existe $R>0$ tel que pour tout $\ell$ il existe un sous-groupe normal $U_\ell$ de $G_\ell$ tel que $U_\ell D_\ell=D_\ell$
et $|D_\ell|/|U_\ell|  < R |G_\ell|^\beta$. 
\end{itemize}
Sous cette condition, la même preuve s'applique sans utiliser Artin, à condition d'appliquer
le théorème de Chebotarev effectif de Lagarias-Odlyzko-Serre (\ref{chefelosd}) au lieu du thèorème~\ref{chef}, et de l'appliquer au groupe $G_m / \prod_{\ell \in m} U_\ell$ au lieu de $G_m$, et aussi de remplacer partout dans la preuve la norme de Littlewood $\lambda$ par la norme $L^1$.
\end{itemize}
\end{remarque}

 \section{Applications}
 
 \subsection{Plus petit nombre premier générateur de $(\Z/q\Z)^\ast$}
 
Soit $\ell$ un nombre premier impair. On sait depuis Gauss que le groupe $(\Z/\ell \Z)^\ast$ est cyclique, et ce résultat suggère immédiatement la question "pratique" suivante:  quel est le plus entier positif $g(\ell)$ tel que $a \pmod{\ell}$ soit un générateur de $(\Z/\ell\Z)^\ast$? Notons aussi $g^\ast(\ell)$ le plus petit nombre {\it premier} $p$ tel que  $a \pmod{\ell}$ soit un générateur de $(\Z/\ell\Z)^\ast$. 
 On conjecture (cf. \cite{bach}) que $g^\ast(\ell) = O(\log \ell (\log \log \ell)^2)$, et bien sûr on conjecture le même résultat pour $g(\ell)$ puisque $g(\ell) \leq g^\ast(\ell)$. On sait par ailleurs ({\it loc. cit.}) que $g(\ell) \neq o (\log \ell \log \log \ell)$.
 
 Le meilleur résultat sous GRH connu semble être le suivant, dû à Shoup:
 \begin{theoreme}[cf. \cite{shoup} or \cite{martin}] \label{thmshoup} Sous GRH, on a
 $g^\ast(\ell) = \bigo {(\omega(\ell-1) \log_1 \omega(\ell-1))^4 (\log \ell)^2}$, où 
 l'on a posé $\log_1 x = \max(1,\log x)$.
 \end{theoreme}
 
Les méthodes de cet article permettent de montrer le résultat suivant:
\begin{prop} \label{primroot} Supposons que $\ell$ soit tel que $\sum_{p \mid \ell-1} 1/p < 1$. Alors
$$g^\ast(\ell) = \bigo { \left( \frac{\omega(\ell-1)}{1-\sum_{p \mid \ell-1} 1/p }\right)^2 (\log \ell)^2}.$$
\end{prop}
\begin{pf} On applique le théorème~\ref{thmpluspetit} à l'extension $L/\Q$ où $L=\Q(\zeta_\ell)$. On a dans ce cas une identification 
canonique $G=(\Z/\ell \Z)^\ast$ qui identifie $\Frob_p$ à $p \mod \ell$. On prend pour $D$ l'ensemble des générateurs de $G$. On a 
$\varphi_G(D) \leq  \bigo{ \frac{1+\omega(\ell-1)}{1 - \sum_{p \mid \ell-1} 1/p }}$ sous notre hypothèse sur $\ell-1$ par le proposition~\ref{phigen}
 et le résultat suit.
\end{pf}

Un principe heuristique bien connu  affirme que les nombres de la forme $\ell-1$ avec $\ell$ premier devrait se comporter statistiquement, en ce qui concerne leur factorisation en facteurs premiers, comme tous les entiers $n$ ayant les mêmes "propriétés locales", i.e. comme tous les entiers pairs. Selon ce principe, il devrait y avoir des nombres premiers $\ell$ avec $\sum_{p \mid \ell-1} 1/p \leq 3/4$ (par exemple) et $\omega(\ell-1)$ arbitrairement grand, puisqu'il est clair
qu'on peut trouver des entiers pair $n$ ayant ces propriétés. Si c'est la cas, la proposition~\ref{primroot} donne pour ces nombres premiers $\ell$ un résultat plus
fort que le théorème~\ref{thmshoup}. 

De plus une analyse de le preuve du théorème~\ref{thmshoup}, par exemple à partir de l'article \cite{martin} qui la généralise quelque peu et en donne une version (nettement moins forte) inconditionnelle,  montre (ou du moins nous semble montrer)
que cette preuve est fondamentalement différente de la notre: si le principe de la preuve de l'existence d'un nombre premier $p$ petit tel que 
$p \mod \ell \in D$, à savoir donner une estimation pour  $\psi_1(f,x)$ et voir pour quel $x$
quand le terne principal, strictement positif, donnee le terme d'erreur, ce principe est appliqué pour la banale fonction caractéristique $f=\un_D$ 
(notée $\gamma$ dans \cite{martin}) et non pour une fonction $f$ cherchant à minimiser le rapport $\lambda_G(f)/\mu_G(f)$ (sous les conditions (a) et (b)
de la définition de $\varphi(D)$). D'un autre côté l'estimation de $\psi_1(\un_D,x)$ utilisée par Martin est bien meilleure que celle qu'on obtient en sommant 
en valeur absolue les termes d'erreurs  des $\psi_1(\chi,x)$ où $\chi$ parcourt les caractères de $(\Z/\ell\Z)$: elle se base sur une application du {\it shifted sieve} 
d'Iwaniec et constitue le point crucial de la preuve.

Ces remarques semblent indiquer qu'une estimation précise de $\varphi(D)$ (nous n'avons utilisé qu'une majoration très grossière et valable sous des hypothèses restrictives) pourrait apporter des résultats nouveaux concernant $g^\ast(p)$, et qu'on pourrait peut-être
la combiner avec la méthode de Martin pour obtenir des résultats encore meilleurs.

  \subsection{Réductions des polynômes à coefficients entiers}
 \label{redpolint}
 \begin{theoreme}
Soit $P$ un polynôme unitaire irréductible à coefficients entiers, de degré $n \geq 1$. Soit $M$ le produit des nombres premiers divisant le discriminant de $P$.
Alors souus (GRH) et la conjecture d'Artin:
\begin{itemize} \item[(a)] Il existe un nombre premier $p  < \ccr{cpp} n^2 (\log M + n \log n)^2 $ ne divisant pas $M$ tel que le polynôme $P(X) \pmod{p}$ admette au moins une racine dans $\F_p$.
\item[(b)] Il existe un nombre premier $p < \ccr{cpp} n^4  (\log M + n \log n)^2$  ne divisant pas $M$ tel que le polynôme $P(X) \pmod{p}$ admette au moins deux racines dans $\F_p$.
\item[(c)] Il existe un nombre premier $p <   \ccr{cpp}n^4 (\log M + n \log n)^2$  ne divisant pas $M$ tel que le polynôme $P(X) \pmod{p}$ n'admette aucine racine dans $\F_p$.
\end{itemize}
\end{theoreme}
\begin{pf} 
Soit $L$ le corps de décomposition de $P$, et $G=\Gal(L/\Q)$ munie de son action naturelle transitive sur l'ensemble des racines de $P$ dans $\Cb$, qu'on identifie à $\{1,\dots,n\}$.
Comme, pour $p$ un nombre premier ne divisant pas $M$, 
le nombre de racines de $P \pmod{p}$ est égal au nombre de points fixés par l'action de $\Frob_p \in G$ sur $\{1,\dots,n\}$, un nombre premier $p$ satisfait la condition du théorème dans le cas (a) (resp. (b), resp. (c))  si et seulement si $\Frob_p$ appartient au sous-ensemble $D$ des éléments de $G$ ayant au moins un point fixe (resp. ayant au moins deux point fixes, resp. n'ayant pas de point fixe.) Par le théorème~\ref{thmpluspetit}
il existe un tel $p$ qui soit $< \ccr{cpp}  \varphi_G(D)^2 (\log(M) + \log |G|)^2)$. Le résultat découle donc de la proposition~\ref{propphisym} et de la majoration triviale $\log |G| \leq \log n! \leq n \log n.$
\end{pf}

\begin{remarque} Dans chacun des cas, la majoration sur $p$ du théorème est polynomiale en $n=\deg P$. Ceci est à comparer avec la majoration qu'on obtiendrait par une application directe du théorème  de Lagarias-Odlyzko, cf. (\ref{ppplo}) (resp. Murty et Murty, cf. (\ref{pppmm})) qui serait en $n (n!)^2 \log n$ 
(resp. $n n! \log n$), donc super-exponentielle.

Une majoration polynomiale dans le cas (a) du théorème ci-dessus avait néanmoins déjà été obtenue sous (GRH) par Weinberger (\cite{weinberger})
et Adleman-Odlyzko (\cite{ao}),  indépendamment. Il sera peut-être intéressant de comparer brièvement leurs preuves (qui sont essentiellement les mêmes) avec la notre (qui est différente, mais possède un point commun avec celles-ci). Reformulée dans le langage de cet article, ces preuves reviennent à étudier a fonction $\pi(f,x)$ (ou une variante comme $\psi_1(f,x)$, peu importe), comptant les nombres premiers $p<x$ avec multiplicité égale à $f(\Frob_{p,L})$, i.e. au nombres de points fixes de $\Frob_{p,L}$ agissant sur $\{1,\dots,n\}$. On retrouve dans ce cas particulier donc la même idée que dans notre preuve du théorème~\ref{thmpluspetit} d'utiliser une fonction auxiliaire $f$ sur $G$ (plutôt que la fonction caractéristique des éléments de $G$ ayant au mons un point fixe). Cependant nous avons dit plus haut (cf. Remarque \ref{remarqueppm}) que cette idée ne menait à rien quand on la combinait à la version du théorème de Chebotarev de Lagarias-Odlyzko, la seule disponible à l'époque. Comment Adleman-Odlyzko et Weinberger concluent-ils donc? En utilisant la simple observation suivante, que $\pi(f,x)$ est égal au nombre
d'idéaux premiers $\p$ du corps {\it de rupture} $K$ du polynôme $P$ de degré $1$ sur $\Q$ avec $N(\p)<x$. Notons que le corps $K$ n'est pas en
général galoisien, mais qu'il est de degré $n$, tandis que le corps de décomposition $L$ de $P$ peut être de degré $n!$. On utilise ensuite une estimation élémentaire
du norme d$\p$ de degré plus grand que $1$ avec $N(\p)<x$ (ils sont absolument négligeables) et le "théorème des nombres premiers" pour $K$ (sous GRH) estimant le nombre de $\p$ avec $N(\p)\leq x$ pour conclure. La démonstration n'utilise donc pas le théorème de Chebotarev effectif, mais seulement
le théorème des nombres premiers pour $K$ (qui est si l'on y tient une forme de Chebotarev effectif, mais pour l'extension triviale $K/K$).

Il ne semble pas qu'on puisse généraliser cette astuce à d'autres situations, comme celles des cas (b) et (c) du théorème ci-dessus.
\end{remarque}

\begin{theoreme}
Soit $P$ un polynôme unitaire irréductible à coefficients entiers, de degré $n \geq 1$. Soit $M$ le produit des nombres premiers divisant le discriminant de $P$.
Supposons (GRH) et Artin. S'il existe un nombre premier $p \nmid M$  tel que le polynôme $P \pmod{p}$ soit irréductible dans $\F_p[X]$, alors il en existe un
qui soit $< \ccr{cpp} 4^{n-1}\frac{ (\log M + n \log n)^2}{n^4} $.
\end{theoreme}
\begin{pf} On garde les notations de la preuve du théorème précédent. Pour $p \nmid M$, $\Frob_{p,G}$ est un $n$-cycle si et seulement su $P \pmod{p}$ est irréductible. L'hypothèse implique donc que l'ensemble $D$ des $n$-cycles dans $G$ est non vide. Comme le centralisateur d'un $n$-cycle dans $S_n$ est le sous-groupe d'ordre $n$ qu'il engendre, il en est de même de son centralisateur dans $G$, et donc la classe de conjugaison d'un $n$-cycle est de cardinal $|G|/n$. On a donc $|D| \geq |G|/n$. Par ailleurs $\lambda_G(D) \leq 2^{n-1}/n$ d'après la proposition~\ref{phisym1} et la proposition~\ref{compres}. 
On a donc $\varphi_G(D) \geq 2^{n-1}/n^2$ et le résultat en découle par le théorème~\ref{thmpluspetit}.
\end{pf}
\subsection{Systèmes compatibles et conjecture de Lang-Trotter généralisée}
\label{syscompLT}
\subsubsection{Position du problème}

On aimerait partir d'un (iso-) motif pur $E$ défini sur $\Q$. Comme, faute d'avoir prouvé les conjectures standard et la conjecture de Hodge, on ne dispose pas encore d'une catégorie des motifs (même en se restreignant aux sommes directes de motifs purs) ayant toutes les propriétés qu'on attend d'elle,
nous partirons directement du système compatible de représentations galoisiennes attachées au motif pur $E$ sans faire référence à $E$, en supposant explicitement les propriétés qu'on attend d'un tel système attachée à un motif.

On se donne donc un {\it système compatible de représentations galoisiennes}, pur, à savoir:
\begin{itemize}
\item[(a)] un entier $n \geq 1$ (la {\it dimension})
\item[(b)] un entier sans facteurs carrés $M\geq 1$ (le {\it discriminant réduit)};
\item[(c)] Pour tout $\ell$, une représentation continue $\rho_\ell : G_{\Q} \rightarrow \Gl_n(\Q_\ell)$ (les {\it représentations galoisiennes} )
\end{itemize}
tels que
\begin{itemize}
\item[(i)] Pour tout $\ell$, $\rho_\ell$ est non ramifiée aux nombres premiers $p$ ne divisant par $M\ell$.
\item[(ii)] Pour tout $p$ ne divisant pas $M$, il existe un polynôme $\Xi_p(X) \in \Z[X],\ \Xi_p(X)=X^n + (-1)^n a_p X^{n-1} + \dots $, 
tel que pour tout $\ell \neq p$, le polynôme caractéristique de $\rho_\ell(\Frob_p)$ soit $\Xi_p$.
\item[(iii)] il existe un entier $w < 0$ tel que pour tout $p$ ne divisant pas $M$,
 les racines de $\Xi_p(x)$ dans $\Cb$ soit toutes de module $p^{- w/2}$.
\end{itemize}
Il est clair que les polynômes $\Xi_p(X)$ et l'entier $w$ sont déterminés uniquement par les données (a), (b), (c). L'entier $w$ s'appèle le 
poids du système compatible.  L'application qui à tout premier $p$ ne divisant pas $M$ associe l'entier $a_p$ (cf. (ii)) 
tel que pour tout $\ell \neq p, a_p = \tr \rho_\ell(\Frob_p)$ jouera un rôle essentiel. Par Chebotarev, elle détermine le système compatible à un isomorphisme près.

Nous supposerons également que le système compatible satisfait la propriété suivante, 
qui est conjecturée pour tout système compatible attachée à un  motif pur 
\begin{itemize}
\item[(iv)] Il existe un entier $N \geq 1$, et un sous-schéma en groupes réductif $G$ de $\Gl_n$ défini sur $\Z[1/N]$, contenant le centre de $\Gl_n$,
 tel que pour tout $\ell$ ne divisant pas $N$, $\rho_\ell (\GQ)$ est contenu dans $G(\Z_\ell)$, Zariski-dense dans $G(\Q_\ell)$, et 
 est un sous-groupe ouvert d'indice borné indépendamment de $\ell$ de $G(\Z_\ell)$.
\end{itemize}
L'entier $N$ et le schéma en groupes $G$ sont essentiellement bien déterminé: si $N'$ est un autre entier $\geq 1$ et $G'$ un autre sous-schéma en groupes
satisfaisant la condition (iv), il existe un $N''$ divisible par $N$ et $N'$ tels que $G \simeq G'$ sur $\spec \Z[1/N'']$. On fixe dorénavant un tel $N$ et $G$,
et on notera $G_\Q$ la fibre générique de $G$ et $G_{\F_\ell}$ pour $\ell \nmid N$ ses fibres spéciales. La dimension des groupes $G_{\F_\ell}$ et de $G_\Q$ est constante, de même que leur rang (dimension d'un tore maximale, non nécessairement déployé), et on les notera $d$ et $r$ respectivement.

La propriété (iv) est une combinaison de plusieurs conjectures faites dans \cite{serremotives} (et qui pour la plupart, faisaient déjà partie du folklore: voir l'introduction de \cite{serremotives}). Plus précisément, l'existence du groupe $G_\Q$ tel que pour tout $\ell$, $\rho_\ell(\GQ)$ est Zarisiki-dense dans $G_\Q(\Q_\ell)$  est 
\cite[3.2?]{serremotives}, le fait que $G_\Q$ contiennent le centre de $\Gl_n$ résulte de l'existence de l'homomorphisme $\bf{w: G}_m \rightarrow G_\Q$
associé à la gradation par le poids (cf. \cite[\S5]{serremotives}) et de notre hypothèse $w < 0$ (en fait, ce qui compte est $w\neq 0$), l'existence de la forme entière $G$ résulte de la théorie des $\Z$-formes exposée
au \cite[\S10]{serremotives} et la dernière assertion est \cite[10.3?]{serremotives}.

\begin{exemple}
\begin{itemize}
\item Si $f = \sum a_n q^n \in \Z[q]]$ est une forme modulaire parabolique de niveau $\Gamma_0(M)$, poids $k \geq 1$, et propre pour les opérateurs de Hecke $T_p$ avec $p \nmid M$, alors il existe d'après Eichler-Shimura (si $k=2$) et Deligne (si $k>2$) un système compatible attaché à $f$, de dimension $n=2$, de poids $w=1-k$, de déterminant réduit le radical de $M$ et dont les $a_p$ sont les coefficients $a_p$ de $f$. Un tel système satisfait toujours  (iv) et le groupe $G$ est $\Gl_2$ si $f$ n'a pas de multiplication complexe. 
\item Si $A$ est une variété abélienne sur $\Q$ ayant bonne réduction en dehors de $M$ (un entier sans facteurs carrés), la famille des représentations sur les modules de Tate de $A$ est un système compatible
associé à $A$ de dimension $2g$, poids $w=-1$, discriminant réduit $M$. On sait que ce ce système satisfait (iv) pour un groupe $G$ tel que $G_\Q$ est le groupe de Mumford-Tate de $A$.
Dans le cas particulier où $A$ est la jacobienne d'une courbe propre lisse et géométriquement connexe $C/\Q$ de genre $G$, les $a_p$ du système compatible ont l'interprétation arithmétique bien connue depuis Weil: $a_p = |C(\F_p)| - 1 - p$.
\end{itemize}
\end{exemple}

Fixons donc un système de représentations galoisiennes comme ci-dessus. Fixons par ailleurs un entier $a \in \Z$. Soit 
$$\pi(a,x) = |\{p \, | \, p<x \text{ et } a_p=a\}|.$$ 
Quand le système compatible est associée à une forme modulaire $f$ ou une variété abélienne $A$, on écrira souvent $\pi_f(a,x)$ ou $\pi_A(a,x)$.

Déterminer l'ordre de grandeur de $\pi(a,x)$ quand $x$ tend vers l'infini est une vaste question ouverte généralisant la conjecture 
de Lang-Trotter, dont on ne connait même pas de réponse conjecturale en général (cf. la discussion de \cite[page 423 et suivantes]{katz}). 
Rappelons seulement la majoration conjecturale bien connue suivante:
\begin{conjecture}[cf. \cite{lt}, \cite{serre}] \label{langtrotter1}
Si le système compatible $(\rho_\ell)$ est attaché à une forme modulaire $f=\sum a_n q^n$ à coefficients entiers, non CM, propre pour les opérateurs de Hecke, de poids $k$, on a, quelque soit $a \in \Z$, 
$$ \pi(a,x) = 
\begin{cases} O(x^{1/2}/\log x) &\text{si } w=-1 \\
		     O(\log\log x) & \text{si }w=-2 \\
		     O(1) & \text{si } w\leq -3 \\
\end{cases} $$
Rappelons que le poids $w$ du système compatible attaché à $f$ est $w=1-k$.
\end{conjecture}
\begin{question} Est-il raisonnable de généraliser la conjecture~\ref{langtrotter1}  à tout système compatible tel que $G$ est connexe? \end{question}
En tout cas, Katz énonce une conjecture qui implique celle-ci dans le cas où le système est attachée à une courbe propre et lisse de genre $g$ de groupe de Mumford-Tate $\GSP_{2g}$.

Quoiqu'il en soit, nous sommes très loin de démontrer ces conjectures. Le but de ce \S\  est seulement de prouver en supposant vraie GRH et la conjecture d'Artin, et sous certaines hypothèses sur le système compatible, des majorations de $\pi(a,x)$ qui sont meilleures que celles obtenues jusqu'ici dans la littérature.  

\subsubsection{Majoration en fonction de la complexité de Littlewood des matrices de trace $a$ dans $G(\F_\ell)$ -- sans le crible}

Notre premier résultat est une application directe du théorème~\ref{denszero}:
\begin{theoreme} \label{syscomp1} Soit $(\rho_\ell)$ un système compatible comme ci-dessus, vérifiant les propriétés {\rm(i)} à {\rm (iv)}. 
On suppose de plus que le groupe $G$ est connexe, et que sa dimension $d$ est $\geq 1$.
Soit $\beta$ un nombre réel tel qu'il existe $R>0$ tel que pour tout $\ell \nmid N$, 
$$\lambda_{G(\F_\ell)}(G(\F_\ell)^\ta)  \leq R |G(\F_\ell)|^{\beta}.$$
Alors, sous la conjecture de Riemann et la conjecture d'Artin, on a
$$\pi(a,x) = \bigo{ x^{\frac{2 \beta d +1}{2 \beta d+2}} (\log x)^\frac{1-\beta d}{1+\beta d}}.$$
\end{theoreme}
\begin{pf}
Pour appliquer le théorème~\ref{denszero}, il nous faut définir un ensemble d'indice $\Lambda$, des groupes de Galois
$G_\nu$ pour $\nu \in \Lambda$, des sous-ensembles invariants par conjugaisons $D_\nu$ de $G_\nu$, et vérifier les hypothèses (a), (b) et (c)
de ce théorème. On prend pour $\Lambda$ l'ensemble des nombre premiers ne divisant pas $N$, et pour $\ell \in \Lambda$, on définit $G_\ell$ comme l'image de
$\rho_\ell(\GQ)$ par l'application $G(\Z_\ell) \rightarrow G(\F_\ell)$. Comme cette application est surjective pour presque tout $\ell$, il 
résulte de l'hypothèse (iv) que $G_\ell$ est d'indice borné dans $G(\F_\ell)$. On prend pour $D_\ell$ l'ensemble $G_\ell^\ta = G(\F_\ell)^{\ta} \cap G_\ell$, ou ce qui revient au même, $G_\ell^\ta$ est l'ensemble des matrices de trace $a \pmod{\ell}$ dans $G_\ell \subset G(\F_\ell) \subset \Gl_n(\F_\ell)$.
L'ensemble $\tilde D$ associée à ces $D_\ell = G_\ell^\ta$ comme en (\ref{tildeDdef}) est alors l'ensemble des $p \nmid M$ tels que $a_p \equiv a \pmod{\ell}$ pour presque tout $\ell$, i.e. l'ensemble des $p \nmid M$ tels que $a_p=a$. En d'autres termes, la quantité $\pi(D,x)$ estimée dans le théorème~\ref{denszero}
n'est autre que $\pi(a,x)$.

Vérifions les hypothèses (a), (b), et (c) du théorème~\ref{denszero}. Puisque $G_\Q$ est de dimension $d$, et $G$ de type fini sur $\Z[1/N]$, 
$|G(\F_\ell)| \sim \ell^d$ quand $\ell$ tend vers l'infini. Comme $G_\ell$ est d'indice borné dans $G(\F_\ell)$, il suit que
$\log |G_\ell| \sim d \log \ell$ quand $\ell$ tend vers l'infini, ce qui implique (a). Pour prouver (b) avec $\alpha=1/d$, il suffit de
voir que le sous-schéma fermé de $G$ définie par l'équation $\tr = a$ a sa fibre générique de dimension (au plus) $d-1$. Comme $G_\Q$ est connexe, donc géométriquement irréductible, il suffit de voir que $G(\Cb)^{\ta} \neq G(\Cb)$. Mais $G(\Cb)$ contient au moins une matrice de trace $n>0$ (l'identité)
donc une matrice de toute trace non nulle (puisque $G$ contient le centre de $\Gl_n$), et donc il est clair que quel que soit la valeur de $a$, on a bien
$G(\Cb)^{\ta} \neq G(\Cb)$, d'où (b). Quant à l'hypothèse (c), elle est vérifiée avec le même $\beta$ que dans l'énoncé du théorème,
puisque $\lambda_{G_\ell}(G_\ell^\ta) \leq \lambda_{G(\F_\ell)}(G(\F_\ell)^{\ta})$ par la proposition~\ref{compres}.

Le théorème résulte maintenant de la formule (\ref{formdenszero}) du théorème~\ref{denszero} en y faisant $\alpha = 1/d$.
\end{pf}
\begin{remarque} On remarquera que sans l'hypothèse $G$ connexe, l'hypothèse (b) du théorème~\ref{denszero} peut être mise en défaut, 
en prenant par exemple pour $G$ le normalisateur d'un tore dans $\Gl_2$ et $a=0$.
\end{remarque}
\begin{cor} \label{corsccs}
 Soit $(\rho_\ell)$ un système compatible comme ci-dessus, vérifiant les propriétés {\rm (i)} à {\rm (iv)}. On suppose de plus que le groupe $G$ est connexe, et de dimension $d \geq 1$. Alors, sous (GRH) et la conjecture d'Artin, on a
 $$\pi(a,x) = \begin{cases} \bigo{ x^{\frac{d}{d+1}} (\log x)^{\frac{3-d}{d+1}} }& \text{ si }a\neq 0 \\
 \bigo{ x^{\frac{d-1}{d}} (\log x)^{\frac{4-d}{d}} }& \text{ si }a =  0 \end{cases}
 $$
 \end{cor}
 \begin{pf}
En effet, on peut prendre dans le théorème précédent $\beta = (1-1/d)/2$ d'après la majoration de Cauchy-Schwarz si $a \neq 0$, $\beta = (1-2/d)/d$ 
si $a=0$ en utilisant la majoration de Cauchy-Schwarz après avoir quotienté par les homothéties. 
\end{pf}
Ce corollaire est prouvé dans \cite{effective} dans le cas d'un système compatible attachée à une forme modulaire parabolique propre pour les opérateurs de Hecke non CM à coefficients entiers de poids $k \geq 2$ (dans ce cas on a $d=4$: cf. \cite[page 254]{effective}).
Bien qu'il ne soit pas énoncé dans cette généralité dans \cite{effective}, il résulte facilement des méthodes employées dans cet article.
Pour que le théorème~\ref{syscomp1} donne un meilleur résultat que son corollaire, c'est-à-dire pour qu'il donne un résultat nouveau par rapport à \cite{effective},
il faut disposer d'une estimation de la complexité de Littlewood  $\lambda_{G(\F_\ell)}(G(\F_\ell)^\ta)$ asymptotiquement meilleure 
que celle donnée par Cauchy-Schwarz. Malheureusement, nous ne savons prouver de telles estimations dans aucun cas.
Dans le seul cas où nous pouvons estimer $\lambda(G(\F_\ell)^\ta)$, i.e.  pour $G=\Gl_n$, la majoration de Cauchy-Schwarz est du bon ordre de grandeur d'après  le théorème~\ref{thmlambdagln}. 

\subsubsection{Majoration en fonction de la complexité de Littlewood des matrices de trace $a$ dans $T(\F_\ell)$ -- sans le crible}

Le but de cette partie est de montrer sur un exemple (celui des courbes de genres 2) 
qu'on peut obtenir des résultats meilleurs que ceux du théorème précédent quand l'on dispose de
réponses à la question~\ref{questiontore}, i.e. d'estimations meilleures que celles de Cauchy-Schwarz de la complexité de Littlewood
de certains sous-ensembles de tores.

\begin{theoreme} \label{thmcompreg}
Soit $(\rho_\ell)$ un système compatible comme ci-dessus, vérifiant les propriétés {\rm(i)} à {\rm (iv)}. 
On pose $\pi^\reg(a,x) = |\{p<x, p \nmid M, a_p=a, \Xi_p(x) \text{ est à racines simples dans $\Cb$} \} |$.
On suppose de plus que le groupe $G$ est connexe, et de dimension $d \geq 1$.  Soit $T$ un tore maximal de $G$ défini sur $\Z[1/N]$. Soit $0 \leq \gamma < 1$ un réel tel qu'il existe $R>0$ tel que pour tout $\ell \nmid N$, on ait $$\lambda_{T(\F_\ell)}(T(\F_\ell)^{\ta,\dr}) < R |T(\F_\ell)|^\gamma.$$
 Alors, sous (GRH) et la conjecture d'Artin, on a
\begin{eqnarray} \label{formcompreg} \pi^\reg(a,x) = \bigo{ x^{\frac{2 \beta d +1}{2 \beta d+2}} (\log x)^\frac{1-\beta d}{1+\beta d}} \text{ 
où l'on a posé } \beta =\frac{ 2 \gamma r + d-r}{2d}.\end{eqnarray}
\end{theoreme}

\begin{pf}
Pour $\ell$ premier,  on pose 
$$\pi(a,\ell,x) = |\{p \leq x, \ p \nmid N,\  a_p=a,\text{ et  le polynôme $\Xi_p(X)$ est scindé à racines simples modulo $\ell$}\}|.$$
Le lemme crucial, inspiré de \cite[lemma 4.4]{effective} mais utilisé ici dans un but différent, est le suivant:

\begin{lemme} Il existe deux constantes $\ccl{clemme}>0$ et $\ccl{clemme2}>0$, dépendant seulement du système compatible considéré, telles que, pour tout $x >2$, $y > \ccr{clemme2}$, et tout  $u$ tel que $y^{1/2} (\log xy) (\log y)^2 \leq u \leq y$, on ait  sous GRH,
$$\pi^\reg(a,x) \leq \ccr{clemme} \max_{\ell \in I} \pi(a,\ell,x),$$
où le maximum est pris sur les nombres premiers $\ell$ appartenant à l'intervalle $I=[y,y+u]$.
\end{lemme}
\begin{pf}
Dans toute cette preuve les constantes $\ccr{cpreuve},$ etc. dépendront du système compatible fixé (et donc de sa dimension $n$ et son poids $w$), 
mais pas de $x$, $y$, ou $u$.

Posons \begin{eqnarray*} \pi_p(I) &=& | \{ \ell \in I, \Xi_p(X) \text{ est scindé à racines simples modulo $\ell$}\}| \\
\pi_p(y) &=& | \{ \ell \leq y, \Xi_p(X) \text{ est scindé à racines simples modulo $\ell$}\}| \end{eqnarray*}
Fixons $x \geq 3$. Soit $p$ un nombre premier $<x$ ne divisant pas $M$, et tel que $\Xi_p$ est à racines simples dans $\Cb$. 
D'après l'hypothèse de pureté (iii) des systèmes compatibles, les racines de $\Xi_p(X)$
sont de modules complexes au plus $x^{-w/2}$, donc leurs différences
sont de modules au plus $2x^{-w/2}$, et le discriminant $d_{\Xi_p}$ du polynôme $\Xi_p(X)$ est donc majoré par $(2x)^{-n(n-1)w/2}$. Comme ce discriminant est non nul, le produit $M_p$ des nombres premiers $\ell$ qui divisent $d_{\Xi_p}$ est tel que $\log M_p \leq \ccl{cpreuve} \log x$.
Soit $L_p$ le corps de décomposition de $\Xi_p(x)$, $G_p$ son groupe de Galois. Tout nombre premier ramifié dans $L_p$ divise $M_p$, 
et un nombre premier $\ell$ est tel que $\Xi_p(X)$ est scindé à racines simples modulo $\ell$ si et seulement si $\ell$ ne divise pas $M_p$ et $\Frob_{\ell,L_p} = 1$.  Par le théorème de Chebotarev effectif de Lagarias-Odlyzko-Serre sous GRH (\ref{chefelosf}), on a donc pour $y \geq 2$
$$|\pi_p(y) - \frac{1}{|G_p|} \Li(y) |  \leq  \ccr{ctlos} y^{1/2} (\log y + \log |G_p| + \log M_p) $$
et donc, comme $\log |G_p| \leq \log n!$, 
$$|\pi_p(y) -  \frac{1}{|G_p|}  \Li(y) | \leq \cc y^{1/2} \log(xy)$$
et de même, pour $u \leq y$
$$|\pi_p(y+u)  - \frac{1}{|G_p|} \Li(y+u) | \leq \cc y^{1/2} \log(xy).$$
 Si $\pi(I)$ est le nombre de nombres premiers dans l'intervalle $I=[y,y+u]$, le théorème des nombres premiers (sous sa forme
 usuelle sous l'hypothèse de Riemann) donne donc pour $y$ assez grand et $u$  comme dans  dans l'énoncé
$$\pi(I) > \ccr{clemme} \pi_p(I)$$

On a 
\begin{eqnarray*} \max_{\ell \in I} \pi(a,\ell,x) &\geq & \frac{1}{\pi(I)} \sum_{\ell \in I} \pi(a,\ell,x)\\
&=& \frac{1}{\pi(I)} \sum_{p \leq x, a_p = a} \pi_p(I), \text{ \ \ \ \ d'après les définitions} \\
& > & \frac{1}{\pi(I)}  \sum_{\substack{p \leq x, a_p = a \\ \Xi_p(X) \text{ à racines simples dans } \Cb}} \ccr{clemme} \pi(I) \\
& = & \ccr{clemme}  \pi^\reg(a,x) \end{eqnarray*}
et le lemme suit.
\end{pf}

Revenons à la preuve du théorème~\ref{thmcompreg}.
Posons $G_\ell^{\ta,\dr} = G(\F_\ell)^{\ta,\dr} \cap G_\ell$.
\begin{eqnarray*} \pi(a,\ell,x) \leq \pi(G_\ell^{\ta,\dr},x). \end{eqnarray*}
En effet le membre de gauche compte les nombres premiers $p<x$, $p \nmid M$, avec $a_p=a$ et $\Xi_p(x)$ scindé à racines simples modulo $\ell$, donc tels que 
$\Frob_{p,G_\ell}$ est  diagonalisable régulier dans $\Gl_n(\F_\ell)$ tandis que le terme de droite compte les nombres premiers satisfaisant le même condition avec $a_p=a$ remplacée par $a_p \equiv a \pmod{\ell}$.

La preuve est alors basée sur le même principe que celle du théorème~\ref{denszero}.
Posons $\alpha=1/d$, et $\beta = \frac{ 2 \gamma r + d-r}{2d} $ afin que les formules qui suivent ressemblent à celles de cette preuve.
On prend $y = y(x) =  x^{\alpha /(2\alpha+2\beta)} (\log x)^{-2\alpha/(\alpha+\beta)}  $ et $u=y$.
Pour $x$ assez grand, $u$ satisfait bien les hypothèses du lemme, et il existe un $\ell=\ell(x) \in [y, 2y]$ tel que
$\pi(a,x) < c' \pi(a,\ell,x) \leq c' \pi(G_\ell^{\ta,\dr},x)$. On supposera dorénavant que $x$ est assez grand pour que $\pi^\reg(a,x) > 0$ (c'est possible sauf si $\pi^\reg(a,x)$ est identiquement nul, auquel cas il n'y a rien à prouver), ce qui implique $\pi(a,\ell,x) > 0$ et donc que $T_{\F_\ell}$ est un tore déployé.

On a alors $|G_{\ell}| \asymp |G(\F_\ell)| \asymp x^{1/(2\alpha+2\beta)} (\log x)^{-2/(\alpha+\beta)}$ quand $x$ tend vers l'infini
(avec $u=y=y(x)$ comme ci-dessus et $\ell$ arbitraire entre $y$ et $2y$) d'après l'hypothèse (iv) des systèmes compatibles. 

Le théorème de Chebotarev~\ref{thmchefD} donne, en admettant (GRH) et la conjecture d'Artin pour les fonctions $L$ d'Artin des représentations $\pi$ dans
 le support spectral de $G_\ell^{\ta,\dr}$:
\begin{eqnarray} \label{tchsr} \pi(G_\ell^{\ta,\dr},x) = \frac{|G_\ell^{\ta,\dr}|}{|G_\ell|} \Li(x) + \bigo{ x^{1/2} (\log x) \lambda_{G_\ell}(G_\ell^{\ta,\dr})}.\end{eqnarray}

d'où
$$\pi^\reg(a,x) = \bigo{  \frac{|G_\ell^{\ta,\dr}|}{|G_\ell|} \Li(x) } +  \bigo{ x^{1/2} (\log x) \lambda_{G_\ell}(G_\ell^{\ta,\dr})}$$
et pour prouver le  théorème il suffit de voir que chacun des deux termes $\bigo{\dots}$ dans la formule ci-dessous est
$ \bigo{ x^{\frac{2 \beta d +1}{2 \beta d+2}} (\log x)^\frac{1-\beta d}{1+\beta d}}.$
Pour le premier, cela résulte simplement de ce que $G_\ell^{\ta,\dr} \subset G_\ell^\ta$ et de la majoration $|G_\ell^\ta| = O(\ell^{d-1})$ obtenue
dans la preuve du théorème~\ref{syscomp1}.
Pour le second, on observe que quand $x$ tend vers l'infini, et $\ell=\ell(x)$ avec lui,
\begin{eqnarray*} \lambda_{G_\ell}(G_\ell^{\ta,\dr}) &\leq& \lambda_{G(\F_\ell)}(G(\F_\ell)^{\ta,\dr}) \\
& \asymp & \ell^{(d-r)/2}\lambda_{T(\F_\ell)}(T(\F_\ell)^{\ta,\dr}) \text{ par le corollaire~\ref{corGT}} \\
& = & \bigo{ \ell^{(d-r)/2} \ell^{\gamma r} } \\
&=& \bigo{ |G(\F_\ell)|^{ \gamma r / d + (d-r)/2d} }\\
&=& \bigo{ |G(\F_\ell)|^\beta}
\end{eqnarray*} 
si bien qu'avec l'estimation de $|G(\F_\ell)|$ donnée plus haut, on voit facilement, que le terme
$$x^{1/2} (\log x)  \lambda_{G(\F_\ell)}(G(\F_\ell)^{\ta,\dr})$$ est $ \bigo{ x^{\frac{2 \beta d +1}{2 \beta d+2}} (\log x)^\frac{1-\beta d}{1+\beta d}}$,
ce qui termine la preuve.
\end{pf}

\begin{remarque} Il n'est pas très difficile de voir que si l'on suppose $\rho_\ell$ surjectif pour $\ell$ assez grand, le théorème
précédent est valable sans supposer vraie conjecture d'Artin. En effet celle-ci n'est utilisée qu'une seule
fois dans le preuve ci-dessus, pour appliquer le théorème de Chebotarev effectif (\ref{tchsr}) à l'extension $L_\ell/\Q$ de groupe de Galois $G_\ell$,
pour le sous-ensemble $G_\ell^{\ta,\dr}$.
Si $\rho_\ell$ est surjectif, et si $G(\Z_\ell) \rightarrow G(\F_\ell)$ l'est aussi, donc pour tout $\ell$ assez grand, 
on a $G_{\ell} = G(\F_\ell)$ et donc $G_\ell^{\ta,\dr} = G(\F_\ell)^{\ta,\dr}$. Il s'agit donc de voir que la formule (\ref{tchsr}) est vraie sans 
avoir recours à la conjecture d'Artin. Mais la preuve du théorème~\ref{tits} montre que $\un_{G(\F_\ell)^{\ta,\dr}}$ peut s'écrire $\sum c_\chi \tr \Ind_{B(\F_\ell)}^{G(\F_\ell)} \chi$, où $B(\F_\ell)$ est un Borel de $G(\F_\ell)$ (on rappelle que pour les $\ell$ considérés ici, $T(\F_\ell)$ est un tore maximal déployé),
et les $\chi$ des caractères abéliens de $B(\F_\ell)$, et qu'on a $\lambda(G(\F_\ell)^{\ta,\dr})=[G(\F_\ell):B(\F_\ell)] \sum_{\chi} |c_\chi|$.
Pour prouver la formule (\ref{tchsr}), il suffit donc par linéarité de prouver la formule analogue où $\un_{G(\F_\ell)^{\ta,\dr}}$ est remplacé par
la fonction $ \tr \Ind_{B(\F_\ell)}^{G(\F_\ell)} \chi$ pour chacund es caractères $\chi$. 
Mais la fonction $L$ attachée à la représentation $\Ind_{B(\F_\ell)}^{G(\F_\ell)} \chi$
satisfait la conjecture d'Artin, puisque celle-ci est connue pour les caractères abélien et stable par induction, et donc la formule voulue est vraie sous GRH seule.
\end{remarque}

\begin{theoreme} \label{thmgenre2} Soit $A$ une variété abélienne de dimension $2$ sur $\Q$ telle que $\End_{\bar \Q}(A)=\Z$.
Alors sous (GRH) et Artin on a
$$\pi_A(0,x) = \bigo{x^{9/10} (\log x)^{-3/5}}.$$
\end{theoreme}
\begin{pf} Soit $(\rho_\ell)_{\ell \text{ premier}}$ le système de représentations galoisiennes attachées à $C$. Il est bien connu depuis Weil que $(\rho_\ell)$ est un système compatible, au sens où il satisfait les conditions (i), (ii) et (iii). Sous l'hypothèse $\End(A)=\Z$, un résultat de Serre (\cite{serrevigneras}) implique
que ce système satisfait également (iv) pour $G=\GSP_4$. Ce groupe est de dimension $d=11$, et de rang $r=3$.

On écrit $\pi_C(0,x)=\pi_C^\reg(0,x) + \pi_C^\nreg(0,x)$,
où le premier terme (resp. le second terme) compte les $p<x$, $p \nmid N$ tels que $a_p=0$ et $\Xi_p(X)$ n'a pas de racine multiple dans $\Cb$ (resp. a au moins une racine multiple).

Pour le premier terme, on applique le théorème~\ref{thmcompreg} en notant qu'on peut prendre $\gamma=0$ dans l'énoncé de ce théorème grâce au théorème~\ref{toresymplectique}, donc $\beta=4/11$ et $2 \beta d = 8$, si bien que 
$\pi_C^\reg(0,x) = \bigo{x^{9/10} (\log x)^{-3/5}}.$

Pour le second, on applique le théorème~\ref{syscomp1} en prenant $D_\ell \subset G(\F_\ell)$ l'ensemble des éléments non réguliers de trace $0$.
Comme $|D_\ell/Z| = O(\ell^8)$, on peut prendre par Cauchy-Schwarz $\beta=4/11$, et on obtient la même estmitation pour $\pi_C^\nreg$.
\end{pf}

\begin{remarque} En appliquant le corollaire~\ref{corsccs}, on obtiendrait dans ce cas $\pi(0,x)=\bigo{x^{10/11} (\log x)^{-7/11}}$.
Nous verrons plus bas qu'on peut encore améliorer ce résultat en utilisant la méthode du crible.
\end{remarque}

\subsubsection{Majoration en fonction de la complexité de Littlewood des matrices de trace $a$ dans $G(\F_\ell)$ -- avec le crible}

Ajoutons une hypothèse d'indépendance pour notre système compatible. Pour $\ell \nmid N$, on note comme auparavant $G_\ell$ l'image de
$\GQ$ dans $G(\F_\ell)$.
\begin{itemize} 
\item[(v)] Il existe $\ell_0$ tel que l'application $\GQ \rightarrow \prod_{\ell > \ell_0} G_\ell$ est surjective.
\end{itemize}
Il est conjecturé (cf. \cite{serremotives}) que si l'on replace $\GQ$ par $\Gal(\bar \Q/K)$ (deux fois: dans la définition de $G_\ell$ et dans la propriété (v))
 ou $K$ est une extension finie de $\Q$ assez grande. On sait par ailleurs que la propriété (v) est vraie (telle qu'énoncée, i.e. pour $K=\Q$) 
par exemple pour les systèmes compatibles attachés aux courbes elliptiques sur $\Q$ sans multiplication complexe.

\begin{theoreme} Soit $(\rho_\ell)$ un système compatible comme ci-dessus, vérifiant les propriétés {\rm(i)} à {\rm (v)}. 
On suppose de plus que le groupe $G$ est connexe, et que sa dimension $d$ est $\geq 1$, et que la représentation $\rho_\ell: \GQ \rightarrow G(\Z_\ell)$
est surjective pour $\ell$ suffisamment grand.  Soit $\beta$ un nombre réel tel qu'il existe $R>0$ tel que pour tout $\ell \nmid N$, 
$$\lambda_{G(\F_\ell)}(G(\F_\ell)^\ta)  \leq R |G(\F_\ell)|^{\beta}.$$
Alors sous GRH et la conjecture d'Artin, on a pour tout $\epsilon>0$:
$$\pi_D(x)=\bigo{ x^{\frac{d\beta+1}{d\beta+2}+\epsilon}} $$
\end{theoreme}
\begin{pf} On applique le théorème~\ref{denszerocrible}, avec $G_\ell$ et $D_\ell$ définie comme dans la preuve du théorème~\ref{syscomp1}. Grâce à 
l'hypothèse de surjectivité que nous faisons ici, on sait que $G_\ell = G(\F_\ell)$ pour $\ell$ suffisamment grand, disons $\ell \nmid N$ quitte à remplacer $N$ par un entier plus grand, et  $D_\ell = G(\F_\ell)^{\ta}$.
Ainsi $|D_\ell|$, resp. $|G_\ell|$ est l'ensemble des points sur $\F_\ell$ d'un schéma de type fini sur $\Z[1/N]$ dont la dimension de la fibre générique est $d-1$ (resp. $d$). L'hypothèse (b') du théorème~\ref{denszerocrible} avec $\alpha=1/d$ en résulte aisément. L'hypothèse (c') est vraie per hypothèse. L'hypothèse (d')
du théorème~\ref{denszerocrible} est conséquence de l'hypothèse (v) pour notre système compatible.
Le résultat suit.
\end{pf}
 On a le corollaire:
\begin{cor} \label{corsyscomp3} Soit $(\rho_\ell)$ un système compatible comme ci-dessus, vérifiant les propriétés {\rm(i)} à {\rm (v)}. 
On suppose de plus que le groupe $G$ est connexe, et que sa dimension $d$ est $\geq 1$, et que la représentation $\rho_\ell: \GQ \rightarrow G(\Z_\ell)$
est surjective pour $\ell$ suffisamment grand. Alors sous (GRH) et la conjecture d'Artin, on a
 $$\pi(a,x) = \begin{cases} \bigo{ x^{\frac{d+1}{d+3}+\epsilon}} & \text{ if }a \neq 0 \\
\bigo{x^{\frac{d}{d+2} + \epsilon}} & \text{ if }a=0. \end{cases}$$
\end{cor}
La preuve est la même que celle du corollaire~\ref{corsccs}.
On obtient en particulier:
\begin{cor}
Soit $A/\Q$ une variété abélienne de dimension $g$ sur $\Q$, ayant bonne réduction en dehors des premiers divisant un entier $M \geq 1$. On suppose que le système de représentations attachée à $A$ satisfait la condition (v) et que  $\rho_\ell: \GQ \rightarrow \GSP_{2g}(\Z_\ell)$  est surjective pour $\ell$ suffisamment grand.
Alors sous (GRH) et la conjecture d'Artin, on a
$$\pi(a,x) = \begin{cases} \bigo{ x^{\frac{2g^2+g+2}{2g^2+g+4}+\epsilon}} & \text{ if }a \neq 0 \\
\bigo{x^{\frac{2g^2+g+1}{2g^2+g+3} + \epsilon}} & \text{ if }a = 0. \end{cases}$$
\end{cor}
La condition de surjectivité de $\rho_\ell$ est satisfaite si $g$ est  impair ou égal à $2$ ou $6$, et 
$\End_{\Q}(A)=\Z$ (cf. \cite{serrevigneras}).

Dans le cas d'une courbe elliptique $E$ sans multiplication complexe, le théorème donne $\pi(a,x) = \bigo{x^{5/7+\epsilon}}$ si $a \neq 0$, 
$\pi(0,x)=\bigo{x^{2/3+\epsilon}}$. Dans les deux cas ce sont les meilleures bornes (conditionnelles) obtenues jusqu'à présent en direction de la conjecture de Lang-Trotter pour les courbes elliptiques (voir l'introduction pour plus de détail). 

Dans la cas où $A$ est de dimension $2$, et pour $a=0$, le résultat qu'on obtient ici est meilleur que celui obtenu au théorème~\ref{thmgenre2}
mais avec une hypothèse d'indépendance en plus. Rappelons que les méthodes utilisées pour démontrer ces deux résultats (qui tous deux améliorent le résultat qu'on obtiendrait avec la méthode de Murty-Murty-Saradha) sont différentes, la première reposant sur un calcul de complexité de Littlewood pour un tore symplectique, la 
seconde sur la méthode du crible. On peut se demander s'il est possible de combiner ces deux méthodes pour obtenir un meilleur résultat. C'est probable, mais je n'ai pas réussi à le faire. 

\subsection{Critère d'isogénie de courbes elliptiques et conjecture d'uniformité de Serre}

\begin{lemme} Soit $G$ un groupe fini, $\chi_1,\chi_2: G \rightarrow \Cb$ deux caractères  de $G$ de dimension $d$, non nécessairement irréductibles. 
Soit $$A=\{g \in G, \chi_1(g) \neq \chi_2(g)\}.$$ Alors ou bien $A = \emptyset$, ou bien $$|A| > \frac{|G|}{ 2d^2}.$$
\end{lemme}
\begin{pf} On a $A=\emptyset$ si et seulement si $\chi_1=\chi_2$. Supposons donc $\chi_1 \neq \chi_2$.
Pau Cauchy-Schwarz, on a $| \langle \chi_1,\chi_2 \rangle |^2 \leq \langle \chi_1,\chi_1 \rangle \langle \chi_2,\chi_2 \rangle$,
et cette inégalité est stricte car l'égalité impliquerait que $\chi_1$ et $\chi_2$ sont proportionnels, donc égaux vu que $\chi_1(1)=d=\chi_2(1)$.
On en déduit que l'un au moins des facteurs du membre de droite, disons le premier par symétrie, est strictement plus grand que $| \langle \chi_1,\chi_2 \rangle |$
On a donc $\langle \chi_1,\chi_1 \rangle > \langle \chi_1, \chi_2 \rangle$ et comme ces deux nombres sont des entiers,
$\langle \chi_1,\chi_1 \rangle - \langle \chi_1, \chi_2 \rangle \geq 1$, soit \begin{eqnarray*}
1 &\leq& \frac{1}{|G|} \sum_{g \in G} \chi_1(g) \overline{(\chi_1(g) - \chi_2(g))}   \\
&=& \frac{1}{|G|} \sum_{g \in A} \chi_1(g) \overline{( \chi_1(g)-\chi_2(g))} \\
& \leq & \frac{1}{|G|} |A| d (d+d), \end{eqnarray*}
d'où le lemme.
\end{pf}

\begin{prop} 
Soit $G$ un sous-groupe de $\Gl_2(\F_\ell) \times \Gl_2(\F_\ell)$. On note $g \mapsto g_1$ et $g \mapsto g_2$ les deux projections de $G$ sur $\Gl_2(\F_\ell)$ et on suppose que pour tout $g \in G$, $\det g_1 = \det g_2$. Soit $A = \{g \in G, \tr g_1 \neq \tr g_2\}$. On suppose que $A$ est non-vide. Alors:
\begin{eqnarray} \label{varphiGA} \varphi_G(A) < \ccl{cpropisogenie} \ell^{5/2},\end{eqnarray}
où $\ccr{cpropisogenie}$ est une constante absolue.

Si de plus l'on suppose que $G$ a un sous-groupe $G'$ d'indice $2$ tel que $g_1=g_2$ pour tout $g \in G'$, alors
$$\varphi_G(A) <8 \sqrt{2} \ell^{2}.$$
\end{prop}

\begin{pf} Par le lemme précédent, on a $|A| \geq |G|/8$. On a
\begin{eqnarray} \label{varphiGAeq} \varphi_G(A) \leq \lambda_G(A) |G|/|A| \leq 8 \lambda_G(A).\end{eqnarray}

Soit $A'$ le complémentaire de $A$ dans $G$, i.e. l'ensemble des $(g_1,g_2)$ tels que $\tr g_1 = \tr g_2$.
Si $(g_1,g_2) \in A'$, il y a au plus $\ell^4$ possibilités pour $g_1$, et $g_2$ a le même polynôme caractéristique que $g_1$,
ce qui laisse au plus $\ell^2+\ell$ possibilités pour $g_2$ une fois que $g_1$ est fixé, d'où $|A'| \leq \ell^4(\ell^2+\ell)$. On a par ailleurs 
$UA'=A'$ où $U$ le sous-groupe normal de $G$ des éléments $(x \Id,x\Id)$ de $G$ avec $x \in \F_\ell^\ast$. On a donc 
$\lambda_G(A) \leq 1+\lambda_G(A') = 1+ \lambda_{G/U}(A'/U) \leq 1+ \sqrt{ \frac{\ell^4(\ell^2+\ell)}{\ell-1}}$, et le résultat suit.

Si de plus $G$ contient un sous-groupe $G'$ d'indice $2$ avec $g_1=g_2$ pour tout $g \in G'$, on voit aisément que $|A| \leq |G| \leq 2 \ell^4$ donc
$\lambda_G(A) \leq \sqrt{2} \ell^2$ et $\phi_G(A) \leq 8 \sqrt{2} \ell^2$.
\end{pf}

\begin{theoreme} \label{isogenieeffective} Soit $E$ et $E'$ deux courbes elliptiques sur $\Q$ non isogènes. Soit $M$ le produit des nombres premiers $p$ tels que $E$ ou $E'$ a mauvaise réduction modulo $p$, et posons  $a_p = |E(\F_p)|-1-p$, $a'_p=|E'(\F_p)|-1-p$ pour tout $p$ ne divisant pas $M$. Alors, sous (GRH) et Arin, il existe un $p$ ne divisant pas $M$ avec
$$ p  < \ccl{cisog} (\log M)^2 (\log \log 2 M)^5$$
tel que $$a_p \neq a'_p.$$
Si de plus $E'$  est obtenue de $E$ par torsion par un caractère quadratique, alors il existe un $p$ ne divisant pas $M$ avec
$$p < \ccl{cisog2} (\log M)^2 (\log \log 2 M)^3$$
 tel que $a_p \neq a'_p$. Les constantes $\ccr{cisog}$ et $\ccr{cisog2}$ sont absolues
\end{theoreme}
\begin{remarque} Ce résultat est légèrement meilleur que celui de \cite{serre} qui affirme qu'on peut prendre $p < \cc (\log M)^2 (\log \log 2M)^{12}$ dans le cas général, et $p <\cc (\log M)^2 (\log \log 2M)^{6}$ dans le cas où  $E'$  est obtenue de $E$ par torsion par un caractère quadratique.
\end{remarque}

\begin{pf} 
(On suit  la méthode de Serre, \cite[pages 192--194]{serre}, à ceci près qu'on utilise la majoration de $\phi_G(A)$ donnée par la proposition précédente,
là ou ce que fait Serre revient à utiliser la majoration $\phi_G(A) \leq |G|$.) Toutes les constantes apparraissant dans le preuve qui suit sont absolues.

Soit $p$ le plus petit nombre premier ne divisant pas $M$ tel que $a_p \neq a'_p$. Soit $\ell$ le plus petit nombre premier ne divisant pas $a_p-a'_p$. D'après le lemme 12 de \cite{serre}, page 192, on sait que $\ell < \cc \log |a_p-a'_p| $ 
donc \begin{eqnarray} \label{ellp} \ell <\ccl{ctruc1} \log p\end{eqnarray} en utilisant les estimées de Hasse $|a_p-a'_p| \leq 4 \sqrt{p}$. Considérons alors la représentation $\rho_E \times \rho_{E'} : G_\Q \rightarrow \Gl_2(\F_\ell) \times \Gl_2(\F_\ell)$, o\`u $\rho_E$ (resp. $\rho_{E'}$) est la représentation sur les points de $\ell$-torsion de $E(\bar \Q)$
(resp. $E'(\bar \Q)$). Soit $G$ l'image de $\rho_E \times \rho_{E'}$ et $A \subset G$ l'ensemble des $g=(g_1,g_2)$ tels que $\tr g_1 \neq \tr g_2$.
Notons que $\tr \rho_E(\Frob_p) = a_p \neq a'_p = \tr \rho_{E'} (\Frob_p)$, si bien que $A$ est non vide. Comme de plus $\det \rho_E = \det \rho_{E'}$, les hypothèses de la proposition précédente sont satisfaites et on a donc $$\phi_G(A) < \ccr{cisog} \ell^{5/2}.$$ 
D'après le théorème~\ref{thmpluspetit}, on a $p < \cc \ell^5 (\log |G| + \log M\ell)^2$.  De plus, $|G| < \ell^8$, donc
\begin{eqnarray} p < \cc \ell^5 (\log \ell + \log M)^2 \end{eqnarray}
Comme $\ell < \ccr{ctruc1} \log p$, on obtient
\begin{eqnarray} p < \cc (\log p)^5 (\log \log p + \log M)^2 \end{eqnarray}
On en déduit facilement (cf. \cite{serre}, bas de la page 194 et deux premières lignes de la page 195 -- noter que l'inégalité $M \geq 2$, qui est vraie car "il n'y a pas de courbes elliptiques sur $\Z$", est utilisée implicitement
{\it loc. cit.} pour déduire ($267_R$) de ($266_R$)) 
que $p < \cc (\log M)^2 (\log \log 2 M)^5$.

Dans le cas où $E'$ est obtenue de $E$ par torsion par un caractère quadratique, on a $ \phi_G(A) < \cc \ell^3$ par la proposition, d'où 
$p < \cc \ell^3 (\log \ell + \log M)^2$, puis $p < \cc (\log p)^3 (\log \log p + \log M)^2$, et finalement $p < \cc (\log M)^2 (\log \log 2 M)^3$.
\end{pf}

Concernant la conjecture d'uniformité de Serre (qui est l'assertion qu'il existe une constante absolue $\ccl{cus}$ telle que pour toute courbe elliptique $E$ sur $\Q$ sans multiplication complexe, et pour tout $\ell > \ccr{cus}$, la représentation $\rho_E: \GQ \rightarrow \Gl_2(\F_\ell)$ est surjective), on obtient: 
\begin{cor} Soit $E$ une courbe elliptique sur $\Q$ sans multiplication complexe ayant bonne réduction en dehors de $M$. 
Sous (GRH) et Artin, la représentation $\rho_E : \GQ \rightarrow \Gl_2(\F_\ell)$,
est surjective pour $\ell > \cc \log M (\log \log 2 M)^{3/2}$.
\end{cor}
\begin{pf}
Cela se déduit du cas "torsion quadratique" du théorème précédent exactement comme le théorème 22 de \cite{serre} se déduit du théorème 21'.
\end{pf}
Ceci améliore très légèrement \cite[théorème 22]{serre} en remplaçant l'exposant $3$ du facteur $\log \log 2M$ par $3/2$.

Voici une variante du théorème~\ref{isogenieeffective}:
\begin{theoreme} Soit $N\geq 1$ et $k \leq 2$ deux entiers, $f=\sum_{n \geq 1} a_n q^n$ et $f'=\sum_{n \geq 1} a'_n q^n$ 
deux formes modulaires cuspidales nouvelles de poids $k$ et niveau $\Gamma_0(N)$. Notons $M$
le produit des nombres premiers divisant $N$.  Si $f \neq f'$, il existe un $p \nmid M$,
$$ p  < \ccl{cmodform} k^5  (\log kM)^2 (\log k + \log \log 2kM)^5 $$ tel que $a_p \neq a'_p$.
\end{theoreme}
\begin{pf} La preuve est la même que celle du théorème~\ref{isogenieeffective}, à ceci près qu'on remplace les estimées de Hasse par les estimées de Hecke 
$|a_p-a'_p| \leq \cc p^{k/2}$, d'où (\ref{ellp}) par $\ell \leq \cc k \log p$, puis la représentation $\rho_E$ (resp. $\rho_{E'}$)
par $\rho_f$ (resp. $\rho_{f'}$) défine comme les réductions modulo $\ell$ d'un réseau stable quelconque pour la représentation $\ell$-adique attachée à $f$ (resp. $f'$).  On obtient donc \begin{eqnarray} p < \cc k^5 (\log p)^5 (\log \log p + \log kM)^2. \end{eqnarray}
En divisant cette inégalité par $ (\log p)^7$, on obtient 
$p/(\log p)^7 < \cc k^5 (1+ \log kM)^2 < \cc k^5(\log 2kM)^2$,
d'où $p^{1/2} < \cc k^5 (\log 2kM)^2$ et donc $\log p < \cc + 10 \log k + 2 \log \log 2kM < \cc (\log k + \log \log 2kM)$ puisque $\log k>0$
et $\log \log 2kM > 0$ (car $2kM \geq 4$). On obtient donc
$$p < \cc k^5  (\log k + \log \log 2kM)^5 (\log \log (\log k + \log \log 2kM)+\log kM)^2,$$
d'où le résultat.
\end{pf}
\begin{remarque} La borne sur $p$ obtenue est polynomiale en $k$ (de degré $5+\epsilon$ -- en utilisant la version de Lagrias-Odlyzko de Chebotarev, on obtiendrait un degré $12+\epsilon$), logarithmique en $M$. Ceci est à comparer avec la borne classique due  à Sturm $p < k [\Sl_2(\Z):\Gamma_0(N)] / 12$ (cf. \cite{sturm}). Cette borne est linéaire en $k$ (et donc meilleure que la nôtre sous cet aspect), 
mais polynomiale en $N$ (en quoi elle est moins bonne que la nôtre, logarithmique en $M \leq N$). La borne de Sturm a aussi les avantages d'être plus ancienne, élémentaire, inconditionelle, et de ne pas faire intervenir de constante non calculée.
\end{remarque}

\end{document}